\documentclass[a4paper,10pt]{amsart}
\usepackage{amsmath,amssymb,amsfonts}
\usepackage{verbatim}
\usepackage[utf8]{inputenc}
\usepackage{graphicx}
\usepackage{caption}
\usepackage{pstricks} 

\allowdisplaybreaks[1]
\numberwithin{equation}{section}

\newcommand{\<}{\langle}
\renewcommand{\>}{\rangle}

\newcommand{\C} {\mathbb{C}}

\newcommand{\tC}{\tilde C}
\newcommand{\cC}{\check C}
\newcommand{\tP}{\tilde P}
\newcommand{\cP}{\check P}

\renewcommand{\d}[1][x]{\,\operatorname{d}\!#1}
\newcommand{\dds}{\fract{\text{d}}{\text{ds}}}
\newcommand{\ddt}{\fract{\text{d}}{\text{dt}}}
\newcommand{\diag}{\operatorname{diag}}
\renewcommand{\div}{\operatorname{div}}

\newcommand{\eps}{\varepsilon}

\newcommand{\F}{\mathcal{F}}
\newcommand{\fract}[2] { {\textstyle\frac{#1}{#2}} }
\newcommand{\fu}{\frac f{f_\infty}}

\renewcommand{\H}{\mathcal{H}}

\newcommand{\Id}{\operatorname{Id}}
\newcommand{\im}{\operatorname{im}}
\newcommand{\ild}{\int\limits_{\R^d}}

\newcommand{\N} {\mathbb{N}}
\newcommand{\norm}[1]{\lVert#1\rVert}

\newcommand{\R}{\mathbb{R}}

\newcommand{\tb}{\textbf}
\newcommand{\tf}{\tilde f}
\newcommand{\tr}{\mathrm{Tr}}

\newcommand{\tL}{\tilde L}
\newcommand{\tQ}{\tilde{\mathcal Q}}

\newcommand{\M}{\mathcal M}
\renewcommand{\P}{\mathcal P}
\renewcommand{\phi}{\varphi}

\newtheorem{Proposition}{Proposition}[section]
\newtheorem{Corollary}[Proposition]{Corollary}
\newtheorem{Lemma}[Proposition]{Lemma}
\newtheorem{Theorem}[Proposition]{Theorem}
\newtheorem{Definition}[Proposition]{Definition}
\newtheorem{Remark}[Proposition]{Remark}
\newenvironment{Proof}[1][]{\textbf{Proof#1: }}{\hfill $\Box$\\}

\title{Sharp entropy decay for hypocoercive and non-symmetric Fokker-Planck equations with linear drift}
\author{Anton Arnold}
\address{Anton Arnold, Institute for Analysis and Scientific Computing, TU Vienna, Wiedner Hauptstraße 8-10, 1040 Vienna, Austria. e-mail: anton.arnold@tuwien.ac.at}
\author{Jan Erb}
\address{Jan Erb, Institute for Analysis and Scientific Computing, TU Vienna, Wiedner Hauptstraße 8-10, 1040 Vienna, Austria. e-mail: jan.erb@tuwien.ac.at}
\date{\today}


\begin{document}

\begin{abstract}
\noindent
We investigate the existence of steady states and exponential decay for hypocoercive Fokker--Planck equations on the whole space with drift terms that are linear in the position variable. For this class of equations, we first establish that hypoellipticity of its generator and confinement of the system is equivalent to the existence of a unique normalised steady state. These two conditions also imply hypocoercivity, i.e.\ exponential convergence of the solution to equilibrium. 

Since the standard entropy method does not apply to degenerate parabolic equations, we develop a new  modified entropy method (based on a modified, non-degenerate entropy dissipation--like functional) to prove this exponential decay in relative entropy (logarithmic till quadratic) -- with a sharp rate. Furthermore, we compute the spectrum and eigenspaces of the generator as well as flow-invariant manifolds of Gaussian functions.

Next, we extend our method to kinetic Fokker--Planck equations with a class of non-quadratic potentials. And, finally, we apply this new method to non-symmetric, uniformly parabolic Fokker-Planck equations with linear drift. At least in 2D this always yields the sharp exponential envelopes for the entropy function. In this case, we obtain even a sharp multiplicative constant in the decay estimate for the non-symmetric semigroup.  
\end{abstract}

\maketitle

\tableofcontents
 
\section{Introduction}
\label{sec:intro}

This paper is concerned with the large-time behaviour of degenerate parabolic Fokker-Planck equations. In applications, the most important model of this class is the kinetic Fokker-Planck equation
\begin{align}
 \label{kinFP}  
 \partial_t f +v\cdot\nabla_x f-\nabla_x V\cdot\nabla_v f&= \nu \div_v (vf)+\sigma\Delta_v f\,;\quad x,\,v\in\R^n;\,t>0,
\end{align}
describing the time evolution of the phase space probability density $f(t,x,v)$, e.g.\ in a plasma  \cite{RiFP89, Vi02}.
Here, $V=V(x)$ is a given confinement potential for the system, and $\nu,\,\sigma$ denote the (positive) friction and diffusion parameters, respectively. For quadratic potentials, \eqref{kinFP} has a linear drift term and its solution can be represented by the Mehler formula \cite{Ho95}. But for non-quadratic potentials its large-time behaviour (i.e.\ exponential convergence towards the steady state) is highly non-trivial.

In the main part of this paper we shall analyse Fokker-Planck equations with linear drift terms. Our objective is to develop a new entropy method for proving the exponential decay of Fokker-Planck solutions towards equilibrium and to understand the structure of their entropy decay --- beyond explicit representation formulas. This new method has the potential to be generalised to non-quadratic operators. To this end we shall also illustrate that it can be extended to certain kinetic Fokker-Planck equations with non-quadratic potentials.\\

We start to consider a Fokker-Planck equation on $(0,\infty)\times\R^d$ of the form
\begin{align}
\label{masterequ} \partial_t f = Lf &:= \div (D\nabla f + F f), \\
\nonumber	f(t=0)&= f_0 \in L^1(\R^d),\\
\nonumber \ild f_0\d=1&, f_0\geq0.
\end{align}
Throughout this paper, we make the assumptions
\begin{itemize}
 \item $D^T=D\in\R^{d\times d}$ is positive semidefinite and constant in $x$,
 \item $F:\R^d\to\R^d$, $x\mapsto Cx$ with $C\in\R^{d\times d}$.
\end{itemize}
So we consider the degenerate parabolic Fokker-Planck equation
\begin{align}
 \label{linmasterequ}  \partial_t f &= Lf := \div (D\nabla f + Cxf)=\div(D\nabla f)+x^TC^T\nabla f + \tr(C)f,
\end{align}
and analyse solutions that satisfy $f(t,\cdot)\in L^1(\R^d)$ along with $\int\limits_{\R^d} f(t,x) \d =1$ for all $t>0$.\\
Since $D$ is symmetric, it can be diagonalised and normalised (all entries 0 or 1) by rescaling the space variable. We can thus always assume $D$ to be a ``defect" identity, i.e.
\begin{align*}
 D=\diag\{\underbrace{1,\dots,1}_{k},\underbrace{0,\dots,0}_{d-k}\},
\end{align*}
where $k:=\operatorname{rank} D$, $1\leq k<d$. The case $k=d$ has been studied extensively, see for example \cite{ArMaToUn01}. In the case $k<d$, the operator $L$ is not elliptic, and classical parabolic results will not apply for (\ref{linmasterequ}). However, if certain conditions on $C$ are met, the solution still retains typical parabolic properties: regularisation, long-term convergence and a maximum principle. This behaviour stems from an interaction between the degenerate dissipative part and the non-symmetric part of $L$. Regularity and maximum principle for solutions to (\ref{linmasterequ}) are due to the hypoellipticity of the operator $\partial_t-L$; the long-term decay of solutions to (\ref{linmasterequ}) is generally connected to \emph{hypocoercivity}.\\

A very good, broad discussion of hypocoercivity can be found in \cite{ViH06}, which also contains a precise definition of hypocoercivity:
\begin{Definition}
\label{hypodefinition}
 Let $H$ be a Hilbert space, $L$ an unbounded operator on $H$ with kernel $\mathcal{K}$. Let $\tilde H$ be another Hilbert space, which is continuously and densely embedded in $\mathcal{K}^\perp$. Then $-L$ is said to be \emph{hypocoercive} on $\tilde H$ if and only if there is $\lambda>0$ and some constant $c\ge1$ such that
 \begin{align*}
  \forall h\in \tilde H,\,\forall t\ge0: \norm{e^{tL}h}_{\tilde H} \leq ce^{-\lambda t}\norm{h}_{\tilde H}.\\
 \end{align*}
\end{Definition}

\cite{ViH06} also establishes a general criterion for exponential convergence of solutions for a class of hypocoercive evolution equations, based on a Lyapunov functional equivalent to a weighted $H^1$-norm. While the main theorem in \cite{ViH06} covers a wide class of problems, the price paid is in the estimate for the decay rate, which is off by orders of magnitude.\\

In the last few years several papers dealt with the large-time behaviour of hypocoercive equations. But to our knowledge, sharp decay rates (so far only in $L^2$) were obtained only via a spectral analysis: in \cite{OPPS12} for parabolic equations associated to hypoelliptic quadratic operators; and in \cite{GaMiS13} for two specific toy models, using the spectral decomposition of their generators. 
In \cite{MoNeQ06} several collisional kinetic models (including the Fokker-Planck equation, linearised Boltzmann and Landau) are analysed on the torus (in the spatial variable): exponential convergence to the steady state is shown in the $H^1$--norm. In \cite{MoH14}, a decay estimate is obtained for a 2-dimensional kinetic Fokker-Planck model using higher order time derivatives of the $L^2$-norm of solutions and their space derivative.
Also \cite{DuH09} and \cite{BaB13} study dissipative kinetic models (i.e.\ with $k = \frac d2$) in $H^1$. While \cite{DuH09} uses a \emph{macro-micro decomposition} of the models, \cite{BaB13} is based on an (augmented) $\Gamma_2$--calculus and local computations (in contrast to the integrated functionals used by most other authors), cf.\ also \cite{BaEmD85}. 
\cite{DoMoScH10} and \cite{BaB13} also analyse much more general hypocoercive equations. Along with \cite{DuH09} they require the following restriction on the interaction between the degenerate dissipative part and the non-symmetric part of $L$: It is assumed that the matrix $C^T$ does not map any subspace of the kernel of $D$ into the kernel of $D$, which is equivalent to using only first order H\"ormander-commutators to span all of $\R^d$ (i.e.\ $\tau=1$ in Lemma \ref{Definiteness} (iii) and Remark \ref{tauremark} below; cf.\ also \S3 in \cite{BaB13}). But this condition is more restrictive than necessary. In this paper we shall impose a weaker condition 
(see the first part of \emph{condition (A)} in Definition \ref{condAdef} below; or \cite{ViH06}). \\

The common approach to study the long-term behaviour of hypocoercive equations has been via a Lyapunov functional - usually on a weighted $H^1$-space, but \cite{ViH06} also contains (in Theorem 28) a Lyapunov functional based on the logarithmic entropy. In \cite{DoMoScH10}, the authors get rid of the $H^1$-regularity restriction on initial states and prove decay towards the steady state using a modified $L^2$-norm. In \cite{ViH06}, it is shown that even for methods based upon $H^1$-functionals, one can often get rid of the regularity assumptions by using the regularisation of the semigroup $e^{tL}$. So far, there is no knowledge on the decay of general entropies ``between'' logarithmic and quadratic, nor on sharp decay rates for equations of type (\ref{linmasterequ}). In this paper we shall modify the entropy method (see \cite{ArMaToUn01}, \cite{BaEmD85}-\cite{BaEmI85}) to achieve all three results for equations of type (\ref{linmasterequ}): no $H^1$-regularity requirement for the initial state, sharp decay 
rates, and decay for a wide class of relative entropies.

\begin{figure}[htbp]
\begin{center}
\includegraphics[width=7.8cm]{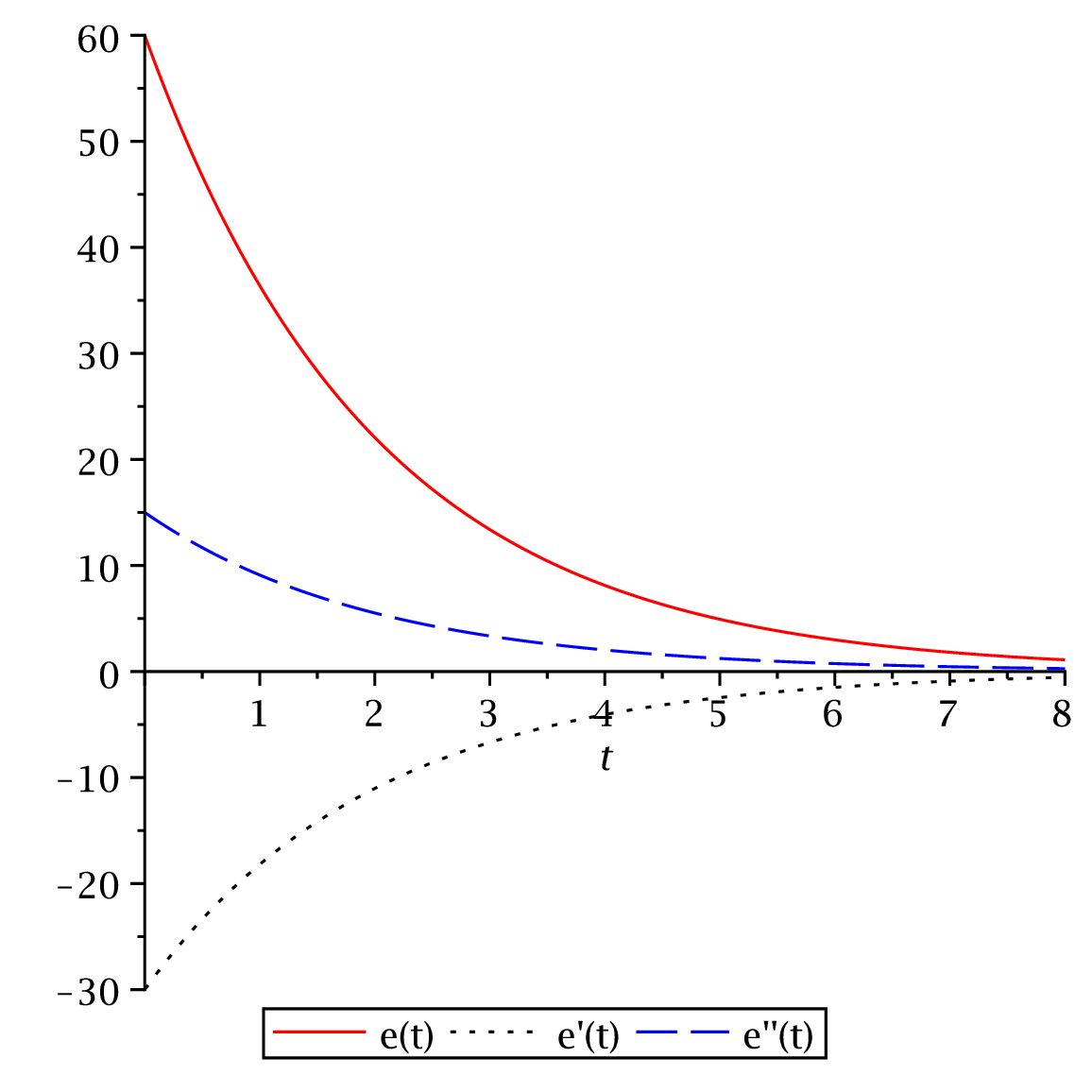}\hfill
\includegraphics[width=7.8cm]{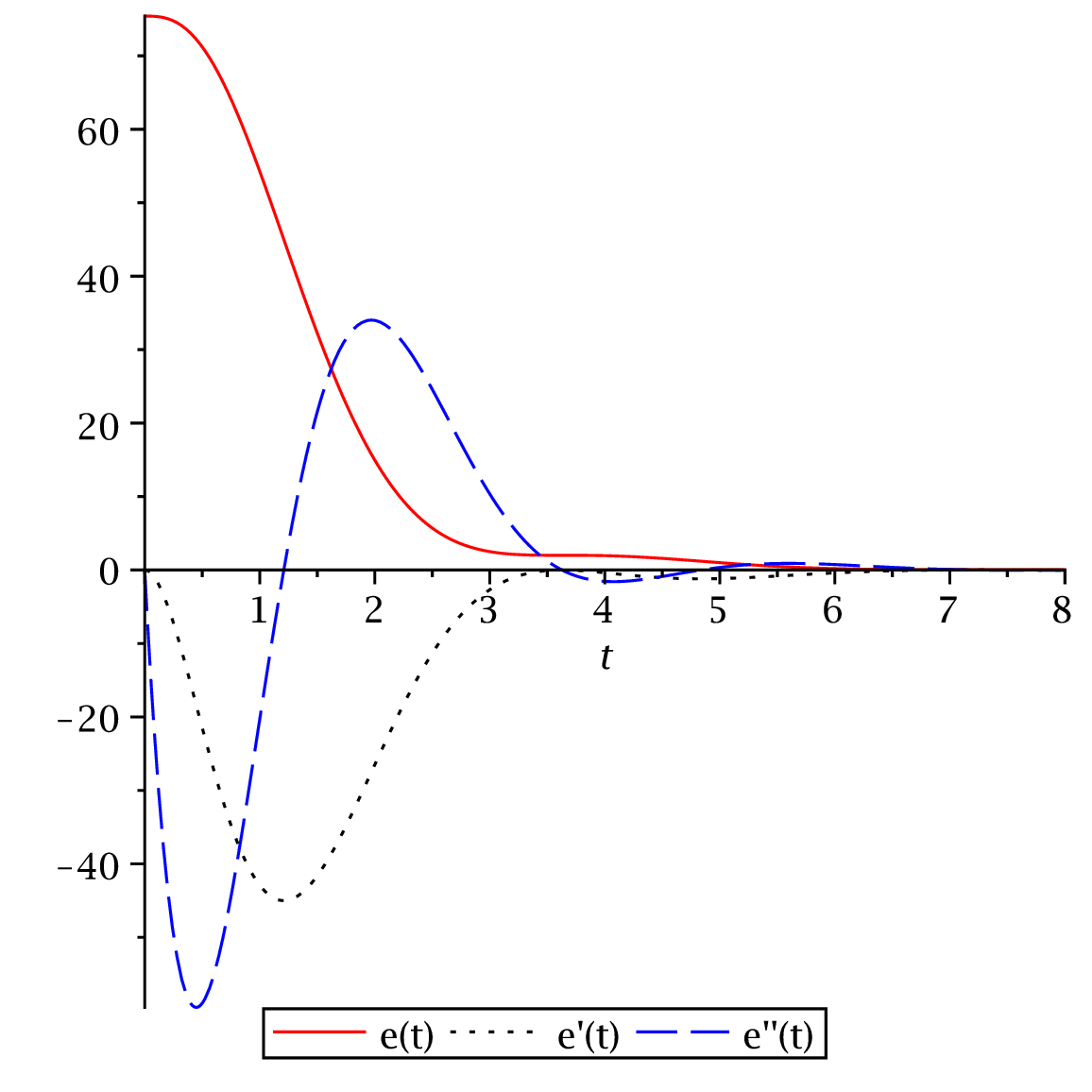}
\end{center}
\vspace{-0.3cm}
\caption{\label{fig:classical} {\footnotesize Prototypical behaviour of the logarithmic relative entropy $e(t)$, its first and second time derivatives. (a) Left: Non-degenerate case: The inequalities $e'\leq-\mu e$, $e''\geq -\mu e'$ can be obtained. (b) Right: Degenerate case; equation \eqref{linmasterequ} with 
 $D=\diag(1,\,0)$, 
 $C=[1\;-1\,;\,1\;0]$ 
 \;\;: The inequalities $e'\leq-\mu e$, $e''\geq -\mu e'$ are wrong, in general.} 
}
\end{figure}

The strategy of the standard entropy method is to derive first a differential inequality between the first and second time derivative of the relative entropy (of the solution w.r.t.\ the equilibrium state). Their time evolution in a prototypic situation is shown in Fig.\ \ref{fig:classical}(a). Integration in time of the inequality then allows to deduce exponential decay of the relative entropy, which is a convex function of time. But this approach is not feasible for degenerate Fokker-Planck equations, since the entropy dissipation can vanish for states other than the equilibrium. Hence, the second time derivative of the entropy may change its sign along a trajectory. So the entropy functional exhibits a ``wavy'' decay in time, see Fig.\ \ref{fig:classical}(b) and Fig.\ \ref{fig:non-symmetric}. This oscillatory behaviour is also known from space-inhomogeneous kinetic equations (cf.\ \S3.7 of \cite{Vi02}; and \cite{FMP06} for a numerical study on the Boltzmann equation).

As a remedy for the analysis, one therefore has to use either some ``modified relative entropies'' (as in \cite{DoMoScH10}) or ``modified entropy dissipations''. Here, we shall introduce an auxiliary functional -- structurally related to the entropy dissipation, but an upper bound for the latter. A Bakry-\'Emery-type estimate then yields exponential decay of this auxiliary functional, and consequently also of the entropy dissipation. A convex Sobolev inequality with the auxiliary functional as its relative Fisher information \cite{ArMaToUn01} finally 
yields the exponential decay of the relative entropy. Initially, this approach shall need an additional regularity assumption for the initial state. But
this can then be removed using the regularisation of the parabolic equation (\ref{linmasterequ}), as in \cite{ViH06}.\\

The novelties of this paper include:
\begin{enumerate}
\item A new modified entropy method for hypocoercive and non-symmetric Fokker-Planck equations with the potential of a generalization to (some) equations with nonlinear drift;
\item sharp exponential decay rates for relative entropies ``between'' logarithmic and quadratic functionals;
\item clarification of \emph{global entropy decay} estimates as envelopes for entropy functionals that are non-convex in time.
\end{enumerate}

This paper is organised as follows: In Section \ref{sec:solpos}, we give a sufficient and necessary condition on the matrix $C$ such that (\ref{linmasterequ}) is hypocoercive. We establish that the solutions will be positive for any $t>0$. Section \ref{sec:steadystate} follows this up by explicitly giving the unique (up to normalisation) steady state $f_\infty$ and discussing the operator $L$ in $L^2(\R^d,f_\infty^{-1})$, the standard space for Fokker-Planck equations. In Section \ref{sec:modentmod} we state our main result in the Theorems \ref{entropydecay}, \ref{convergencerate}: a modified entropy method allows to compute an explicit decay rate for solutions of (\ref{linmasterequ}) in relative entropy. In Section \ref{sec:spectrum} we compute the spectrum of $L$ on the weighted space $L^2(\R^d,f_\infty^{-1})$ as well as flow-invariant manifolds (the eigenspaces of $L$ and Gaussian manifolds). 
Sharpness of the decay rate and the multiplicative constant will be shown in Theorem \ref{sharpdecayrate} and Proposition \ref{sharp-constant} of Section \ref{sec:sharprate}.
In Section \ref{sec:kinFP} we illustrate the 
extension of our new method to kinetic 
Fokker-Planck equations with nonlinear drift terms.
Finally, in Section \ref{sec:nonsymmFP} we show how the presented method improves known entropy decay rates also for non-symmetric Fokker-Planck equations that are non-degenerate. In this context we shall distinguish between \emph{sharp local} and \emph{sharp global decay rates}. For the latter, we derive an exponential function that is the global envelope for the entropy functional.
\\


\section{Existence of solutions and positivity}
\label{sec:solpos}

If $D$ is not regular, the operator $L$ is neither coercive nor elliptic. In general, such an operator does not have a unique normalised steady state. We thus need additional assumptions on the parameters in $L$, which shall be assumed throughout \S2-6 of the paper:
\begin{Definition}
\label{condAdef}
 The operator $L$ from (\ref{linmasterequ}) fulfils \tb{\emph{condition (A)}} if and only if
\begin{itemize}
 \item there is no non-trivial $C^T$-invariant subspace of $\ker D$,
 \item the matrix $C\in\R^{d\times d}$ is positively stable\footnote{A matrix is positively stable iff all eigenvalues have real part greater than zero.}.
\end{itemize}
\end{Definition}

The first part of condition (A) is equivalent to the hypoellipticity of $\partial_t-L$ (cf.\ \S1 of \cite{Ho67}), and it allows for smooth solutions to \eqref{linmasterequ} (see Proposition \ref{solexistence} below). Due to the special form of $D$, $C$ cannot be diagonal (under condition (A)) unless $k=d$.\\ 
The second condition, positive stability of $C$, means that there is a confinement potential. While there are solutions even without a confinement potential, there would be no steady state. Indeed, Theorems \ref{ssexistence} and \ref{convergencerate} will show that condition (A) is both sufficient and necessary for the existence of a unique normalised steady state and exponential convergence of solutions to the steady state. So for equations of type \eqref{linmasterequ}, hypoellipticity and confinement are equivalent to hypocoercivity.\\

\begin{Proposition}
\label{solexistence}
 Let $f_0\in L^1(\R^d)$. Then there is a unique solution $f\in C^\infty(\R^+\times\R^d)$ of (\ref{linmasterequ}) iff no non-trivial subspace of $\ker D$ is invariant under $C^T$.
\end{Proposition}
\begin{Proof} 
 See page 148 of \cite{Ho67}.
\end{Proof}\\
If the hypoellipticity condition in Proposition \ref{solexistence} does not hold, (\ref{linmasterequ}) clearly also has a unique solution, but it would be less regular. 

A heuristic explanation of this condition is that the solution cannot stay in the kernel of the dissipative part of $L$, and therefore the evolution under (\ref{linmasterequ}) acts dissipative in all space directions: If one considers merely the drift part of the equation,
\begin{align}
\label{driftequation}
 f_t &= (Cx)\cdot \nabla f,
\end{align}
the solution is $f(t,x)=f_0(e^{Ct}x)$. So, for the dissipative part to ``extend" to the whole space, one needs that $e^{Ct}x$ reaches the whole space $\R^d$ for all $x\in\im D$ ($\im D$ being the image of $D$). Conversely, this means that $e^{C^Tt}x$ evolves into $\im D$ for all $x\in\ker D$ as shown in Lemma \ref{Definiteness} (iv) below.\\ 

In the following lemma we give four equivalent characterisations of the hypoellipticity of $L$.


\begin{Lemma}
\label{Definiteness}
 The following four statements are equivalent:
 \begin{itemize}
  \item[(i)] No non-trivial subspace of $\ker D$ is invariant under $C^T$.
  \item[(ii)] No eigenvector $v$ of $C^T$ fulfils $Dv=0$.
  \item[(iii)] There exist constants $\tau\in\{1,\dots,d-k\}$ and $\kappa>0$ such that
\begin{align}
\label{sumTdefinite} \sum\limits_{j=0}^{\tau} C^jD(C^T)^j \geq \kappa\Id,
\end{align}
where $k=\operatorname{rank} D$.
  \item[(iv)] For any $t\in\R$, $h>0$, it holds that
\begin{align*}
 \forall 0\neq\xi\in\ker D\ \exists s\in[t,t+h]\ \exists \eta\in\im D:\ & \<e^{C^Ts}\xi,\eta\> = 1.
\end{align*}
 \end{itemize}
\end{Lemma}
\begin{Proof}
(i)$\Rightarrow$(ii): 
Each eigenvector (or pair of complex conjugated eigenvectors) of $C^T$ with $Dv=0$ spans a $C^T$--invariant subspace of $\ker D$. \\
(ii)$\Rightarrow$(i): Any $C^T$--invariant subspace of $\ker D$ contains a (possibly complex) eigenvector of $C^T$.\\
(i)$\Rightarrow$(iii): All matrices $C^jD(C^T)^j$ are symmetric and positive semidefinite, since $D$ is symmetric and positive semidefinite. It suffices to show that for any vector $v\neq0$, there exists $j\leq d-k$ with $D^{\frac12}(C^T)^jv=D(C^T)^jv\neq 0$, since then $\sum\limits_{j=0}^{\tau} C^jD(C^T)^j$ is regular for $\tau:=\max\limits_{v\neq 0}\min\limits_{j\in\N} \{j| D^{\frac12}(C^T)^jv\neq0\}$.\\
If $v\notin\ker D$, we choose $j=0$, and hence $Dv\neq0$. So let now $0\neq v\in \ker D$. Then either $C^Tv\notin\ker D$, in which case $DC^Tv\neq 0$, or $C^Tv\in\ker D$. Repeating this procedure, we see that either there is $j\leq d-k$ such that $(C^T)^jv\notin\ker D$ or $\forall 0\leq j\leq d-k: (C^T)^jv\in\ker D$. Assume the latter. Since the dimension of $\ker D$ is $d-k$, the $d-k+1$ vectors $(C^T)^jv$, $0\leq j\leq d-k$ are not linearly independent. Thus, 
$\exists\,l\in\{1,...,d-k\}$ such that
$\operatorname{span}\{C^Tv,\dots,(C^T)^{l}v\}=\operatorname{span}\{v,\dots,(C^T)^{l-1}v\}$. 
Hence, $\operatorname{span}\{v,\dots,(C^T)^{l-1}v\}$ is a $C^T$-invariant subspace of $\ker D$, which has to be trivial due to condition (A). But then $v=0$, which is a contradiction.\\
(iii)$\Rightarrow$(i): If $0\neq v\in\ker D$, then by (iii) there is a $j\in\{1,\dots,\tau\}$ such that $D^\frac12(C^T)^jv\neq0$, i.e. $(C^T)^jv\not\in\ker D$. Thus, no non-trivial subspace of $\ker D$ can be invariant under $C^T$.\\
(i)$\Rightarrow$(iv):
Let $0\neq\xi\in\ker D$, $t\in\R$, and $h>0$. To proceed by contradiction we assume
\begin{align}
\label{kernelpart}
 \forall s\in[t,t+h]\ \forall \eta\in\im D:\<e^{C^Ts}\xi,\eta\> &= 0.
\end{align}
This implies
\begin{align*}
 \forall s\in[t,t+h]\ :e^{C^Ts}\xi\in\ker D,
\end{align*}
and therefore in particular $\nu:=e^{C^Tt}\xi\in\ker D$.
Differentiating (\ref{kernelpart}) with respect to $s$ yields
\begin{align}
\label{equationALPHA}
 \forall s\in[t,t+h]\ \forall \eta\in\im D:\<e^{C^Ts}C^T\xi,\eta\> &= 0.
\end{align}
But this implies $C^T\nu\in\ker D$. Differentiating (\ref{equationALPHA}) repeatedly with respect to $s$ yields $(C^T)^j\nu\in\ker D$ for any $0\leq j \leq d-1$. Hence, $\operatorname{span} \{\nu,\dots,(C^T)^{d-1}\nu\} \subset \ker D$ is a $C^T$-invariant subspace of $\ker D$. This contradicts (i).\\
(iv)$\Rightarrow$(i):
Let $\xi\ne0$ be in a $C^T$--invariant subspace of $\ker D$, i.e.\ $(C^T)^j\xi\in\ker D$ for all $j\in\N_0$. Since $C^T\in\R^{d\times d}$, $e^{C^T s}$ is a polynomial in $C^T$. Hence, $e^{C^T s}\xi\in\ker D$ $\forall\,s\in[t,t+h]$ which contradicts (iv).
\end{Proof}\\

\begin{Remark}
\normalfont
If $\tau$ is the minimal constant for which \eqref{sumTdefinite} holds, then $L$ fulfils the \emph{finite rank Hörmander condition} of order $\tau$ (see \cite{Ho67}, Theorem 1.1). Using $\tau=d-k$ in (\ref{sumTdefinite}) is the worst-case scenario. But in many examples, $\sum\limits_{j=0}^\tau (C)^jD(C^T)^j$ with $\tau<d-k$ is already positive definite. This is the case in the kinetic equations discussed in  \cite{BaB13} and \cite{DuH09}, which require $\tau=1$ and $k=\fract d2$. Also in \cite{DoMoScH10}, $\tau=1$ is assumed.  
\label{tauremark}
\end{Remark}

We shall now discuss further the connection between restrictions on $\tau$ and the first part of conditon (A). Several approaches from the literature require a stricter condition than in Proposition \ref{solexistence}: ``That no subspace of the kernel of $D$ be mapped into the kernel of $D$ by $C^T$,'' which is equivalent to requiring $\tau=1$. To illustrate this restriction, we consider the examples
\begin{align*}
 D_1=\left(\begin{array}{cccc} 1 & 0 & 0 & 0 \\ 0 & 1 & 0 & 0 \\ 0 & 0 & 0 & 0 \\ 0 & 0 & 0 & 0 \end{array}\right); &\quad C_1^T=\left(\begin{array}{cccc} 1 & 0 & -1 & 0 \\ 0 & 1 & 0 & -1 \\ 1 & 0 & 0 & 0 \\ 0 & 1 & 0 & 0 \end{array}\right)
 \qquad\mbox{which implies }\tau=1,
\end{align*}
and
\begin{align*}
 D_2=\left(\begin{array}{cccc} 1 & 0 & 0 & 0 \\ 0 & 1 & 0 & 0 \\ 0 & 0 & 0 & 0 \\ 0 & 0 & 0 & 0 \end{array}\right); &\quad C_2^T=\left(\begin{array}{cccc} 1 & 0 & 0 & 0 \\ 0 & 1 & -1 & 0 \\ 0 & 1 & 0 & -1 \\ 0 & 0 & 1 & 0 \end{array}\right)
 \qquad\mbox{which implies }\tau=2.
\end{align*}

In both cases, \eqref{linmasterequ} has a unique normalised steady state and all solutions converge exponentially to it. For the case of $D_1$ and $C_1$, the condition given in \cite{DoMoScH10}, \cite{DuH09} and \cite{BaB13} holds - no subspace of $\ker D$ is mapped into $\ker D$ by $C^T$. In the case of $D_2$ and $C_2$, that condition does not hold, but condition (A) holds. 
The difference can be seen as follows: consider a vector of the form $(0,0,0,a)^T$. If we apply $C_1^T$ to this vector, it is moved out of the kernel of $D$. However, if we apply $C_2^T$, it is not. In order to move it out of the kernel of $D$, we need to apply $C_2^T$ twice (i.e. multiply by $(C_2^T)^2$). Invariance of a subspace $U$ under $C^T$ means that $C^TU\subset U$. So the condition given in \cite{ViH06}, \cite{Ho67} and in this paper is less strict. As will be shown in §3, condition (A) is equivalent to the existence of a unique normalised steady state.\\

Let us recall from \S1 of \cite{Ho67} the Green's function for (\ref{linmasterequ}) (see also Lemma 1.5 of \cite{Erb14} for a short proof):
\begin{Lemma}
\label{greenfunction}
 Let the first part of condition (A) hold.
 Then the Green's function $g$ to (\ref{linmasterequ}) is given by 
 \begin{align}\label{green}
  g(t,x) &= \frac1{(2\pi)^{\fract d2}\det(W(t))} \exp(-x^TW(t)^{-1}x),
 \end{align}
 where
 \begin{align*}
  W(t) &= \int\limits_0^t e^{C(s-t)}De^{C^T(s-t)} \d[s]
 \end{align*}
 is positive definite for all $t>0$.
\end{Lemma}

We now state an existence result on solutions in $L^p$, which is similar to Corollary 3.1 from \cite{CaDoMaSpF04}:
\begin{Corollary}
\label{Lpsolutions}
 Let $f_0\in L^1(\R^d)\cap L^p(\R^d)$, for a $p\in[1,\infty]$. Then there exists a unique classical solution $f$ to (\ref{linmasterequ}) with $f\in C([0,\infty),L^p(\R^d))\cap C^\infty(\R^+\times \R^d)$. If $\int\limits_{\R^d} f_0 \d =1$, it follows that $\int\limits_{\R^d} f(t) \d =1$ for all $t>0$.
\end{Corollary}
\begin{Proof}
 Proposition \ref{solexistence} already yields a smooth solution $f$ for any $t>0$. With the Green's function from Lemma \ref{greenfunction} we obtain
 \begin{align*}
  f(t,\cdot) &= g(t,\cdot)*f_0.
 \end{align*}
 Applying Young's inequality yields
 \begin{align*}
  \norm{f(t)}_{L^p(\R^d)} &= \norm{g(t)*f_0}_{L^p(\R^d)} \leq \norm{f_0}_{L^p(\R^d)}\norm{g(t)}_{L^1(\R^d)},
 \end{align*}
 where $\norm{g(t)}_{L^1(\R^d)} = 1$. The claimed mass conservation then follows from the divergence form of the operator $L$.
\end{Proof}\\

Next we pass to the positivity of solutions. For a \emph{non-degenerate} Fokker-Planck equation, the solution for a non-trivial $f_0\ge0$ is globally positive for any positive time. This follows from a strong maximum principle supplied by the fully parabolic operator. In our degenerate case a (standard) strong maximum principle does not hold. However, global positivity still holds and it is important for the computations in the entropy method in \S\ref{sec:modentmod}.
\begin{Theorem}
\label{globalpos}
Let the first part of assumption (A) hold and $f_0\in L^1_+(\R^d)$. 
The solution $f$ to (\ref{linmasterequ}) then satisfies
\begin{align*}
 \forall t>0\ \forall x\in\R^d: f(t,x)>0.
\end{align*}
\end{Theorem}

This theorem follows directly from the strict positivity of the Green's function $g$ from Lemma \ref{greenfunction}. However, we give a second proof via a \emph{sharp maximum principle} for degenerate elliptic-parabolic equations, cf.\ \cite{Hi70}. This second approach is more general and will be used in \S\ref{sec:kinFP} for the kinetic Fokker-Planck equation with a nonlinear drift coefficient. To this end we introduce some notation.\\
First, we rewrite our operator in degenerate elliptic form:
\begin{align*}
 \tilde Lf &:= \left[\left(\begin{array}{c} \partial_t \\ \nabla  \end{array}\right)^T\tilde D\left(\begin{array}{c} \partial_t \\ \nabla \end{array}\right)\right]f + b\cdot\left(\begin{array}{c} \partial_t \\ \nabla \end{array}\right)f,
\end{align*}
where
\begin{align*}
 \tilde D &:= \left(\begin{array}{cc} 0 & 0 \\ 0 & D \end{array}\right)\in \R^{(d+1)\times(d+1)},\\
 b(x) &:= \left(\begin{array}{c} -1 \\ Cx \end{array}\right) \in \R^{d+1}.
\end{align*}
Comparing this with our original operator $L$, we have
\begin{align}
 \label{tildecomparison}
 \tilde Lf &= Lf - f_t - \tr(C)f.
\end{align}
Due to the special form of $D$, the columns $d_j$ of $\tilde D$ are of the form
\begin{align*}
  (d_j)_l&=\delta_{jl}
\end{align*}
for $2\leq j\leq k+1$ ($k=\operatorname{rank} D$) and $d_j=0$ for $j=1$; $k+1<j\leq d+1$. With this notation, we shall now introduce drift and diffusion trajectories:
\begin{Definition}\label{DDtrajectories}
 Let $\Omega$ be a connected open set in $\R^{d+1}$, $p_0\in\Omega$.\\
\begin{itemize}
 \item If $p(s)$ is the solution to
\begin{align*}
 \dds p(s) &= d_j, \\
   p(0) &= p_0,
\end{align*}
 with some $2\leq j\leq k+1$ and $p(s)\in\Omega$ for $s_1\leq s \leq s_2$ with some $s_1<0<s_2$, then we call $\Gamma:=\{p(s)\,|\,s_1\leq s \leq s_2\}$ a \emph{diffusion trajectory} running through $p_0$.\\
 \item If $p(s)$ is the solution to
\begin{align}\label{p-ODE}
 \dds p(s) &= b(p(s)), \\
   p(0) &= p_0, \nonumber
\end{align}
 with $b(p(s))\neq 0$ and $p(s)\in\Omega$ for $0\leq s \leq s'$ with some $s'>0$, then we call $\Gamma:=\{p(s)\,|$ $0\leq s \leq s'\}$ a \emph{drift trajectory} starting at $p_0$.\\
\end{itemize}
\end{Definition}
\textbf{Remark:} Drift trajectories are oriented at $p_0$ in the direction $b(p_0)$; they do not run both ways. Diffusion trajectories are not oriented, they run in both directions. In our special case of a diagonal $D$, each diffusion trajectory moves along one of the canonical unit vectors in $\im D$.\\
\\
Next, we introduce the propagation set:
\begin{Definition}
 Let $\Omega\subset\R^{d+1}$. Two points $p,q\in\R^{d+1}$ are connected by a diffusion trajectory in $\Omega$ iff there is some diffusion trajectory $\Gamma\subset\Omega$ with $p,q\in\Gamma$. $q$ is connected to $p$ by a drift trajectory in $\Omega$ iff there is a drift trajectory $\Gamma\subset\Omega$ starting at $p$ with $q\in\Gamma$.\\
 For any point $p\in\Omega$, the \emph{propagation set} $S(p,\Omega)$ consists of all $q\in\Omega$ that are connected to $p$ by a finite series of drift and diffusion trajectories.
\end{Definition}
Again, note that drift trajectories are oriented and can only connect to points backward in time. Therefore, it is possible that $q\in S(p,\Omega)$ while $p\notin S(q,\Omega)$.\\
With this notation, we can restate the \emph{interior maximum principle} from Theorem 1 of \cite{Hi70}:
\begin{Theorem}
 \label{CDHmax}
Let $p=(t,x)\in\Omega\subset\R^+\times\R^{d}$. Let the function $f\in C^2(\Omega)$ satisfy $\tilde Lf\leq0$ on the propagation set $S(p,\Omega)$ and
\begin{align*}
 \inf\limits_{S(p,\Omega)} f \geq 0.
\end{align*}
If $f(p)=0$, then $f=0$ in $\overline{S(p,\Omega)}$.
\end{Theorem}
The propagation set corresponding to equation (\ref{linmasterequ}) can be characterised as follows:
\begin{Lemma}
 \label{PropSet}
Let $p=(t,x)\in\Omega:=\R^+\times\R^{d}$. Then $S(p,\Omega)=[0,t)\times\R^d\cup\ \{(t,x_0)\}\times\R^k$, where $x_0$ is the orthogonal projection of $x$ onto the kernel of $D$ (restricted to $\R^{d-k}$).
\end{Lemma}
While elementary, the proof of Lemma \ref{PropSet} is somewhat lengthy and deferred to the appendix.\\

To see that the solution $f$ of (\ref{linmasterequ}) fulfils $f\geq0$, one can employ the same method used for the (standard) weak maximum principle for non-degenerate parabolic equations. Now we give the proof of Theorem \ref{globalpos} via the sharp maximum principle from \cite{Hi70}:\\
\\
\begin{Proof}[ (of Theorem \ref{globalpos})]
Let $f\in C^2(\R^+\times\R^d)$ be a solution to (\ref{linmasterequ}) for some $f_0\geq0$ with $\int\limits_{\R^d} f_0 \d =1$. Then
\begin{align*}
 g(t)&:=e^{-\tr(C)t}f(t)
\end{align*}
solves
\begin{align}
\label{nonnegativefequ} g_t - \div(D\nabla g) - (Cx)^T\nabla g  &= 0,\\
\nonumber g(t=0)=f_0\geq 0\,.
\end{align}
The standard maximum principle then shows that $g(x,t)\geq0$ for $t\geq0$. {}From (\ref{tildecomparison}), we have $\tilde L g=0$.
Assume $g(t',x')=0$ for some $t'>0$, $x'\in\R^d$. Then Theorem \ref{CDHmax} gives $g=0$ on $[0,t']\times\R^d$ and in particular $f_0=0$. But this is a contradiction to
\begin{align*}
 \int\limits_{\R^d} f_0(x)\d &= 1.
\end{align*}
Hence, $f(t,x)>0$ for all $t>0, x\in\R^d$.
\end{Proof}

\section{Existence of a steady state, decomposition of the generator $L$}
\label{sec:steadystate}

In light of Theorem \ref{globalpos}, we are looking for a steady state $f_\infty$ of (\ref{linmasterequ}) that fulfils the conditions
\begin{align}
\label{ssconditions} \int\limits_{\R^d} f_\infty(x) \d &= 1,\\
\nonumber	f_\infty &> 0.
\end{align}
In fact, the existence of such a steady state is equivalent to condition (A):
\begin{Theorem}(Existence of a steady state)
\label{ssexistence}
There exists a unique steady state $f_\infty\in L^1(\R^d)$ of (\ref{linmasterequ}) fulfilling (\ref{ssconditions}) iff condition (A) holds.\\
Moreover, this steady state is of the (non-isotropic) Gaussian form
\begin{align*}
 f_\infty(x)&=c_K\exp(-\fract{x^TK^{-1}x}2),
\end{align*}
where $K$ is the unique, symmetric, and positive definite solution to the continuous Lyapunov equation
\begin{align}
\label{Dequ} 2D=CK+KC^T,
\end{align}
and $c_K=(2\pi)^{-\frac{d}{2}} (\det K)^{-\frac{1}{2}}$ is the normalisation constant.
\end{Theorem}
For the proof of Theorem \ref{ssexistence} we consider the Fourier transform of (\ref{linmasterequ}):
\begin{align}
 \label{masterFequ} \hat f_t(t,\xi) &= -(\xi^TD\xi) \hat f(t,\xi) - (C^T\xi)\cdot(\nabla_\xi \hat f(t,\xi)), \\
\nonumber \hat f(t=0)&= \hat f_0.
\end{align}
A steady state $f_\infty\in L^1(\R^d)$ implies $\hat f_\infty\in C_0(\R^d)$. Also note that
\begin{align*}
 \hat f_\infty(0)=\int\limits_{\R^d} f_\infty(x)\d = 1
\end{align*}
for the normalised steady state.\\
\\
Thus, the steady state equation in Fourier space reads
\begin{align}
\label{FourierCharEqu} 0 &= (\xi^TD\xi)\hat f_\infty(\xi) + (C^T\xi)\cdot\nabla_\xi \hat f_\infty(\xi),\\
 \nonumber \hat f_\infty(0) &= 1.
\end{align}

The problem at hand is closely related to the stationary Fokker-Planck equation in §2.2 in \cite{ArCaMa10}. But for $k<d$, the singularity of $D$ requires a more careful analysis.\\

We will split the proof of Theorem \ref{ssexistence} into three lemmas: in Lemmas \ref{LemmaSS1} and \ref{LemmaSS2} we establish that existence of a steady state is equivalent to condition (A). Lemma \ref{LemmaSS3} establishes that the steady state is Gaussian.  
\begin{Lemma}
\label{LemmaSS1}
 Let (\ref{FourierCharEqu}) have a unique solution $\hat f_\infty\in C_0(\R^d)$. Then condition (A) holds.
\end{Lemma}
\begin{Proof}
First, we shall show that $C^T$ is regular: if $C^T$ had a non-trivial kernel, (\ref{FourierCharEqu}) restricted to the kernel of $C^T$ would read
\begin{align}
\label{restrequ}
  \forall \xi\in\ker C^T:(\xi^TD\xi)\hat f_\infty(\xi) &=0.
\end{align}
Now, either $\ker C^T\subset \ker D$, which would mean that both drift and diffusion in (\ref{linmasterequ}) only act on a proper subspace of $\R^d$ and there would be no unique steady state; or (\ref{restrequ}) implies
\begin{align*}
 &\exists v\in\R^d:\ \forall s\in\R\backslash\{0\}:\hat f_\infty(sv)=0.
\end{align*}
Hence, $f_\infty(0)=0$ by continuity, which is a contradiction to $\hat f_\infty(0)=1$. So $C^T$ is regular.\\

Next, we will show that $C$ is positively stable, i.e. that all eigenvalues have a strictly positive real part. The characteristic equations for (\ref{FourierCharEqu}) are
\begin{align}
\label{characteristicequations}
 \dot \xi(s) &= C^T\xi(s),\quad s\in\R,\\
 \nonumber \dot z(s) &= -(\xi(s)^TD\xi(s))z(s),\quad s\in\R, \\
  \nonumber (z(0),\xi(0)) &= (z_0,\xi_0)\in\R^{d+1}.
\end{align}
The solutions to these equations are
\begin{align*}
 \xi(s) &= e^{C^Ts}\xi_0,\\
 z(s) &= z_0\exp\left(-\int\limits_0^s \xi(\tau)^TD\xi(\tau) \d[\tau]\right).
\end{align*}
Assume that $C$ has an eigenvalue $\lambda$ with $\Re\{\lambda\}<0$. Let $v$ be a corresponding eigenvector of $C^T$, i.e. $C^Tv=\lambda v$, chosen such that $v\notin i\R^d$. Consider the characteristic curve starting at $\xi_0:=v+\bar v\neq 0$:
\begin{align*}
 \xi(s)&=e^{\lambda s}v+e^{\bar \lambda s}\bar v\\
\Rightarrow |\xi(s)| &= e^{\Re\{\lambda\}s} |v+e^{2i\Im\{\lambda\}s}\bar v| \to \infty,\ s\to-\infty\\
\Rightarrow \forall s\leq0:|z(s)| &= |z_0|\exp\left(\int\limits_s^0 \xi(\tau)^TD\xi(\tau) \d[\tau]\right) \geq |z_0|,
\end{align*}
due to $D$ being positive semidefinite. If $z_0\neq0$, this is a contradiction to $|\hat f_\infty(\xi)|\to0$, $|\xi|\to\infty$. If $z_0=0$, we can take the limit $s\to\infty$ and obtain a contradiction to $\hat f_\infty(0)=1$ and the continuity of $\hat f_\infty$. So $C$ cannot have eigenvalues with negative real part.\\

Now assume that $C$ has a purely imaginary eigenvalue. Then there exist characteristics $\xi(s)$ which form circles. Due to
\begin{align*}
 z(s) &= z_0\exp\left(-\int\limits_0^s \xi(\tau)^TD\xi(\tau) \d[\tau]\right)
\end{align*}
and the continuity of $\hat f_\infty$, one of the following statements has to hold on any such characteristic curve:
\begin{itemize}
 \item[(a)] $\forall s\in\R: \xi(s)\in\ker D$,
 \item[(b)] $z_0=0$.
\end{itemize}
If (a) holds, then we have $z(s)=z_0$ on this characteristic. Since the characteristic is closed, there will be no uniqueness of $\hat f_\infty$.\\
So (b) holds, and for any $\eps$ we can find such a characteristic starting at a vector $\xi_0$ with $|\xi_0|<\eps$. But then $\hat f_\infty(\xi(s))=z_0=0$, which is a contradiction to the continuity of $\hat f_\infty$ at $0$.\\

This shows that $C$ has to be positively stable. It remains to show the first part of condition (A). We employ the reformulation of Lemma \ref{Definiteness} (ii). So assume $C^T$ has an eigenvector $v$ with $Dv=0$. Then $D\bar v=0$, and for the characteristic starting at $\xi(0)=v+\bar v$, we have
\begin{align*}
 \xi(s)&=e^{\lambda s}v+e^{\bar \lambda s}\bar v,\\
 z(s)&= z(0)\exp\left(-\int\limits_0^s (e^{\lambda \tau}v+e^{\bar \lambda \tau}\bar v)^TD(e^{\lambda \tau}v+e^{\bar \lambda \tau}\bar v) \d[\tau]\right) = z(0).
\end{align*}
This means that $z$ is constant on the characteristic $\xi$. Now, since $C$ is positively stable,
\begin{align*}
 \lim\limits_{s\to\infty} |\xi(s)| &= \infty,\\
 \lim\limits_{s\to-\infty} |\xi(s)| &= 0.\\	
\end{align*}
So we would need $z(0)=1$ because of the continuity in $0$, and $z(0)=0$ because $\hat f_\infty\in C_0(\R^n)$. That is a contradiction, so there can be no eigenvector $v$ of $C^T$ with $Dv=0$.
\end{Proof}

\begin{Lemma}
\label{LemmaSS3}
 Let $C$ be positively stable. Then, the function
\begin{align*}
 \hat f_\infty(\xi) &:= \exp(-\fract{\xi^TK\xi}2)
\end{align*}
is a solution to (\ref{FourierCharEqu}), where $K\geq0$ is the unique solution of (\ref{Dequ}).\\
Furthermore, $K$ is regular iff no eigenvector $v$ of $C^T$ satisfies $Dv=0$. In this case, $f_\infty$ is Gaussian and hence in $L^1(\R^d)$.
\end{Lemma}
\begin{Proof}
 We insert the ansatz
\[
 \hat f_\infty(\xi) = \exp(-\fract{\xi^T K\xi}2)
\]
with a symmetric matrix $K\in\R^d$ into (\ref{FourierCharEqu}) and get
\begin{align*}
  \forall \xi\in\R^d: 0 &= \Big(\xi^TD\xi  - (C^T\xi)\cdot(K\xi)\Big) \hat f_\infty \\
\Leftrightarrow \forall \xi\in\R^d:0 &= \xi^T (D - CK)\xi \\
\Leftrightarrow &D - CK \text{ is antisymmetric}\\
\Leftrightarrow &D - CK = KC^T-D,\
\end{align*}
and this is equivalent to (\ref{Dequ}). This continuous Lyapunov equation has a unique, symmetric and positive semidefinite solution $K$ since $C$ is positively stable (see, e.g., Theorem 2.2 in \cite{SnZaN70}, Theorem 2.2.3 in \cite{HoJoT91}).\\

Now assume that $K$ is not regular. Then there is a $v\neq0$ with $Kv=0$ and (\ref{Dequ}) implies
\begin{align*}
 2v^TDv = v^TCKv+v^TKC^Tv = 0\,.
\end{align*}
Due to $D=D^T=D^2$ and \eqref{Dequ} this implies
$$
  0 = 2Dv = CKv + KC^Tv = KC^Tv, 
$$
so $C^Tv$ is also an eigenvector of $K$ to the eigenvalue $0$. Since $v\neq 0$ and $C^T$ is regular, $C^Tv\neq0$. Repeating this calculation with $C^Tv$ instead of $v$, we can see that $C^Tv$ is in the kernel of $D$, and thus $(C^T)^2v$ is in the kernel of $K$. A proof by induction then gives $(C^T)^kv\in \ker D\cap\ker K$ for all $k\in\N$. Therefore, the space
\begin{align*}
 V &:= \operatorname{span}[v,\dots,(C^T)^{d-1}v]
\end{align*}
is a $C^T$-invariant subspace of $\ker D$. So $K$ is regular if there is no eigenvector $v$ of $C^T$ with $Dv=0$.\\
For the reversed implication, assume that there is an eigenvector $v$ of $C^T$ (corresponding to the eigenvalue $\lambda_v$) with $Dv=0$. This implies
\begin{align*}
 0 &= 2v^TDv = v^TCKv+v^TKC^Tv = \overline{\lambda_v} v^TKv + \lambda_v v^TKv\\
   &= 2\Re\{\lambda_v\} v^TKv.
\end{align*}
Since $\Re\{\lambda_v\}>0$ for all eigenvalues of $C^T$, it follows that $v^TKv=0$ and thus $K$ is not regular.
\end{Proof}

\begin{Lemma}
\label{LemmaSS2}
 Let condition (A) hold. Then the steady state $f_\infty$ from Lemma \ref{LemmaSS3} is unique.
\end{Lemma}
\begin{Proof}
We will show that the characteristic equations (\ref{characteristicequations}) have a unique solution fulfilling (\ref{ssconditions}). As the starting manifold for the characteristics, we take $\Gamma:=\{\xi_0\in\R^d:|\xi_0|=1\}$, which is admissible since $C$ is positively stable. The characteristic curve starting at $\xi_0$ is
\begin{align*}
 \xi(s)=e^{C^Ts}\xi_0.
\end{align*}
Since $C^T$ is positively stable, we have
$$
 \lim\limits_{s\to\infty} |\xi(s)| = \infty,\quad
 \lim\limits_{s\to-\infty} |\xi(s)| = 0,
$$
and the characteristic curves cover all of $\R^d$.\\
The value of solutions along the characteristics is
\begin{align*}
 z(s) &= z(0)\exp\left(-\int\limits_0^s \xi(\tau)^TD\xi(\tau)\d[\tau]\right).
\end{align*}
So taking
\begin{align*}
 z(0)&=\exp\left(-\int\limits_{-\infty}^0 \xi(\tau)^TD\xi(\tau)\d[\tau]\right)
\end{align*}
as initial condition implies $1=\lim\limits_{s\to-\infty} z(s)=\hat f_\infty(0)$. Since $\xi(s)$ decays exponentially for $s\to-\infty$, $z(0)$ is always finite and there is a unique solution $z(s)$.
\end{Proof}
\\
This lemma completes the proof of Theorem \ref{ssexistence}.\\

In analogy to the entropy method for linear, non-degenerate Fokker-Planck equations presented in \cite{ArMaToUn01}, we now consider (\ref{linmasterequ}) in the weighted space $L^2:=L^2(\R^d,f_\infty^{-1})$ with inner product $\<\cdot,\cdot\>$. On this space, the operator $L=\div(D\nabla \cdot +Cx\cdot)$ can be decomposed very naturally.\\
\begin{Theorem}
\label{OPsplit}
 Let (\ref{linmasterequ}) fulfil condition (A). Consider $L$ on the weighted space $L^2$. Then $L$ can be decomposed into its symmetric part $L_s$ and its antisymmetric part $L_{as}$ as
\begin{align}
\nonumber L_sf 
&= \div(D\nabla(\fract{f}{f_\infty})f_\infty), \\
\label{Lasdefinition} L_{as}f 
&= \div(R\nabla(\fract{f}{f_\infty})f_\infty).
\end{align}
Here, $R:=\frac12(CK-KC^T)$ is antisymmetric, $K$ is the covariance matrix of $f_\infty$ from Theorem \ref{ssexistence}. 
\end{Theorem}
\tb{Remarks:} \begin{itemize}
\item[(i)] Note that the steady state $f_\infty$ fulfils both $L_sf_\infty=0$ and $L_{as}f_\infty=0$.
\item[(ii)] $R\neq0$ and hence (\ref{linmasterequ}) is non-symmetric in $L^2$. Otherwise (\ref{Dequ}) would imply $D=KC^T$ and $\ker D=\ker C^T$, which contradicts condition (A).\\
\end{itemize}

\begin{Proof}[ (of Theorem \ref{OPsplit})]\\
 We compute
\begin{align*}
 \<Lf,g\>&=\ild (Lf) g \exp(\fract{x^TK^{-1}x}2)\d \\
	 &=-\ild [D\nabla f + Cxf]\cdot[\nabla g+K^{-1}xg]\exp(\fract{x^TK^{-1}x}2)\d \\
	 &=\ild f\div [(D\nabla g+DK^{-1}xg)\exp(\fract{x^TK^{-1}x}2)]\d - \ild f x^TC^T [\nabla g+K^{-1}xg] \exp(\fract{x^TK^{-1}x}2)\d \\
	 &=\ild f [\div(D\nabla g+DK^{-1}xg)+x^TK^{-1}D\nabla g+x^TK^{-1}DK^{-1}xg-x^TC^T\nabla g-x^TC^TK^{-1}xg]\\
	 &\;\;\;\times \exp(\fract{x^TK^{-1}x}2)\d. \\
\end{align*}
Using (\ref{Dequ}), we have $K^{-1}DK^{-1}-C^TK^{-1}=K^{-1}RK^{-1}$. Since $R$ is antisymmetric it follows that $x^T(K^{-1}DK^{-1}-C^TK^{-1})x = 0$ and hence
\begin{align*}
   L^*g&= \div(D\nabla g+DK^{-1}xg)+x^TK^{-1}D\nabla g-x^TC^T\nabla g.
\end{align*}
Furthermore, $\tr(DK^{-1}-C)=\tr((D-CK)K^{-1})=0$, since $D-CK=-R$ is antisymmetric and $K^{-1}$ is symmetric. Thus we can write (using (\ref{Dequ}) in the last step)
\begin{align*}
 L^*g&=\div(D\nabla g+DK^{-1}xg+(DK^{-1}-C)xg)\\
     &=\div(D\nabla g+(2DK^{-1}-C)xg)\\
     &= \div(D\nabla g+(KC^TK^{-1}x)g).
\end{align*}
So we get, again using (\ref{Dequ}),
\begin{align*}
 L_sf &= \fract{L+L^*}2f \\
      &= \div(D\nabla f+\fract12(C+KC^TK^{-1})xf)= \div(D\nabla f+DK^{-1}xf)\\
      &= \div(D\nabla(\fract{f}{f_\infty})f_\infty);\\
 L_{as}f&=\fract{L-L^*}2f \\
      &=\div(\fract12 (C-KC^TK^{-1})xf)=\div(RK^{-1}xf) \\
      &=\div(R\nabla(\fract{f}{f_\infty})f_\infty),\\
\end{align*}
where we have used $\div(R\nabla f)=0$ for the last equality.
\end{Proof}


\section{Entropy method, explicit decay rate}
\label{sec:modentmod}

In this section, we will prove an explicit decay rate for the solution $f$ of (\ref{linmasterequ}) under condition (A). To do so, we consider relative entropies, as in \cite{ArMaToUn01}. We will see that, unlike in the fully parabolic case, a direct entropy-entropy dissipation estimate cannot be obtained. Instead, we prove exponential decay of an auxiliary functional that bounds the entropy dissipation. This still implies a decay rate for the relative entropy, initially at the price of additional regularity requirements on the initial state $f_0$. A regularisation result adapted from \cite{ViH06} is then employed to obtain the sharp decay rate for solutions with finite initial entropy. The sharpness of this rate will be shown in the next section.\\

With the notations of §3 we introduce the relative entropy:
\begin{Definition}
\label{admissableentropy}
Let $0\not\equiv\psi\in C(\R^+_0)\cap C^4(\R^+)$, $\psi(1)=\psi'(1)=0$, $\psi''\geq0$ on $\R^+$, $(\psi''')^2\leq\fract12\psi''\psi^{IV}$ on $\R^+$. Let $f\in L^1_+(\R^d)$ with $\int f\d=1$. Then
\begin{align*}
 e_\psi(f|f_\infty)&:=\int\limits_{\R^d} \psi(\fract{f}{f_\infty})f_\infty \d
\end{align*}
is called an \emph{admissible relative entropy} with \emph{generating function} $\psi$.
\end{Definition}

The most important examples are the logarithmic entropy $e_1(f|f_\infty)$ with $\psi_1(s) = s\ln s-s+1$ and the quadratic entropy $e_2(f|f_\infty)$ with $\psi_2(s)=(s-1)^2$. For the latter we admit $f\in L^1(\R^d)$ and hence we consider $\psi_2$ on $\R$. $e_1$ and $e_2$ are also the limiting cases of admissible relative entropies (cf.\ \S2.2 of \cite{ArMaToUn01}).

The entropy method is based on computing a bound on the first two time-derivatives of the relative entropy $e_\psi(f(t)):=e_\psi(f(t)|f_\infty)$ with $f$ the solution to (\ref{linmasterequ}). Formally,
\begin{align}
\label{entropyderivative} \ddt e_\psi(f(t)) &= -\int\limits_{\R^d} \psi''(\fu) \nabla(\fu)^TD\nabla(\fu) f_\infty \d =: -I_\psi(f) \leq 0.
\end{align}
However, there may be a technical problem if $f(t,x)=0$ (which can happen at the initial state $f_0$). For example, $\psi_1''(s)=\fract1s$ 
, and this would lead to a division by zero. For this reason, we use a trick from \cite{ArMaToUn01} (see Remark 2.12) to rewrite (\ref{entropyderivative}):

\begin{Definition}
\label{wdefinition}
Let $\psi$ generate an admissible entropy, and let $f_0\in L^1_+(\R^d)$ (or $f_0\in L^1(\R^d)$ for quadratic $\psi$) with $\int\limits_{\R^d} f_0 \d =1$. Define
\begin{align}
\label{wdef}
 w(x):=\int\limits_{1}^{\fract{f_0}{f_\infty}(x)} \sqrt{\psi''(s)}\d[s].
\end{align}
Then we call $f_0$ a \emph{$\psi$-compatible initial state} iff $\nabla w\in L^2(\R^d,f_\infty)$.
\end{Definition}

With Definition \ref{wdefinition}, (\ref{entropyderivative}) can be written as
\begin{align}
\label{Irewritten} I_\psi(f) &=\int\limits_{\R^d} (\nabla w)^TD(\nabla w) f_\infty \d.
\end{align}
Whenever $f\neq0$, this is equivalent to (\ref{entropyderivative}); however, now there is no longer a problem when $f=0$. So the assumption of Definition \ref{wdefinition} clearly implies that the initial state has finite entropy dissipation. It also has finite relative entropy, as we shall prove in Proposition \ref{Sfiniteefinite} below.\\
\\
\tb{Remark:} The integral in Definition \ref{wdefinition} can be calculated explicitly for the most common entropies:\\
For the quadratic entropy, $\psi_2(s)=\alpha(s-1)^2$ for some $\alpha>0$, and thus
\begin{align}
\label{quadraticw} w &= \sqrt{2\alpha}(\fract{f_0}{f_\infty}-1).
\end{align}
For the logarithmic entropy, with
\begin{align}
\label{logentropydefinition} \psi_1(s)=\alpha(s+\beta)\ln(\fract{s+\beta}{1+\beta})-\alpha(s-1)
\end{align}
for some $\alpha>0$, $\beta\geq0$, we have
\begin{align}
\label{logw} w &= 2\sqrt{\alpha}(\sqrt{\fract{f_0}{f_\infty}+\beta}-\sqrt{1+\beta}).
\end{align}
For the $p$-entropies, $1<p<2$, $\psi_p(s)=\alpha[(s+\beta)^p - (1+\beta)^p - p(1+\beta)^{p-1}(s-1)]$ for some $\alpha>0$, $\beta\geq0$, and thus
\begin{align*}
 w &= 2\sqrt{\fract{\alpha(p-1)}p}((\fract{f_0}{f_\infty}+\beta)^{\fract p2}-(1+\beta)^{\fract p2}).
\end{align*}
\\
There is another, in fact systematic problem with the entropy dissipation (\ref{entropyderivative}): Since $D$ is singular for $k<d$, this functional is `lacking information' on some partial derivatives of $\fu$. But this information would be vital for the (standard) entropy method to work. More precisely, the functional $I_\psi$ vanishes \emph{not only} for $f=f_\infty$. As shown in Corollary \ref{ZeroTangentStart}, for any $t^*\geq0$ there are initial conditions such that $I_\psi(f(t^*))=0$. Also, due to the monotonicity of $e_\psi(f(t))$, $I_\psi(f(t^*))=0$ for some $t^*\geq0$ implies $I_\psi'(f(t^*))=0$. So, for degenerate Fokker-Planck equations, $e_\psi(f(t))$ is not a convex function of $t$ -- in contrast to the non-degenerate case from \cite{ArMaToUn01}. The possibility of having $I_\psi(f(t^*))=I_\psi'(f(t^*))=0$ for $f(t^*)\neq f_\infty$ also shows that the standard entropy method cannot be carried over directly to the degenerate case in (\ref{linmasterequ}).\\
We therefore introduce the modified functional 
\begin{align}
\label{Sdefinition} S_\psi(f) &:= \int\limits_{\R^d} (\nabla w)^TP(\nabla w)f_\infty \d = \int\limits_{\fu>0} \psi''(\fu) \nabla(\fu)^TP\nabla(\fu)f_\infty\d, 
\end{align}
where we replace the matrix $D$ in $I_\psi$ with a symmetric, positive definite matrix $P$. $P$ will be chosen in such a way that it allows for an estimate between $\ddt S_\psi(f(t))$ and $S_\psi(f(t))$ for solutions $f$ to (\ref{linmasterequ}), as shown later in this section. Moreover, since $P$ is positive definite, there is a constant $c_P>0$ with $P\geq c_PD$, and hence $S_\psi\geq c_PI_\psi$.\\
\\
\textbf{Remark:} 
Introducing the functional $S_\psi$ differs from the modified entropy dissipation approach in \cite{DeH06}. There one considers an ``intermediate functional'' $K(f)$, which measures the distance of $f$ to the set of stationary states of the symmetric part ($L_{s}$ in our case).  
\\

Choosing the matrix $P$ is the crucial ingredient for the definition of our modified entropy dissipation $S_\psi$:
\begin{Lemma}
\label{Pdefinition}
 Let $Q:=KC^TK^{-1}$, with $K$ from \eqref{Dequ}. Let $\mu:=\min\left\{\Re\{\lambda\}|\lambda\text{ is an eigenvalue of } C\right\}$. Due to condition (A), $\mu>0$. Let $\{\lambda_{m}|1\leq m\leq m_0\}$ be all the eigenvalues of $C$ with $\mu=\Re\{\lambda_m\}$, only counting their geometric multiplicity.\\
 \begin{itemize}
 \item[(i)] If $\lambda_m$ is non-defective\footnote{An eigenvalue is defective if its geometric multiplicity is strictly less than its algebraic multiplicity.} for all $m\in\{1,\dots,m_0\}$, then there exists a symmetric, positive definite matrix $P\in\R^{d\times d}$ with
\begin{align}
\label{matrixestimate} QP+PQ^T &\geq 2\mu P.
\end{align}
\item[(ii)] If $\lambda_m$ is defective for at least one $m\in\{1,\dots,m_0\}$, then for any $\eps>0$ there exists a symmetric, positive definite matrix $P=P(\eps)\in\R^{d\times d}$ with
\begin{align}
\label{degeneratematrixestimate} QP+PQ^T &\geq 2(\mu-\eps) P.
\end{align}
\item[(iii)] For any such matrix $P$, and for any $\psi$-compatible function $f_0$, $S_\psi(f_0)<\infty$.
 \end{itemize}
\end{Lemma}
\begin{Proof}
The idea behind the construction of $P$ is the following: If $Q$ is not defective (and hence diagonalizable) and $w_1,\dots,w_d$ are its eigenvectors, then one can choose $P$ as the weighted sum of the following rank 1 matrices:
\begin{align}
\label{simpleP} P:= \sum\limits_{j=1}^d b_j \,w_j\otimes \overline{w_j}^T\,,
\end{align}
with $b_j\in\R^+$, $j=1,\dots,d$. As $\{w_j\}_{j=1,\dots,d}$ is a basis of $\C^d$, $P$ is positive definite. If any $w_j$ is complex, its complex conjugate $\overline{w_j}$ is also an eigenvector of $Q$, since $Q$ is real. By taking the same coefficient $b_j$ for both, we obtain a real matrix $P$. Apart from this restriction, the choice of $b_j>0$ is arbitrary.
For $P$ from (\ref{simpleP}), we obtain
\begin{align*}
 QP+PQ^T &= 
 \sum\limits_{j=1}^d b_j (\lambda_j+\overline{\lambda_j}) w_j\otimes \overline{w_j}^T 
         \geq 2\mu \sum\limits_{j=1}^d b_j \,w_j\otimes \overline{w_j}^T = 
				 2\mu P\,.
\end{align*}

If at least one of the eigenvalues of $Q$ is defective, one can still construct $P$ in a similar fashion to \eqref{simpleP}, but including now the generalised eigenvectors of $Q$.
To this end, we consider the Jordan normal form $J$ of $Q$, given by the similarity transformation $A^{-1}QA=J$ with some $A\in\C^{d\times d}$. Let $J$ have $N$ Jordan blocks, each of length $l_n$; $n=1,\dots,N$.\\
(i) By assumption, all Jordan blocks corresponding to eigenvalues with $\Re\{\lambda_n\}=\mu$ are trivial, i.e.\ of length 1. Corresponding to the structure of $J$, we define the positive diagonal matrix
\begin{align*}
 B := \diag(B_1,\dots,B_N),
\end{align*}
with
\begin{align*}
  B_n := \diag(b^{l_n}_n,\dots,b^1_n)\,,\quad n=1,\dots,N\,.
\end{align*}
Its entries are defined as
\begin{align}
\label{bnjdefinition} 
  b_n^1:=1;\quad b_n^j:=c_j\,(\tau_n)^{2(1-j)}\,;\quad j=2,\dots,l_n\,,
\end{align}
where  $c_1:=1$, $c_j:=1+(c_{j-1})^2;\,j=2,\dots,l_n$, 
and $\tau_n:=2(\Re\{\lambda_n\}-\mu)\geq0$ for $n=1,\dots,N$. This yields for the $n$-th Jordan block $J_n$ in the case $l_n=1$: $B_n=1$ and
\begin{align*}
 J_nB_n + B_nJ_n^H = (\lambda_n + \overline{\lambda_n}) B_n \geq 2\mu B_n.
\end{align*}
Here, $J_n^H$ denotes the Hermitian adjoint of $J_n$. In the case $l_n>1$, we have $\tau_n>0$ and
\begin{align*}
 J_nB_n + B_nJ_n^H-2\mu B_n&= \left(\begin{array}{cccc}  2(\Re\{\lambda_n\}-\mu)b_n^{l_n} & b^{l_n-1}_n & & \\ b^{l_n-1}_n & 2(\Re\{\lambda_n\}-\mu)b_n^{l_n-1} & \ddots & \\ & \ddots  & \ddots & b^1_n \\ & & b^1_n & 2(\Re\{\lambda_n\}-\mu)b^1_n \end{array}\right) \geq 0.
\end{align*}
The last inequality follows from
\begin{align*}
 M_m:= \left(\begin{array}{cccc}  \tau^{3-2m} c_{m} & \tau^{4-2m}c_{m-1} & & \\ \tau^{4-2m}c_{m-1} & \tau^{5-2m} c_{m-1} & \ddots & \\ & \ddots  & \ddots & c_1 \\ & & c_1 & \tau c_1 \end{array}\right) \geq 0,\quad m=1,\dots,\max_n(l_n),
\end{align*}
for any $\tau>0$, which can be verified by the principal minor test and the recursion
$$
  \det M_{m} = \tau^{3-2m} c_{m}\det M_{m-1} - (\tau^{4-2m}c_{m-1})^2\det M_{m-2} 
  = \tau^{m(2-m)}>0\,\quad\mbox{for }m\ge3\,.
$$

In total, we have $JB+BJ^H \geq 2\mu B$, and hence
\begin{align*}
 A^{-1}QA B+B A^HQ^T(A^{-1})^H \geq 2\mu B,
\end{align*}
which implies
\begin{align*}
 QABA^H+ABA^HQ^T &\geq 2\mu ABA^H.
\end{align*}
The claim then follows with $P:=ABA^H$.\\
(ii) In this case, there exists a non-trivial Jordan block $J_{\tilde n}$ corresponding to an eigenvalue with $\Re\{\lambda_{\tilde n}\}=\mu$. In \eqref{bnjdefinition} of the above construction, we then choose (instead of $\tau_n$) $\tau_{\tilde n}:=2(\Re\{\lambda_{\tilde n}\}-\mu+\eps)>0$ for some $\eps>0$. 
Hence, $J_{\tilde n}B_{\tilde n}+B_{\tilde n}J_{\tilde n}^H\geq 2(\mu-\eps)B_{\tilde n}$ and the result follows. However, in this case $P$ depends on $\eps$.\\
(iii) Using $P\le c\Id$ with some $c>0$, this is clear from (\ref{Sdefinition}).
\end{Proof}\\
\tb{Remarks:}
\begin{itemize}
\item[(i)] The matrix $P$ in Lemma \ref{Pdefinition} is not uniquely determined (in general not even up to a multiplicative factor; see the construction \eqref{simpleP}).
\item[(ii)] {}From (\ref{bnjdefinition}) with $\tau_{\tilde n}:=2(\Re\{\lambda_{\tilde n}\}-\mu+\eps)$, we see that for a defective eigenvalue $\lambda_{\tilde n}$,
\begin{align*}
  \forall 1< j\leq l_{\tilde n}:\quad \lim\limits_{\eps\to0} b^j_{\tilde n} &= \infty\,.
\end{align*}
 With this ``scaling'' of $P=P(\eps)$ we thus have (for general $f_0$)
 \begin{align*}
 \lim\limits_{\eps\to0} S_\psi(f_0,\eps) &= \infty,
\end{align*}
 with $S_\psi(f_0,\eps):=\ild (\nabla w)^T P(\eps) \nabla w f_\infty\d$.
 An alternative ``scaling'' of $P$ would be to multiply with an appropriate power of $\eps$, to keep $S_\psi(f_0,\eps)$ bounded. But then, $P$ would be singular in the $\eps\to0$ limit.
\item[(iii)] To appreciate the matrix inequality (\ref{matrixestimate}) we multiply it with $\sqrt{P}^{\,-1}$ from both sides:
$$
  \sqrt{P}^{\,-1} Q \sqrt{P} + \sqrt{P} Q^T \sqrt{P}^{\,-1}\ge 2\mu \Id\,.
$$
With the similarity transformation $\tilde Q:=\sqrt{P}^{\,-1} Q \sqrt{P}$ we have 
$\mu=\min\,\{\Re\{\lambda\}\,|\,\lambda\in\sigma(\tilde Q)\}$, and the above inequality reads
$$
  \tilde Q_s \ge\mu\Id\,.
$$
But note that, in general, we would have the opposite inequality for the smallest eigenvalue of the symmetric part of a matrix. This motivates that the choice $P=\Id$ will not work in general.
\item[(iv)] (\ref{matrixestimate}) can be rewritten as $(Q-\mu)P+P(Q^T-\mu)\geq 0$, which bears a close resemblance to the continuous Lyapunov equation from Theorem \ref{ssexistence}. If we assume equality in (\ref{matrixestimate}) and if $Q-\mu$ were positively stable, then there would be a unique solution $P=0$, see e.g. \cite{HoJoT91}. But since $\mu$ is the real part of an eigenvalue of $Q$, $Q-\mu$ is not positively stable. This explains why we can find a non-trivial solution of (\ref{matrixestimate}) at the price of uniqueness. 

There is equality in (\ref{matrixestimate}) iff all eigenvalues of $Q$ have the same real part $\mu$ and are non-degenerate. For additional details, we refer to \cite{HoJoT91}, \cite{SnZaN70}.\\
\end{itemize}

Next we show that any $\psi$-compatible $f$ (or equivalently $S_\psi(f)<\infty$) also has finite relative entropy generated by $\psi$:
\begin{Proposition}
\label{Sfiniteefinite}
 Let $f$ be $\psi$-compatible. Then it holds that
 \begin{align}\label{convSobinequ}  
  e_\psi(f|f_\infty) &\leq \frac1{2\lambda_P} S_\psi(f) < \infty\,,
 \end{align}
where $\lambda_P>0$ is the largest possible constant in the matrix inequality $K^{-1}\ge \lambda_P P^{-1}$.
\end{Proposition}
\begin{Proof}
 Consider the Fokker-Planck operator
 \begin{align*}
  L_Pf &:= \div(P\nabla(\fu)f_\infty)
 \end{align*}
 on $L^2$. Then $L_P$ is symmetric due to the symmetry of $P$, and $f_\infty$ spans the kernel of $L_P$. One easily checks that
 \begin{align*}
  \ddt e_\psi(\tf(t)|f_\infty) &= - S_\psi(\tf(t))
 \end{align*}
 for a solution $\tf(t)$ to $\tf_t = L_P\tf$. As shown in Corollary 2.17, \cite{ArMaToUn01}, this symmetric, non-degenerate Fokker-Planck equation leads to an exponential decay of the relative entropy, and in parallel to a convex Sobolev inequality: Using the notation $f_\infty(x):=c_Ke^{-V(x)}$, $V(x):=\frac{x^TK^{-1}x}2$, we have the Bakry-Émery condition
 \begin{align*}
  \frac{\partial^2 V}{\partial x^2} &= K^{-1} \geq \lambda_P P^{-1}\,.
 \end{align*}
Hence, all $g\in L^1_+(\R^d)$ with $\ild g \d =1$ satisfy the convex Sobolev inequality
 \begin{align*}
e_\psi(g|f_\infty) &\leq \frac1{2\lambda_P} S_\psi(g),
 \end{align*}
 where both sides may be infinite. Since $f$ is $\psi$-compatible, we have $S_\psi(f)<\infty$, $f\in L^1_+(\R^d)$ and $\int\limits_{\R^d} f \d =1$. This completes the proof.
\end{Proof}\\
\\
The strategy of the standard entropy method is to prove first the exponential decay of the entropy dissipation. In analogy, we shall prove first the decay of the modified entropy dissipation $S_\psi$. Afterwards, this will yield the decay of $f(t)$ in relative entropy.
\begin{Proposition}
 \label{Sconvergence}
 Assume condition (A). Let $\psi$ generate an admissible entropy and let $f$ be the solution to (\ref{linmasterequ}) with a $\psi$-compatible initial state $f_0$, $\mu:=\min\left\{\Re\{\lambda\}|\lambda\text{ is an eigenvalue of } C\right\}$. Let $P$, $S_\psi(f_0)$ be defined as in Lemma \ref{Pdefinition}, $\{\lambda_m|1\leq m\leq m_0\}$ be the eigenvalues of $C$ with $\mu=\Re\{\lambda_m\}$.
 \begin{itemize}
  \item[(i)] If all $\lambda_m$, $1\leq m\leq m_0$, are non-defective, then
\begin{align*}
 S_\psi(f(t)) &\leq S_\psi(f_0)e^{-2\mu t},\quad t\geq 0.
\end{align*}
  \item[(ii)] If $\lambda_m$ is defective for at least one $m\in\{1,\dots,m_0\}$, then
\begin{align*}
 S_\psi(f(t),\eps) &\leq S_\psi(f_0,\eps)e^{-2(\mu-\eps)t},\quad t\geq 0,
\end{align*}
for any $\eps\in(0,\mu)$.
 \end{itemize}
\end{Proposition}
\tb{Remark:} This result holds for all matrices $P$ chosen according to Lemma \ref{Pdefinition}. Clearly, the rate $\mu$ is independent of the choice of $P$.\\
\\
\begin{Proof}[ (of Proposition \ref{Sconvergence})]Due to Proposition \ref{solexistence}, the solution $f$ is sufficiently smooth to allow the following computations. They are inspired by the decay estimate for the entropy dissipation in non-degenerate Fokker-Planck equations (cf. Lemma 2.13 of \cite{ArMaToUn01}). Due to the global positivity shown in §2, the solution remains $\psi$-compatible for all $t>0$. Let $u:=\nabla\fract{f}{f_\infty}$ and $S_\psi$ be given as in Lemma \ref{Pdefinition}. Then
\begin{align*}
 u_t&=-\frac{\partial^2 V}{\partial x^2}(D+R)u-\frac{\partial u}{\partial x}(D-R)\nabla V+(\nabla^TD\nabla) u.
\end{align*}
We compute
\begin{align*}
Z_\psi(f(t))&:= \ddt S_\psi(f(t)) \\
    &= \underbrace{2\int\limits_{\R^d} \psi''(\fract{f}{f_\infty}) (u_t)^TPuf_\infty \d}_{=:(I)}+ \underbrace{\int\limits_{\R^d} \psi'''(\fract{f}{f_\infty}) u^TPu f_t \d}_{=:(II)},
\end{align*}
where we have used the symmetry of $P$. We have
\begin{align*}
 (I)&= 2\int\limits_{\R^d} \psi''(\fract{f}{f_\infty}) (u_t)^TPuf_\infty \d \\
    &= -2\int\limits_{\R^d} \psi''(\fract{f}{f_\infty}) u^T(D-R)\frac{\partial^2 V}{\partial x^2}Pu f_\infty\d - 2\int\limits_{\R^d} \psi''(\fract{f}{f_\infty}) (\nabla V)^T(D+R)\frac{\partial u}{\partial x}Pu f_\infty \d\\
    &\quad + 2\int\limits_{\R^d} \psi''(\fract{f}{f_\infty}) ((\nabla^TD\nabla) u)^TPu f_\infty \d.\\
\end{align*}
For the last term, we compute (using the summation convention over double indices)
\begin{align*}
    &\quad  2\int\limits_{\R^d} \psi''(\fract{f}{f_\infty}) (\nabla^TD\nabla u)^TPu f_\infty \d \\
    &= 2 \int\limits_{\R^d} \psi''(\fract{f}{f_\infty}) D_{lk}u_{j,kl}P_{jr}u_r f_\infty \d \\
    &= -2\int\limits_{\R^d} D_{lk}P_{jr} u_{j,k} (\psi''(\fract{f}{f_\infty})u_rf_\infty)_{,l} \d \\
    &= -2\int\limits_{\R^d} \psi'''(\fract{f}{f_\infty}) D_{lk}u_{j,k} u_lP_{jr}u_rf_\infty \d - 2\int\limits_{\R^d} \psi''(\fract{f}{f_\infty}) D_{lk}u_{j,k} P_{jr}u_{r,l} f_\infty \d\\
    &\quad + 2\int\limits_{\R^d} \psi''(\fract{f}{f_\infty}) D_{lk}u_{j,k} V_{,l}P_{jr}u_r f_\infty \d \\
    &= -2\int\limits_{\R^d} \psi'''(\fract{f}{f_\infty}) u^TD\frac{\partial u}{\partial x}Pu f_\infty \d - 2\int\limits_{\R^d} \psi''(\fract{f}{f_\infty}) \tr(D\frac{\partial u}{\partial x}P\frac{\partial u}{\partial x}) f_\infty \d \\
    &\quad + 2\int\limits_{\R^d} \psi''(\fract{f}{f_\infty}) (\nabla V)^TD\frac{\partial u}{\partial x}Pu f_\infty\d,
\end{align*}
where we have used $u_{k,j}=u_{j,k}$ in the last equality. We obtain
\begin{align*}
(I)&= -2\int\limits_{\R^d} \psi''(\fract{f}{f_\infty}) (\nabla V)^TR\frac{\partial u}{\partial x}Pu f_\infty \d  - 2\int\limits_{\R^d} \psi''(\fract{f}{f_\infty}) \tr(D\frac{\partial u}{\partial x}P\frac{\partial u}{\partial x}) f_\infty \d  \\
&\quad -  2\int\limits_{\R^d} \psi''(\fract{f}{f_\infty}) u^T(D-R)\frac{\partial^2 V}{\partial x^2}Pu f_\infty\d - 2\int\limits_{\R^d} \psi'''(\fract{f}{f_\infty}) u^TD\frac{\partial u}{\partial x}Pu f_\infty \d.
\end{align*}
Next, we rewrite the first term of this formula:
\begin{align*}
 &\quad \int\limits_{\R^d} \psi''(\fract{f}{f_\infty}) (\nabla V)^TR\frac{\partial u}{\partial x}Pu f_\infty \d \\
 &= \int\limits_{\R^d} \psi''(\fract{f}{f_\infty}) V_{,l}R_{lk}u_{j,k}P_{jr}u_r f_\infty \d \\
 &= -\int\limits_{\R^d} R_{lk}u_jP_{jr}(V_{,l}u_rf_\infty\psi''(\fract{f}{f_\infty}))_{,k} \d \\
 &= -\int\limits_{\R^d} \psi''(\fract{f}{f_\infty})R_{lk}V_{,lk} u_jP_{jr}u_r f_\infty \d - \int\limits_{\R^d} \psi''(\fract{f}{f_\infty}) V_{,l}R_{lk}u_{r,k}P_{jr}u_j f_\infty \d \\
 &\ +\int\limits_{\R^d} \psi''(\fract{f}{f_\infty}) V_{,l}R_{lk}V_{,k} u_jP_{jr}u_r f_\infty\d -\int\limits_{\R^d} \psi'''(\fract{f}{f_\infty}) V_{,l}R_{lk}u_k u_jP_{jr}u_r f_\infty \d\\
 &= -\int\limits_{\R^d} \psi''(\fract{f}{f_\infty}) (\nabla V)^TR\frac{\partial u}{\partial x}Pu f_\infty \d - \int\limits_{\R^d} \psi'''(\fract{f}{f_\infty}) [(\nabla V)^TRu] [u^TPu] f_\infty \d.\\
\end{align*}
Here we have used the skew-symmetry of $R$ to conclude $R_{lk}V_{,lk}=0$, $V_{,l}R_{lk}V_{,k}=0$. Hence,
\begin{align*}
  -2\int\limits_{\R^d} \psi''(\fract{f}{f_\infty}) (\nabla V)^TR\frac{\partial u}{\partial x}Pu f_\infty \d &= \int\limits_{\R^d} \psi'''(\fract{f}{f_\infty}) [(\nabla V)^TRu] [u^TPu] f_\infty \d.
\end{align*}
So we arrive at
\begin{align*}
(I)&= -2\int\limits_{\R^d} \psi''(\fract{f}{f_\infty}) \tr(D\frac{\partial u}{\partial x}P\frac{\partial u}{\partial x}) f_\infty \d - 2\int\limits_{\R^d} \psi''(\fract{f}{f_\infty}) u^T(D-R)\frac{\partial^2 V}{\partial x^2}Pu f_\infty\d \\
&\quad - 2\int\limits_{\R^d} \psi'''(\fract{f}{f_\infty}) u^TD\frac{\partial u}{\partial x}Pu f_\infty \d + \int\limits_{\R^d} \psi'''(\fract{f}{f_\infty}) [(\nabla V)^TRu][u^TPu] f_\infty \d.
\end{align*}
Next, we compute
\begin{align*}
 (II)&= \int\limits_{\R^d} \psi'''(\fract{f}{f_\infty}) u^TPu f_t \d\\
      &= \int\limits_{\R^d} \psi'''(\fract{f}{f_\infty}) u^TPu \div(f_\infty(D+R)u) \d\\
      &= \int\limits_{\R^d} \psi'''(\fract{f}{f_\infty}) u_rP_{rj}u_j (f_\infty(D_{lk}+R_{lk})u_k)_{,l} \d \\
      &= -\int\limits_{\R^d} \psi'''(\fract{f}{f_\infty}) u_rP_{rj}u_j V_{,l}(D_{lk}+R_{lk})u_k f_\infty\d + \int\limits_{\R^d} \psi'''(\fract{f}{f_\infty}) u_rP_{rj}u_j (D_{lk}+R_{lk})u_{k,l} f_\infty \d.
\end{align*}
Take a closer look at
\begin{align*}
      &\quad \int\limits_{\R^d} \psi'''(\fract{f}{f_\infty}) u_rP_{rj}u_j (D_{lk}+R_{lk})u_{k,l} f_\infty\d\\
      &= \int\limits_{\R^d} \psi'''(\fract{f}{f_\infty}) u_rP_{rj}u_jD_{lk}u_{k,l} f_\infty\d\\
      &= -\int\limits_{\R^d} u_kD_{lk}P_{rj}(u_ru_jf_\infty \psi'''(\fract{f}{f_\infty}))_{,l} \d\\
      &= -\int\limits_{\R^d} \psi'''(\fract{f}{f_\infty}) u_kD_{lk}u_{r,l}P_{rj}u_j f_\infty \d - \int\limits_{\R^d} \psi'''(\fract{f}{f_\infty}) u_kD_{lk}u_{j,l}P_{rj}u_r f_\infty\d \\
      &\quad - \int\limits_{\R^d} \psi^{IV}(\fract{f}{f_\infty}) u_kD_{lk}u_l u_rP_{rj}u_j f_\infty\d + \int\limits_{\R^d} \psi'''(\fract{f}{f_\infty}) u_kD_{lk}V_{,l} u_rP_{rj}u_j f_\infty \d,
\end{align*}
and it follows that
\begin{align*}
(II)&= -2\int\limits_{\R^d} \psi'''(\fract{f}{f_\infty}) u^TD\frac{\partial u}{\partial x}Pu f_\infty \d - \int\limits_{\R^d} \psi^{IV}(\fract{f}{f_\infty}) [u^TDu] [u^TPu] f_\infty \d \\
   &\quad - \int\limits_{\R^d} \psi'''(\fract{f}{f_\infty}) [(\nabla V)^TRu] [u^TPu] f_\infty \d.
\end{align*}
With this, we obtain
\begin{align*}
 Z_\psi(f(t))&=-2\int\limits_{\R^d} \psi''(\fract{f}{f_\infty}) u^T(D-R)\frac{\partial^2 V}{\partial x^2}Pu f_\infty \d -2 \int\limits_{\R^d} \psi''(\fract{f}{f_\infty}) \tr(D\frac{\partial u}{\partial x}P\frac{\partial u}{\partial x}) f_\infty \d \\
&\quad -4 \int\limits_{\R^d} \psi'''(\fract{f}{f_\infty}) u^TD\frac{\partial u}{\partial x}Pu f_\infty \d - \int\limits_{\R^d} \psi^{IV}(\fract{f}{f_\infty}) [u^TPu][u^TDu] f_\infty \d \\
	     &= -2\ild \tr(XY) f_\infty \d - \ild \psi''(\fract{f}{f_\infty}) u^T[(D-R)\frac{\partial^2 V}{\partial x^2}P+P\frac{\partial^2 V}{\partial x^2}(D+R)]u f_\infty \d.
\end{align*}
Here, the matrices $X$, $Y$ are given as (cf. Lemma 2.13 in \cite{ArMaToUn01})
\begin{align*}
 X = \left(\begin{array}{cc} \psi''(\fu) & \psi'''(\fu) \\ \psi'''(\fu) & \frac12 \psi^{IV}(\fu) \end{array}\right),&\quad Y=\left(\begin{array}{cc} \tr(D\frac{\partial u}{\partial x}P\frac{\partial u}{\partial x}) & u^TD\frac{\partial u}{\partial x}Pu \\ u^TD\frac{\partial u}{\partial x}Pu & (u^TPu)(u^TDu) \end{array}\right).
\end{align*}
Due to the assumptions on $\psi$ (cf. Definition \ref{admissableentropy}), $X\geq0$. To see $Y\geq0$, we use the Cauchy-Schwarz inequality for the Hilbert-Schmidt norm and the symmetry of $D$, $P$ to obtain
\begin{align*}
 (u^TD\frac{\partial u}{\partial x}Pu)^2 &= \tr\left(\sqrt{P}uu^T\sqrt{D}\cdot \sqrt{D}\frac{\partial u}{\partial x}\sqrt{P}\right)^2 \\
		    &\leq \tr\left(\sqrt{P}uu^T\sqrt{D}\sqrt{D}uu^T\sqrt{P}\right)\tr\left(\sqrt{D}\frac{\partial u}{\partial x}\sqrt{P}\sqrt{P}\frac{\partial u}{\partial x}\sqrt{D}\right) \\
		    &= [u^TDu][u^TPu] \;\tr\left(D\frac{\partial u}{\partial x}P\frac{\partial u}{\partial x}\right).
\end{align*}
This implies $\tr(XY)\geq0$, and thus
\begin{align}
\label{ZSinequality} Z_\psi(f(t)) &\leq - \int\limits_{\R^d} \psi''(\fract{f}{f_\infty}) u^T[(D-R)\frac{\partial^2 V}{\partial x^2}Pu f_\infty+P\frac{\partial^2 V}{\partial x^2}(D+R)]u f_\infty\d.
\end{align}
We can now use Lemma \ref{Pdefinition} to establish the decay rate for $S_\psi$. First, compute
\begin{align*}
 (D-R)\frac{\partial^2V}{\partial x^2} &= (D-R)K^{-1} = \fract12(CK+KC^T-CK+KC^T)K^{-1}\\
	&= KC^TK^{-1}=Q, \\
 \frac{\partial^2V}{\partial x^2}(D+R) &= K^{-1}CK=Q^T\,,
\end{align*}
with $Q$ from Lemma \ref{Pdefinition}. So the right hand side of (\ref{ZSinequality}) reads $-\ild\psi''(\fu)u^T(QP+PQ^T)uf_\infty\d$. In Lemma \ref{Pdefinition} we proved 
\begin{align}\label{P-ineq}
  QP+PQ^T\geq 2\kappa P\,, 
\end{align}
where $\kappa=\mu$ for case (i), and $\kappa=\mu-\eps$ for case (ii) with some $\eps\in(0,\mu)$. Thus
\begin{align*}
 \ddt S_\psi(f(t)) &\leq -2\kappa S_\psi(f(t)),\\
\end{align*}
and applying Gronwall's lemma completes the proof.
\end{Proof}\\
\\
\tb{Remark:} Inequality \eqref{P-ineq} is the key ingredient of the above proof, and it is a direct generalization of the well known Bakry-\'Emery condition from the standard entropy method. Indeed, for $D=\Id$ and $C=C^T\ge\mu>0$, \eqref{linmasterequ} is a symmetric Fokker-Planck equation. With $K^{-1}=Q=C$ one can choose $P=\Id$. Then, \eqref{P-ineq} reads
$$
  \frac12 (QP+PQ^T) = C =\frac{\partial^2 V}{\partial x^2} \ge \mu\Id\,,
$$
and it is the Bakry-\'Emery condition in its simplest form (cf.\ (A2) in \cite{ArMaToUn01}).

For $D=\Id$, and $C\ne C^T$ normal and positively stable, \eqref{linmasterequ} is a non-symmetric Fokker-Planck equation with $K^{-1}=C_s:=(C+C^T)/2$ and $Q=C^T$. Here, the Bakry-\'Emery condition reads $K^{-1}=C_s\ge \lambda_K\Id$, while inequality \eqref{P-ineq} yields the improvement
$$
  C^TP+PC\ge 2\mu P\,,
$$
with $\mu=\min\{\Re(\lambda)\;|\,\lambda\in\sigma(C)\}$. We always have $\mu\ge\lambda_K=\min \lambda(C_s)$ and the strict inequality holds in many examples.
We shall return to this comparison for non-symmetric Fokker-Planck equations in \S\ref{sec:nonsymmFP}.\\

In the standard entropy method for fully parabolic equations, one derives decay of the relative entropy from the decay of the entropy dissipation by integrating the inequality $\frac{\operatorname{d}^2}{\operatorname{ds}^2}e_\psi(f(s))\geq -\dds e_\psi(f(s))$ over $(t,\infty)$. This requires a-priori knowledge that $e_\psi(f(t=\infty))=0$, which, as shown in \cite{ArCaJuL08}, can be derived from the decay of $S_\psi$ (which is the entropy dissipation functional for fully parabolic equations). Since the decay estimate in Proposition \ref{Sconvergence} depends on the modified entropy dissipation of the initial data, this time integration does not work for hypocoercive equations.

Still, the convex Sobolev inequality from Proposition \ref{Sfiniteefinite} already implies exponential decay of the relative entropy under the assumption that $S_\psi(f_0)<\infty$:
\begin{Theorem}\label{entropydecay}
 Assume condition (A). Let $\psi$ generate an admissible entropy and let $f$ be the solution to (\ref{linmasterequ}) with a $\psi$-compatible initial state $f_0$, $\mu:=\min\left\{\Re\{\lambda\}|\lambda\text{ is an eigenvalue of } C\right\}$. Let $P$, $S_\psi(f_0)$ be defined as in Lemma \ref{Pdefinition}, $\{\lambda_m|1\leq m\leq m_0\}$ be the eigenvalues of $C$ with $\mu=\Re\{\lambda_m\}$.
 \begin{itemize}
  \item[(i)] If all $\lambda_m$, $1\leq m\leq m_0$, are non-defective, then
\begin{align}\label{entropydecay-nondeg}
 e_\psi(f(t)|f_\infty) &\leq \frac{1}{2\lambda_P} S_\psi(f_0)e^{-2\mu t},\quad t\geq 0.
\end{align}
  \item[(ii)] If $\lambda_m$ is defective for at least one $m\in\{1,\dots,m_0\}$, then
\begin{align}\label{entropydecay-deg}
 e_\psi(f(t)|f_\infty) &\leq \frac{1}{2\lambda_P} S_\psi(f_0,\eps)e^{-2(\mu-\eps)t},\quad t\geq 0,
\end{align}
for any $\eps\in(0,\mu)$.
 \end{itemize}
\end{Theorem}

\begin{Remark}[sharpness of the constants in Theorem \ref{entropydecay}]
\normalfont
\label{const-remark}
\
\begin{itemize}
\item[(i)] 
While the l.h.s.\ of \eqref{entropydecay-nondeg} is independent of $P$, the r.h.s.\ clearly depends on $P$. So, \emph{for each fixed $f_0$}, the multiplicative constant $\frac{1}{2\lambda_P} S_\psi(f_0)$ can be optimised w.r.t.\ the admissible matrices $P$ from Lemma \ref{Pdefinition} (in the family \eqref{simpleP}, e.g.). The same statement applies to the defective case of \eqref{entropydecay-deg} (for each fixed $\eps>0$).

\item[(ii)]
But \emph{for each fixed $P$}, the leading multiplicative constant $\frac{1}{2\lambda_P}$ in \eqref{entropydecay-nondeg} and \eqref{entropydecay-deg} is sharp for the logarithmic and quadratic entropies. This is understood in the sense that (for each $P$) there exists an \emph{optimal function} rendering the convex Sobolev inequality \eqref{convSobinequ} an equality. (The coordinate transformation $x=\sqrt{P}y$ changes $\tilde f_t=L_P\tilde f$ into a Fokker-Planck equation with $\Id$ as diffusion matrix such that the sharpness results from \S3.5 in \cite{ArMaToUn01} apply.) Using this optimal function as $f_0$ hence makes \eqref{entropydecay-nondeg} and \eqref{entropydecay-deg} an equality at $t=0$.

\item[(iii)]
The sharpness of the exponential rates in \eqref{entropydecay-nondeg} and \eqref{entropydecay-deg}  (for all admissible entropies) will be proved in \S\ref{sec:sharprate} below.

\item[(iv)]
For the case $d=2$, the combined optimality of rate and multiplicative constant will be shown in Proposition \ref{sharp-constant}.\\
\end{itemize}
\end{Remark}

Using the regularisation of \eqref{linmasterequ} we shall next generalise the entropy decay to initial states with (only) finite relative entropy. 
The basic concept is that evolutions with hypoelliptic operators regularise, though in a weaker sense than non-degenerate parabolic equations. Local estimates of this sort first appeared in the proof by Hörmander \cite{Ho67} as well as in \cite{KoH73}, \cite{RoStH77}. Our result generalises Theorems A.12, A.15 in \cite{ViH06} (expressed for quadratic and logarithmic entropies) to all admissible $\psi$-entropies. Those results, in turn, used an idea developed by Hérau \cite{HeF07}. The regularisation depends on the order $\tau$ of the \emph{finite rank Hörmander condition} for $L$ (cf. Remark \ref{tauremark}).

\begin{Theorem}
\label{entropyregularisation}
 Let condition (A) hold, $f_0\in L^1_+(\R^d)$ with $\ild f_0 \d =1$ and $e_\psi(f_0|f_\infty)<\infty$. Let $f(t)$ be the solution of (\ref{linmasterequ}) with initial condition $f_0$, and let $\tau$ be the minimal constant such that \eqref{sumTdefinite} holds. Then there is a positive constant $c_r>0$ such that
 \begin{align}
 \label{regularityestimate} \forall t\in(0,1]:\ S_\psi(f(t)) &\leq c_rt^{-(2\tau+1)}e_\psi(f_0|f_\infty).
 \end{align}
\end{Theorem}
\begin{Proof}
 The idea of the proof is to construct a decaying-in-time functional $\F$ that is a (positive) linear combination of both sides of (\ref{regularityestimate}) -- multiplied by $t^{2\tau+1}$.\\
 \underline{Step 1 (construction of $\F$):} With $Q=KC^TK^{-1}$ from Lemma \ref{Pdefinition}, we define the matrices
 \begin{align*}
  M_j := Q^jD(Q^T)^j\geq 0,&\quad N_j := Q^jD(Q^T)^{j+1} + Q^{j+1}D(Q^T)^j;\quad j=0,\dots,\tau+2\,.
 \end{align*}
 Since $Q^T=K^{-1}CK=2K^{-1}D-C^T$, we can apply \eqref{sumTdefinite} to $\sum\limits_{j=0}^{\tau} M_j$ and obtain
 \begin{align*}
  \sum\limits_{j=0}^{\tau} M_j \geq c_0\Id  
 \end{align*}
 for some $c_0>0$. Thus there is $c_1>0$ such that
 \begin{align}
  \label{Mtaup1bounded} M_{\tau+2} &\leq c_1 \sum\limits_{j=0}^{\tau} M_j.
 \end{align}
 We compute
 \begin{align}
\label{QedMj}  QM_j+M_jQ^T &= N_j,\\
\label{QedNj}  QN_j+N_jQ^T &= 2M_{j+1} + Q^jD(Q^T)^{j+2} + Q^{j+2}D(Q^T)^j.
 \end{align}
 Using $D^2=D$, we have for any $\eps>0$:
 \begin{align}
\nonumber  0&\leq \left(\frac{1}{\sqrt\eps}Q^jD \pm \sqrt\eps Q^{j+2}D\right)\left(\frac1{\sqrt\eps}D(Q^T)^j\pm\sqrt\eps D(Q^T)^{j+2}\right) \\
\label{MQestimate}  &= \frac1\eps M_j + \eps M_{j+2} \pm \left(Q^jD(Q^T)^{j+2} + Q^{j+2}D(Q^T)^j\right).
 \end{align}
 Then (\ref{QedMj}) and the analogue of (\ref{MQestimate}) with $j+2$ replaced by $j+1$ yield the estimate
 \begin{align}
\label{Njestimate}  \pm N_j &\leq \frac1\eps M_j + \eps M_{j+1}.
 \end{align}
 Further, (\ref{MQestimate}) yields
 \begin{align}
  \label{Mjestimate}  \pm \left(Q^jD(Q^T)^{j+2} + Q^{j+2}D(Q^T)^j\right) &\leq \frac1\eps M_j + \eps M_{j+2}.
 \end{align}
Now we define the matrix-valued polynomial in $t$:
 \begin{align*}
  P(t) := \sum\limits_{j=0}^{\tau+1} \left(a_jt^{2j+1}M_j\right) &+ \sum\limits_{j=0}^{\tau} \left(b_jt^{2j+2}N_j\right),
 \end{align*}
 with $P(0)=0$. As (positive) coefficients, we first choose  $a_{\tau+1} := \frac1{c_1}$,
 \begin{align*}
   b_{\tau} := \frac23\left[1+a_{\tau+1}(2\tau+4)\right],&\quad  a_{\tau} := 2\frac{b_{\tau}^2}{a_{\tau+1}}.
 \end{align*}
 Then we choose iteratively, starting with $j=\tau$ and finishing with $j=1$:
 \begin{align}
 \label{Ptcoefficientdefinition} b_{j-1} :=\frac23\left[2+c_1+a_j(2j+1)+b_j^2+\frac{2(b_j(2j+2)-a_j)^2}{b_j}\right],&\quad a_{j-1} := 8\frac{b_{j-1}^2}{a_j}.
 \end{align}
 Using (\ref{Njestimate}) with $\eps=\frac{2b_jt}{a_j}$, $0\leq j\leq \tau$, we obtain
 \begin{align*}
  \forall j=0,\dots,\tau:\quad b_jt^{2j+2}N_j &\geq  -\frac{a_j}2t^{2j+1}M_j - \frac{2b_j^2}{a_j} t^{2j+3}M_{j+1},
 \end{align*}
 and thus
 \begin{align*}
  \sum\limits_{j=0}^\tau b_jt^{2j+2}N_j &\geq -\frac{a_0}2tM_0- \sum\limits_{j=1}^\tau \left([\frac{a_j}2+\frac{2b_{j-1}^2}{a_{j-1}}]t^{2j+1}M_j\right) - \frac{2b_{\tau}^2}{a_\tau}t^{2\tau+3}M_{\tau+1} \\
	    &= -\frac{a_0}2tM_0- \sum\limits_{j=1}^\tau \left(\frac{3a_j}4t^{2j+1}M_j\right) - a_{\tau+1}t^{2\tau+3}M_{\tau+1},
 \end{align*}
 where we have used (\ref{Ptcoefficientdefinition}). Inserting this into $P(t)$ yields
 \begin{align*}
  P(t) &\geq \frac{a_0}2tM_0+\sum\limits_{j=1}^{\tau} \frac{a_j}4t^{2j+1}M_j\,.
 \end{align*}
 Writing $c_3:=\min\{\frac{a_0}2,\frac{a_1}4,\dots,\frac{a_\tau}4\}$, this implies for $t\in[0,1]$:
 \begin{align}
 \label{Ptdefinite} P(t) &\geq  t^{2\tau+1} c_3 \sum\limits_{j=0}^{\tau}  M_j \geq c_0c_3 t^{2\tau+1}\Id.
 \end{align}
 So $P(t)$ is positive definite for all $t>0$, and we define the functional
 \begin{align*}
  \F(t) &:= \gamma e_\psi(f(t)|f_\infty) + \ild \psi''(\fu) u^TP(t)u f_\infty \d\geq0,
 \end{align*}
 with some $\gamma>0$ to be chosen later.\\
 \\
 \underline{Step 2 (decay of $\F$):} For $\F$, we can repeat all the computations in the proof of Proposition \ref{Sconvergence} and arrive at
 \begin{align*}
  \ddt\F(t) &\leq -\gamma I_\psi(f(t)|f_\infty) + \ild \psi''(\fu) u^T\left[\dot P(t) - \left(QP(t)+P(t)Q^T\right)\right]u f_\infty \d \\
  &= \ild \psi''(\fu) u^T\left[\dot P(t) - \left(QP(t)+P(t)Q^T\right) - \gamma M_0\right]u f_\infty \d\,,
 \end{align*}
where $\dot P(t)$ denotes the time derivative of $P(t)$. We compute
 \begin{align*}
  \dot P(t) &= \sum\limits_{j=0}^{\tau+1} \left(a_j(2j+1)t^{2j}M_j\right) + \sum\limits_{j=0}^{\tau} \left(b_j(2j+2)t^{2j+1}N_j\right)
 \end{align*}
 and further, using (\ref{QedMj}), (\ref{QedNj}), and (\ref{Mjestimate}) with $\eps:=\frac{t^2}{b_j}$:
 \begin{align*}
  -\left(QP(t)+P(t)Q^T\right) &=   -\sum\limits_{j=0}^{\tau+1} \left(a_jt^{2j+1}N_j\right) - 2\sum\limits_{j=0}^{\tau} \left(b_jt^{2j+2}M_{j+1}\right)   \\
    &\quad -\sum\limits_{j=0}^{\tau} \left(b_jt^{2j+2}[Q^jD(Q^T)^{j+2} + Q^{j+2}D(Q^T)]\right) \\
  &\leq -\sum\limits_{j=0}^{\tau+1} \left(a_jt^{2j+1}N_j\right) - 2\sum\limits_{j=0}^{\tau} \left(b_jt^{2j+2}M_{j+1}\right) \\
			 &\quad + \sum\limits_{j=0}^{\tau} b_jt^{2j+2}\left(\frac{b_j}{t^2}M_j+\frac{t^2}{b_j}M_{j+2}\right) \\
			 &= -\sum\limits_{j=0}^{\tau+1} \left(a_jt^{2j+1}N_j\right) - 2\sum\limits_{j=0}^{\tau} \left(b_jt^{2j+2}M_{j+1}\right) \\
			 &\quad + \sum\limits_{j=0}^{\tau} \left(t^{2j}b_j^2M_j\right) + \sum\limits_{j=2}^{\tau+2} \left(t^{2j}M_j\right).
 \end{align*}
 This implies
 \begin{align*}
  &\quad \dot P(t) - \left(QP(t)+P(t)Q^T\right) - \gamma M_0 \\
  &\leq  \left(a_0+b_0^2-\gamma\right)M_0 + \left(3a_1+b_1^2-2b_0\right)t^{2}M_1 \\
  &\quad + \sum\limits_{j=2}^{\tau} \left([a_j(2j+1)+1+b_j^2-2b_{j-1}]t^{2j}M_j\right) \\
  &\quad + (a_{\tau+1}(2\tau+3)+1-2b_{\tau})t^{2\tau+2}M_{\tau+1} \\
  &\quad + \sum\limits_{j=0}^{\tau+1} \left(\alpha_jt^{2j+1}N_j\right) +t^{2(\tau+2)}M_{\tau+2},
 \end{align*}
 where $\alpha_j:=-a_j+b_j(2j+2)$, $0\leq j\leq\tau$; $\alpha_{\tau+1}:=-a_{\tau+1}$. Using
 \begin{align*}
  \forall j=0,\dots,\tau:\quad \pm N_j &\leq \frac{2|\alpha_j|}{b_jt}M_j+\frac{b_jt}{2|\alpha_j|}M_{j+1},\\
  N_{\tau+1} &\leq \frac1tM_{\tau+1}+tM_{\tau+2},
 \end{align*}
 we obtain
 \begin{align*}
  \sum\limits_{j=0}^{\tau+1} \left(\alpha_jt^{2j+1}N_j\right) &\leq  \sum\limits_{j=0}^{\tau}\left(\frac{2\alpha_j^2}{b_j}t^{2j}M_j+\frac{b_j}2t^{2j+2}M_{j+1} \right) +  a_{\tau+1}t^{2\tau+2}M_{\tau+1} + a_{\tau+1}t^{2\tau+4}M_{\tau+2} \\
  &= \frac{2\alpha_0^2}{b_0}M_0 + \sum\limits_{j=1}^{\tau} \left(\frac{2\alpha_j^2}{b_j}+\frac{b_{j-1}}2\right)t^{2j}M_j + (\frac{b_\tau}2+a_{\tau+1})t^{2\tau+2}M_{\tau+1} + a_{\tau+1}t^{2\tau+4}M_{\tau+2}.
 \end{align*}
 Thus, we finally arrive at
 \begin{align*}
   &\quad \dot P(t) - \left(QP(t)+P(t)Q^T\right) - \gamma M_0 \\
  &\leq  \left(a_0+b_0^2+\frac{2\alpha_0^2}{b_0}-\gamma\right)M_0 + \left(3a_1+b_1^2+\frac{2\alpha_1^2}{b_1}+\frac{b_0}{2}-2b_0\right)t^{2}M_1 \\
  &\quad + \sum\limits_{j=2}^{\tau} \left(a_j(2j+1)+1+b_j^2+\frac{2\alpha_j^2}{b_j}+\frac{b_{j-1}}2-2b_{j-1}\right)t^{2j}M_j \\
  &\quad + (a_{\tau+1}(2\tau+4)+1+\frac{b_\tau}2-2b_{\tau})t^{2\tau+2}M_{\tau+1} + (a_{\tau+1}+1)t^{2\tau+4}M_{\tau+2}.
 \end{align*}
 We use (\ref{Mtaup1bounded}) and obtain for sufficiently large $\gamma$ and $t\in[0,1]$:
 \begin{align*}
   &\quad \dot P(t) - \left(QP(t)+P(t)Q^T\right) - \gamma M_0 \\
  &\leq  \left(c_1(a_{\tau+1}+1)+a_0+b_0^2+\frac{2\alpha_0^2}{b_0}-\gamma\right)M_0 + \left(c_1(a_{\tau+1}+1)+3a_1+b_1^2+\frac{2\alpha_1^2}{b_1}-\frac{3b_0}2\right)t^{2}M_1 \\
  &\quad + \sum\limits_{j=2}^{\tau} \left(\left[c_1(a_{\tau+1}+1)+a_j(2j+1)+1+b_j^2+\frac{2\alpha_j^2}{b_j}-\frac{3b_{j-1}}2\right]t^{2j}M_j\right) \\
  &\quad + (a_{\tau+1}(2\tau+4)+1-\frac{3b_\tau}2)t^{2\tau+2}M_{\tau+1} \leq 0,
 \end{align*}
 where we have used (\ref{Ptcoefficientdefinition}).\\
 
 This implies that $\F(t)$ is monotonously decreasing, and thus $\F(t)\leq \F(0)=\gamma e_\psi(f_0|f_\infty)$ for all $t$ in $[0,1]$. Together with (\ref{Ptdefinite}), we obtain
 \begin{align*}
  c_0c_3 t^{2\tau+1} \ild \psi''(\fu) |u|^2 f_\infty \d &\leq \gamma e_\psi(f_0|f_\infty),
 \end{align*}
 which completes the proof using Lemma \ref{Pdefinition} (iii).
\end{Proof}

With this regularisation result, we can finally prove exponential decay of the relative entropy:

\begin{Theorem}
 \label{convergencerate}  Assume condition (A). Let $\psi$ generate an admissible relative entropy and let $f$ be the solution to (\ref{linmasterequ}) with initial state $f_0\in L^1_+(\R^d)$ such that $e_\psi(f_0|f_\infty)<\infty$. Let $\mu:=\min\{\Re\{\lambda\}|\lambda\text{ is an}$ $\text{eigenvalue of } C\}$. Let $\{\lambda_m|1\leq m\leq m_0\}$ be the eigenvalues of $C$ with $\mu=\Re\{\lambda_m\}$, and let
 \begin{align*}
  e(t) &:= e_\psi(f(t)|f_\infty).
 \end{align*}
 Then
 \begin{itemize}
  \item[(i)] If all $\lambda_m$, $1\leq m\leq m_0$, are non-defective, then there is a constant $c>1$ such that
  \begin{align}\label{entropy-decay-nondeg}
   \forall t\geq 0:\quad e(t) &\leq ce^{-2\mu t}e_\psi(f_0|f_\infty).
  \end{align}
  \item[(ii)] If $\lambda_m$ is defective for at least one $m\in\{1,\dots,m_0\}$, then for all $\eps\in(0,\mu)$, there is $c_\eps>1$ such that
  \begin{align}\label{entropy-decay-deg}
   \forall t\geq0:\quad e(t) &\leq c_\eps e^{-2(\mu-\eps)t}e_\psi(f_0|f_\infty).
  \end{align}
 \end{itemize}
\end{Theorem}
\begin{Proof}
Let $P$, $S_\psi(f_0)$ be defined as in Lemma \ref{Pdefinition}. Let $\delta>0$, and let $\kappa:=\mu$ in case (i), and $\kappa:=\mu-\eps$ in case (ii). Using (\ref{convSobinequ}), Proposition \ref{Sconvergence} and Theorem \ref{entropyregularisation}, we compute for $t\geq\delta$:
\begin{align}
\nonumber e_\psi(t) &\leq \frac1{2\lambda_P}S_\psi(f(t))\leq \frac1{2\lambda_P}S_\psi(f(\delta))e^{-2\kappa(t-\delta)} \\
\label{tgdestimate} &\leq e^{2\kappa\delta}\frac{c_r}{2\lambda_P\delta^{2\tau+1}} e_\psi(0)e^{-2\kappa t}.
\end{align}
For $t\leq \delta$, it follows from the monotonicity of $e_\psi$ (cf.\ (\ref{entropyderivative})) that
\begin{align}
\label{tldestimate} e_\psi(t) &\leq e_\psi(0)\,. 
\end{align}
Writing $c_\delta:=e^{2\kappa\delta}\max\{1,\frac{c_r}{2\lambda_P\delta^{2\tau+1}}\}$ and combining (\ref{tgdestimate}), (\ref{tldestimate}) yields
\begin{align*}
\forall t\geq 0:\quad e_\psi(t) &\leq c_\delta e_\psi(0)e^{-2\kappa t}.
\end{align*}
$c_\delta$ can now be optimized for $\delta>0$, completing the proof.
\end{Proof}\\
\tb{Remark:} 
In contrast to the standard entropy method for symmetric Fokker-Planck equations \cite{ArMaToUn01}, the decay estimates \eqref{entropy-decay-nondeg} and \eqref{entropy-decay-deg} have leading multiplicative constants $c,\,c_\eps>1$. This is typical for non-symmetric Fokker-Planck equations, and it is due to the non-orthogonality of the eigenfunctions of $L$ (cf.\ \S\ref{sec:spectrum} below; and \cite{Ni14} for a closely related discussion of $L^2$-estimates for semigroups). 
Due to the applied regularisation and the above proof, we expect that the multiplicative constants in \eqref{entropy-decay-nondeg} and \eqref{entropy-decay-deg} are not sharp.


\section{Spectral analysis and flow-invariant manifolds}
\label{sec:spectrum}

In this section we shall characterise the spectrum of $L$ in $L^2$ and the corresponding eigenspaces, which are of course flow-invariant. Moreover, we also find flow-invariant manifolds that consist of Gaussian functions. In Section \ref{sec:sharprate} we shall need these two types of manifolds to prove the sharpness of decay rates for the quadratic and logarithmic entropy, respectively. 

The main difficulty in the spectral analysis of $L$ is the fact that the eigenfunctions of $L$ are not orthogonal, in contrast to the symmetric, fully parabolic case. They do, however, generate finite dimensional, $L$--invariant and mutually orthogonal subspaces of $L^2$. And this fact will be a crucial ingredient for the computation of the spectrum, cf.\ the proof of Theorem \ref{SpectrumL} below (for a closely related situation see also \cite{GaMiS13}, \cite{ArLaMaW14}).\\

First we introduce some notation.
Let $\mathcal P(\R^d)$ denote the polynomials over $\R^d$ (with complex coefficients) and let $\mathcal Q:=\mathcal P(\R^d)f_\infty$. $\mathcal Q$ is dense in $L^2(\R^d,f_\infty^{-1})$, and it is the natural space for eigenfunctions of the Fokker-Planck operator (see for example \cite{HeNiFP05} or \cite{RiFP89}). 

Let $\alpha\in\N_0^d$ be a multi-index. We write $|\alpha|=\sum\limits_{j=1}^d \alpha_j$, $\nabla^\alpha:=\sum\limits_{j=1}^d \partial_j^{\alpha_j}$. We also introduce the notation $\alpha_{l-}$ and $\alpha_{l+}$:
\begin{align*}
 (\alpha_{l+})_j := \alpha_j\ (j\neq l),&\quad (\alpha_{l+})_l := \alpha_l+1,\\
 (\alpha_{l-})_j := \alpha_j\ (j\neq l),&\quad (\alpha_{l-})_l := \alpha_l-1\quad \text{  if  } \alpha_l\ge1,\\
 \alpha_{l-} &:=0\in \N_0^d \quad \text{  if  } \alpha_l=0.
\end{align*}
So $\alpha_{l-}$, $\alpha_{l+}$ denote the multi-indices that one obtains by lowering or raising the $l$-th entry of $\alpha$ by 1. Analogously we define iterated vector shifts like, e.g., $(\alpha_{l-})_{m-}$ .

To establish the orthogonal decomposition of $L^2$, we introduce a change of coordinates. Let
\begin{align*}
 y&:=\sqrt K^{\,-1}x,\\
 g_0(y)&:=f_\infty(\sqrt K y)=c_K\exp(-\frac{|y|^2}2).
\end{align*}
Now let
\begin{align*}
g_\alpha(y)&:=\nabla_y^\alpha \,g_0(y), \quad \alpha\in\N_0^d,\\
\tilde V_m &:=\operatorname{span}\big\{g_\alpha\big|\ |\alpha|= m\big\}\subset \tQ:=\P(\R^d)g_0,\quad m\in\N_0\,.
\end{align*}
Note that the polynomial part of $g_\alpha$ has degree $|\alpha|$ and its (unique) leading monomial is $(-1)^{|\alpha|}y^\alpha$. From \cite{HeNiFP05}, \cite{RiFP89} we know that $\{g_\alpha\}_{\alpha\in\N_0^d}$ forms an orthogonal basis of $L^2(\R^d,g_0^{-1})$. Hence, the subspaces $\tilde V_m$ are also mutually orthogonal.
With the inverse coordinate transformation we see that the subspaces 
$$
  V_m :=\operatorname{span}\big\{f_\alpha\big|\ |\alpha|= m\big\}\subset \mathcal Q\,;\;\;m\in\N_0\,,\quad
\mbox{with } \;f_\alpha(x):=g_\alpha(\sqrt K^{\,-1}x)
$$
are mutually orthogonal in $L^2$. With this discussion we already obtain the first part of 
\begin{Proposition}\label{Lsplit}
 \begin{itemize}
  \item[(i)] $L^2$ has a decomposition in mutually orthogonal subspaces:
  $$
    L^2= \bigoplus_{m\in\N_0}V_m\,.
  $$
    \item[(ii)]
For every $m\in\N_0$, $V_m$ is invariant under $L$, its adjoint $L^*$ (w.r.t.\ $L^2$), and hence the semigroup $e^{tL}$; $t\ge0$.
 \end{itemize}
\end{Proposition}
\begin{Proof}
(ii) Using the above transformation we shall actually prove the equivalent invariance of the subspaces $\tilde V_m$.
Acting on the transformed function $g(y):=f(\sqrt Ky)\in L^2(\R^d,g_0^{-1})$, $L$ has the form
\begin{align*}
 \tL g &:= \div_y[ (\tilde D + \tilde R) (\nabla_y g + yg) ],\\
 \tilde D &:= \sqrt K^{\,-1}D\sqrt K^{\,-1},\\
 \tilde R &:= \sqrt K^{\,-1}R\sqrt K^{\,-1}.
\end{align*}
Note that the following properties of $D$, $R$, and $C$ also hold for the transformed matrices
(with $\tilde C := \sqrt K^{\,-1}C\sqrt K^{\,-1}$):
\begin{align*}
 2\tilde D = \tilde C K + K \tilde C^T,\quad& \tilde R^T = - \tilde R.
\end{align*}
The adjoint of $\tL$ has the form
\begin{align*}
 \tL^* g &= \div[ (\tilde D - \tilde R) (\nabla g + yg) ].
\end{align*}
Now we compute for some $l\in\{1,...,d\}$:
\begin{align*}
 \partial_{y_l} g_\alpha(y) &= \nabla^\alpha \partial_{y_l} g_0(y) = -\nabla^\alpha (y_l g_0(y)) \\
			&= -\alpha_l g_{\alpha_{l-}}(y) - y_lg_\alpha(y).
\end{align*}
So we have, writing $h_\alpha:=(\alpha_l g_{\alpha_{l-}}(y))_{l=1,\dots,d}$,
\begin{align*}
 \nabla g_\alpha(y) &= -h_\alpha(y) - yg_\alpha(y).
\end{align*}
Inserting this into $\tL$ gives
\begin{align*}
\tL g_\alpha &= \div[ (\tilde D + \tilde R) (-h_\alpha(y) -yg_\alpha(y) + yg_\alpha(y)) ]\\
					&= -\div(\tilde D h_\alpha(y)) - \div(\tilde R h_\alpha(y)),\\
\tL^* g_\alpha &= -\div(\tilde D h_\alpha(y)) + \div(\tilde R h_\alpha(y)).
\end{align*}
Further we compute
\begin{align*}
 \div(\tilde Dh_\alpha) &= \sum\limits_{j,l=1}^d \partial_{y_j} \left(\tilde D_{jl} \alpha_l g_{\alpha_{l-}}\right)(y) 
  = \sum\limits_{j,l=1}^d \alpha_l \tilde D_{jl} g_{(\alpha_{l-})_{j+}}(y),\\
 \div(\tilde Rh_\alpha) &= \sum\limits_{j,l=1}^d \alpha_l \tilde R_{jl} g_{(\alpha_{l-})_{j+}}(y)\,.\\
\end{align*}
Thus we obtain, using $R=\fract12 (CK-KC^T)$ and $D=\fract12(CK+KC^T)$,
\begin{align*}
 \tL g_\alpha &= -\sum\limits_{j,l=1}^d \alpha_l (\tilde D + \tilde R)_{jl} g_{(\alpha_{l-})_{j+}}(y) \\
					 &= - \sum\limits_{j,l=1}^d \alpha_l (\sqrt K^{\,-1}C\sqrt K)_{jl} g_{(\alpha_{l-})_{j+}}(y),\\
\tL^* g_\alpha &= - \sum\limits_{j,l=1}^d \alpha_l (\sqrt KC^T\sqrt K^{\,-1})_{jl} g_{(\alpha_{l-})_{j+}}(y).
\end{align*}
We see that $\tL g_\alpha$, $\tL^* g_\alpha$ are linear combinations only of terms $g_\beta$,  $\beta\in\N_0^d$, with $|\beta|=|\alpha|$. This completes the proof.
\end{Proof}\\

For non-degenerate Fokker-Planck equations this orthogonal decomposition of $L^2$ into invariant subspaces (or equivalently, the block-diagonal structure of the semigroup $e^{tL}$) is well known, cf.\ (57) in \cite{LNP13}. So, Proposition \ref{Lsplit} is its generalization to degenerate Fokker-Planck equations.\\

{}From the orthogonal decomposition of $L^2$ we immediately have
\begin{equation}\label{sigma-inclusion}
  \bigcup_{m\in\N_0} \sigma(L\big|_{V_m}) \subset \sigma(L)\,.
\end{equation}
First we note that the r.h.s.\ cannot include any additional eigenvalue. Otherwise, the orthogonal projection of a corresponding eigenvector to some $V_m$ would be non-trivial. Hence it would already be an eigenvector of $L\big|_{V_m}$. To have equality in \eqref{sigma-inclusion} we have to rule out that eigenvalues of  $L\big|_{V_m}$ accumulate. To this end we now prove the compactness of the resolvent of $L$:
\begin{Lemma}
\label{compactresolvent}
 Under condition (A), the operator $L$ has a compact resolvent on $L^2(\R^d,f_\infty^{-1})$.
\end{Lemma}
The technical proof is deferred to the appendix.\\

As an immediate consequence we have $\sigma(L)=\sigma_p(L)$. Moreover, the eigenvalues have no accumulation point, and all eigenspaces are finite dimensional. \\

For the following spectral analysis, let $\lambda_1,\dots,\lambda_d$ be the eigenvalues of $C$, counted with their algebraic multiplicity.

\begin{Theorem}
\label{SpectrumL}
 Assume condition (A). Then it holds:
 \begin{itemize}
  \item[(i)] The spectrum of $L$ in $L^2$ is given by
 \begin{align}\label{sigmaL}
  \sigma(L) &= \sigma_p(L) = \Big\{-\sum\limits_{j=1}^d \alpha_j\lambda_j\Big|\alpha=(\alpha_j)\in\N_0^d\Big\}\subset \{0\}\cup\left(\R^-\times i\R\right).
 \end{align}
  \item[(ii)] The eigenspace to $0$ is one-dimensional and spanned by $f_\infty$. 
  \item[(iii)] If $C$ is not defective, then the eigenfunctions of $L$ form a basis of $\mathcal Q$.
  \item[(iv)] If $C$ is defective, then the eigenfunctions and generalised eigenfunctions of $L$ form a basis of $\mathcal Q$.
 \end{itemize}
\end{Theorem}
\tb{Remark:} Formula \eqref{sigmaL} is well known for non-degenerate Fokker-Planck equations with linear drift (cf.\ \S1.4 of \cite{LNP13} and references therein). We show here that this formula carries over to degenerate diffusion matrices.

Moreover, the following proof shows that all eigenfunctions and generalised eigenfunctions of $L$ can be computed explicitly.\\

\noindent
\begin{Proof}[ (of Theorem \ref{SpectrumL})]\\
(i): Due to the orthogonal decomposition of $L^2$ and the $L$-invariance of $V_m$, it only remains to prove that 
\begin{equation}\label{sigma-Vm}
  \sigma(L\big|_{V_m}) = \Big\{-\sum\limits_{j=1}^d \alpha_j\lambda_j\Big||\alpha|=m\Big\}\,,
\end{equation}
for each $m\in\N$.
Eigenfunctions of $L\big|_{V_m}$ have the form $\phi(x)=q(x)f_\infty\in V_m\subset\mathcal Q$, where $q$ is a polynomial of order $m$.
Using $D+R=CK$ (see \S\ref{sec:steadystate}) we obtain:
\begin{align*}
L\phi			 &= \div(f_\infty(D+R)\nabla q) = f_\infty \div((D+R)\nabla q)  - f_\infty(x^TK^{-1}(D+R)\nabla q) \\
			 &= f_\infty\left[\div(D\nabla q) - x^TK^{-1}CK\nabla q\right].
\end{align*}
Hence, $\phi\in\mathcal Q$ is an eigenfunction of $L$ iff $\phi\in\P(\R^d)$ is an eigenfunction of $L^{\mathcal P}$ with:
\begin{align*}
L^{\mathcal P}q(x) &:= \nabla^TD\nabla q(x) - x^TK^{-1}CK\nabla q(x) = \nu q(x).
\end{align*}
Since the eigenvalues of $C$ (and thus of $Q$) may be complex, we shall consider the polynomial $q$ in  the space $\P(\C^d)$ in the sequel. Similar to Lemma \ref{Pdefinition}, we shall now use the Jordan normal form $J$ of $Q^T=K^{-1}CK$, with $A_1JA_1^{-1}=Q^T$ for some regular $A_1\in\C^{d\times d}$.\\ 
We introduce the (complex) coordinate transformation
\begin{align}\label{L-polynomial}
 y &:= A_1^Tx, \text{ with } y\in\C^d,\\
 p(y) &:= q((A_1^{-1})^Ty)=q(x) \in \P(\C^d).
\end{align}
So we obtain the following equation for the (transformed) eigenfunctions of $L^{\P}$:
\begin{align}
\label{eigenequ} 
  \tilde L^{\mathcal P}p(y) := 
  \nabla_y^TA_1^TDA_1\nabla_y p(y) - y^TJ\nabla_y p(y) = \nu p(y).
\end{align}
A basis of the polynomials (over $\C$) of degree $m$ or lower is given by the monomials $\{y^{\alpha}|\alpha\in\N_0^d,\ |\alpha|\leq m\}$. We order this basis by increasing degree, and in decreasing lexicographic order for monomials of the same degree. Next, we compute the matrix representation $M_\P$ of $\tilde L^\P$ with respect to this basis.
Let $e_l$ denote the $l$-th unit vector in $\C^d$, and $I_{def}$ be the set of all $l\in\{1,\dots,d\}$ for which $e_l$ is not an ordinary eigenvector of $J$. 
We compute
\begin{align}
\nonumber \tilde L^{\P}y^\alpha 
&= [\nabla^TA_1^TDA_1-y^TJ]\sum\limits_{l=1}^d \left(\alpha_l e_l y^{\alpha_{l-}}\right) \\
\nonumber &= \sum\limits_{l,m=1}^d \left([\alpha_m-\delta_{lm}]\alpha_l e_m^TA_1^TDA_1e_l y^{(\alpha_{l-})_{m-}}\right) - \sum\limits_{l=1}^d \left(\alpha_l\lambda_l y^\alpha\right) - \sum\limits_{l\in I_{def}} \alpha_l y^Te_{l-1} y^{\alpha_{l-}} \\
\label{LPpalpha}&= \sum\limits_{l,m=1}^d \left(d_{lm}(\alpha)y^{(\alpha_{l-})_{m-}}\right) + \nu_{\alpha}y^\alpha - \sum\limits_{l\in I_{def}} \alpha_l y^{(\alpha_{l-})_{(l-1)+}},
\end{align}
where $\nu_\alpha := -\sum\limits_{l=1}^d \alpha_l\lambda_l$ and $d_{lm}(\alpha):=[\alpha_m-\delta_{lm}]\alpha_l e_m^TA_1^TDA_1e_l$. The first term of the r.h.s.\ has degree $\max(|\alpha|-2,0)$. The second and the third term both have degree $|\alpha|$,  but the exponents of the third term come ``earlier'' in lexicographic order. Due to our ordering of the basis $\{y^\alpha|\ |\alpha|\leq m\}$, this implies that $M_\P$ is an upper triangular matrix. The entries on the diagonal are just the $\nu_\alpha$, which are hence the eigenvalues of $\tilde L^\P$ and hence of $L^\P$. 

All elements of $V_m$ (except of $0$) have a polynomial factor of order $m$. Hence \eqref{sigma-Vm} follows. \\
(ii) was already established in Theorem \ref{ssexistence}. \\
(iii, iv): This is a simple consequence of the decomposition of $L$ into its action on the finite dimensional, orthogonal subspaces $V_m$: The (generalised) eigenfunctions of $L\big|_{V_m}$ form a basis of $V_m$. 

As described before, the representation of $L\big|_{V_m}$ on the polynomial factor of such functions (cf.\ \eqref{L-polynomial}) can be transformed to an upper triangular matrix.
If $C$ is not defective, it is even diagonal since the last term of the r.h.s.\ of \eqref{LPpalpha} drops. Hence, the eigenfunctions of $L\big|_{V_m}$ already form a basis of $V_m$. If $C$ is defective, the generalised eigenfunctions of $L\big|_{V_m}$ have to be added to obtain a basis of $V_m$.
\end{Proof}

\medskip

Next we turn to the flow-invariant manifolds consisting of Gaussian functions, i.e., shifted and (anisotropically) stretched versions of the steady state $f_\infty(x)=c_K \exp\left(-\frac{x^TK^{-1}x}{2}\right)$ : 
\begin{Proposition}\label{manifolds}
Assume condition (A). 
\begin{itemize}
  \item[(i)] The manifold $\M_1:=\{f_\infty(x-v)\,|\,v\in\R^d\}$ is invariant under the semi-flow of \eqref{linmasterequ}.
  \item[(ii)] Let $f_0(x)=f_\infty(x-v_0)$ with some $v_0\in\R^d$. Then, 
  \begin{align}\label{v-decay}
    f(t,x):=f_\infty(x-v(t))\quad\mbox{with} \; v(t)=e^{-Ct}v_0
  \end{align}
  is the unique solution to \eqref{linmasterequ} with initial condition $f_0$.
  The logarithmic relative entropy $e_1(t):=e_1(f(t)|f_\infty)$ (with $\psi_1(s)=s\ln s-s+1$ in Definition \ref{admissableentropy}) then satisfies
  $$
    e_1(t)=\frac{v(t)^T K^{-1} v(t)}{2}\,,\quad t\ge0,
  $$
  and hence it decays at least like $\mathcal O\left(\|e^{-Ct}\|_2^2\right)$\,.
  \item[(iii)] The manifold $\M_2:=\{(2\pi)^{-\frac{d}{2}} (\det A)^{-\frac{1}{2}} \exp\left(-\frac{x^TA^{-1}x}{2}\right)\,|\,0<A=A^T\in\R^{d\times d}\}$ is invariant under the semi-flow of \eqref{linmasterequ}.
  \item[(iv)] Let $f_0(x)=(2\pi)^{-\frac{d}{2}} (\det A_0)^{-\frac{1}{2}} \exp\left(-\frac{x^TA_0^{-1}x}{2}\right)$ with some $0<A_0=A_0^T\in\R^{d\times d}$ (i.e., $f_0(x)=\sqrt{\frac{\det K}{\det A_0}} f_\infty(\sqrt K\sqrt{A_0}^{\,-1}x)$). Then, 
  \begin{align}
    f(t,x)&:=(2\pi)^{-\frac{d}{2}} (\det A(t))^{-\frac{1}{2}} \exp\left(-\frac{x^TA(t)^{-1}x}{2}\right)\quad\mbox{with} \label{A-decay}\\
    A(t)&=K + e^{-Ct} (A_0-K) e^{-C^Tt} \label{A-decay2}
  \end{align}
  is the unique solution to \eqref{linmasterequ} with initial condition $f_0$.
  The logarithmic relative entropy then satisfies
  $$
    e_1(t)=\frac12 \tr\left( \sqrt K^{\,-1} A(t) \sqrt K^{\,-1}\right)
    -\frac12 \tr\ln\left( \sqrt K^{\,-1} A(t) \sqrt K^{\,-1}\right)-\frac{d}{2}\,,\quad t\ge0,
  $$
  and hence it decays asymptotically (as $t\to\infty$) at least like $\mathcal O\left(\|e^{-Ct}\|_2^4\right)$\,.
  \item[(v)] The manifold $\M_3:=\{(2\pi)^{-\frac{d}{2}} (\det A)^{-\frac{1}{2}} \exp\left(-\frac{(x-v)^TA^{-1}(x-v)}{2}\right)\,|\,0<A=A^T\in\R^{d\times d},\,v\in\R^d\}$ is invariant under the semi-flow of \eqref{linmasterequ}.
\end{itemize}
\end{Proposition}

\tb{Remarks:} \begin{itemize}
\item[(i)] Proposition \ref{manifolds} also holds for non-degenerate diffusion matrices $D>0$.
\item[(ii)] At least in special cases (e.g., symmetric Fokker-Planck equations, 1D case), the special solutions \eqref{v-decay}, \eqref{A-decay} are well known (cf.\ \cite{ArMaToUn01}; Ex. 13 in \S11.4 of \cite{GKZ04}). We include them here in full generality, as we shall need the explicit formulas in \S\ref{sec:sharprate}.\\
\end{itemize}

\begin{Proof}[ (of Proposition \ref{manifolds})]\\
(i, ii): We insert $f(t,x):=f_\infty(x-v(t))$ into (\ref{linmasterequ}) and obtain
\begin{align*}
 f_t &= (x-v(t))^T K^{-1} \dot v(t)\,f(t),\\ 
 \div(f_\infty(D+R)\nabla\fu) 
			      &= \div[f(t) Cv(t)] = (\nabla f(t))\cdot Cv(t) \\
			      &= -(x-v(t))^T K^{-1} Cv(t) \,f(t),
\end{align*}
where we have used $D+R=CK$ and the symmetry of $K$. Hence, $\dot v=-Cv$ follows.\\

The logarithmic entropy satisfies:
$$
  e_1(t) = \ild \ln\frac{f(t)}{f_\infty}\,f(t) \d = \ild v(t)^TK^{-1}\left([x-v(t)]+\frac{v(t)}{2}\right)\,f(t) \d 
  = \frac{v(t)^TK^{-1}v(t)}2\,.
$$

\noindent
(iii, iv): Inserting \eqref{A-decay} into (\ref{linmasterequ}), an easy computation (using the formulas $\frac{d}{dt}\,\det A=(\det A)\,\tr(A^{-1}\dot A)$, $\frac{d}{dt}\,A^{-1}=-A^{-1}\dot A A^{-1}$) yields:
\begin{align}\label{A-evol}
  -\frac12\tr(A^{-1}\dot A)+\frac{x^TA^{-1}\dot A A^{-1}x}{2}=
  -\tr(D A^{-1})+x^TA^{-1}D A^{-1}x-x^TC^TA^{-1}x+\tr C\quad\forall\,x\in\R^d\,.
\end{align}
The $x$-dependent part of this equation yields
$$
  A^{-1}\dot A A^{-1}=2A^{-1}D A^{-1}-2\left(C^TA^{-1}\right)_s\,,
$$
where $(\cdot)_s$ denotes the symmetric part of the matrix. Hence,
\begin{align}\label{A-ODE}
  \dot A = 2D-CA-AC^T\,,
\end{align}
and subtracting the Lyapunov equation \eqref{Dequ} yields the evolution equation for $A(t)$:
$$
  \frac{d}{dt}\,(A-K) = -C(A-K)- (A-K)C^T\,.
$$
Multiplying 
\eqref{A-ODE} with $A^{-1}$ and taking traces, shows that also the $x$-independent part of \eqref{A-evol} is commensurate with \eqref{A-ODE}.\\

For $f(t)$ from \eqref{A-decay}, the logarithmic entropy satisfies:
\begin{align*}
  e_1(t) &= \frac12 (2\pi)^{-\frac{d}{2}}(\det A(t))^{-\frac{1}{2}}  \ild  \exp\left(-\frac{x^TA^{-1}(t)\,x}{2}\right)\,\left[x^T(K^{-1}-A^{-1}(t))x
   + \ln\det K-\ln\det A(t)\right]\d\\  
  &= -\frac12 \ln\det\left(\sqrt K^{\,-1} A(t)\sqrt K^{\,-1}\right) 
  +\frac12 (2\pi)^{-\frac{d}{2}}
  \ild  \exp\left(-\frac{|y|^2}{2}\right)\,y^T(\sqrt{A(t)} K^{-1}\sqrt{A(t)} -\Id)y\d[y]\\
  &=-\frac12\tr\ln\left(\sqrt K^{\,-1} A(t) \sqrt K^{\,-1}\right)
  +\frac12 \tr\left(\sqrt K^{\,-1} A(t) \sqrt K^{\,-1}\right)-\frac{d}{2}\,,
\end{align*}
where we used the coordinate transformation $x=\sqrt A y$. Using the expansion $\ln s-s+1\approx-(s-1)^2/2$ and (from \eqref{A-decay2})
$$
  \sqrt K^{\,-1} A(t) \sqrt K^{\,-1} = \Id + \mathcal O\left(\|e^{-Ct}\|_2^2\right)\,,
$$
we obtain the claimed decay of $e_1(t)$.

\noindent
(v): Since the evolutions of $v(t)$ and $A(t)$ turn out to be independent within $\M_3$, this result follows just as for (i) and (iii).
\end{Proof}


\section{Sharpness of the decay rate}
\label{sec:sharprate}

In this section, we investigate the sharpness of the decay rate obtained in Theorem \ref{convergencerate} under condition (A). In particular, we show that the rate is optimal for both the quadratic entropy $e_2$ and the logarithmic entropy $e_1$. As shown in \cite{ArMaToUn01}, all admissible entropies are bounded below by the logarithmic entropy and above by the quadratic one. Thus, the rate we obtained is optimal for all admissible entropies.

\begin{Theorem}
\label{sharpdecayrate}
 Let $\mu:=\min\{\Re\{\lambda\}|\lambda\in\sigma(C)\}$, where $\sigma(C)$ denotes the spectrum of $C$.
 \begin{itemize}
  \item[(i)] If $\mu$ is a (real) eigenvalue of $C$, then there exist initial conditions $f_0$, $g_0$ (different from $f_\infty$) such that the corresponding solutions $f(t)$, $g(t)$ of (\ref{linmasterequ}) satisfy
  \begin{align*}
   e_1(f(t)) = e^{-2\mu t}e_1(f_0),&\quad e_2(g(t)) = e^{-2\mu t}e_2(g_0),\quad t\geq 0.
  \end{align*}
  \item[(ii)] If $C$ has a complex conjugate eigenvalue pair with $\Re\{\lambda_{1,2}\}=\mu$, then there are initial conditions $f_0$, $g_0$ (different from $f_\infty$) such that the corresponding solutions $f(t)$, $g(t)$ of (\ref{linmasterequ}) satisfy
  \begin{align}
  \label{complexeigenvaluesharpdecay} e_1(f(t)) \leq ce^{-2\mu t}e_1(f_0),&\quad e_2(g(t)) \leq ce^{-2\mu t}e_2(g_0),\quad t\geq 0,
  \end{align}
  with some $c\ge1$, and equality holds for $t=t_0+n\tau$, $t_0\geq 0$, $\tau>0$, $n\in\N_0$. So the right hand sides of \eqref{complexeigenvaluesharpdecay} are the sharp exponential envelope functions for the entropy decay.
  \item[(iii)] If $C$ has a defective eigenvalue $\lambda$ with $\Re\{\lambda\}=\mu$, then there are initial conditions $f_0$, $g_0$ (different from $f_\infty$) such that the corresponding solutions $f(t)$, $g(t)$ of (\ref{linmasterequ}) satisfy
  \begin{align}\label{defectiveeigenvaluesharpdecay}
   e_1(f(t)) = c_0e^{-2\mu t}(e_1(f_0)+\frac{c_1}2t+\frac{c_2}2t^2),&\quad e_2(g(t)) = c_0e^{-2\mu t}(e_2(g_0)+c_1t+c_2t^2),\quad t\geq 0
  \end{align}
  for some $c_0,c_2>0$, $c_1\in\R$.
 \end{itemize}
 In all cases, $f_0$ is $\psi_1$-compatible and $g_0$ is $\psi_2$-compatible.
\end{Theorem}

\tb{Remark:} In the defective case (iii), the right hand sides of \eqref{defectiveeigenvaluesharpdecay} can also be of the form $e^{-2\mu t}(e_1(f_0)+P_{2n}(t))$ or $e^{-2\mu t}(e_2(g_0)+P_{2n}(t))$ (where $P_{2n}$ is some polynomial of degree $2n$), if $\lambda$ corresponds to a Jordan block of size $n+1$. In all of these cases the exponential decay rate is indeed reduced to $2(\mu-\eps)$ for an arbitrarily small $\eps>0$, as announced in Theorem \ref{convergencerate}. But this estimate will never be sharp.\\

The proof of Theorem \ref{sharpdecayrate} is based on special solutions of (\ref{linmasterequ}), and it is inspired by Theorem 3.11 in \cite{ArMaToUn01}. There, the sharpness of the convex Sobolev inequality (\ref{convSobinequ}) was discussed.
For the optimal decay of the logarithmic entropy we shall consider here shifted Gaussians, whose evolution was already computed in Proposition \ref{manifolds}(ii). For the quadratic entropy we shall consider a second family that consists of trajectories in $f_\infty+V_1$ (defined in \S\ref{sec:spectrum}). Their evolution is computed in the next lemma. \\

\begin{Lemma}
\label{explicitquadentropydecay}
  Let $v_0\in\R^d$. Then
 \begin{itemize}
 \item[(i)]  \begin{align*}
  g_0(x) &:=(1+x^T K^{-1} v_0)f_\infty
 \end{align*}
 is in $L^1(\R^d)$ with $\ild g_0\d=1$. Furthermore, $g_0$ is $\psi$-compatible for the quadratic entropy with $\psi_2(s)=(s-1)^2$.
 \item[(ii)] 
 The function
 \begin{align*}
  g(t,x) &:= (1+x^TK^{-1}v(t))f_\infty\quad\mbox{with } v(t)=e^{-Ct}v_0
 \end{align*}
 is the unique solution to (\ref{linmasterequ}) with initial condition $g_0$.
 \item[(iii)] The quadratic relative entropy $e_2(t):=e_2(g(t)|f_\infty)$ 
 satisfies
 \begin{align*}
 e_2(t) &= v(t)^TK^{-1}v(t)\,,\quad t\ge0.
 \end{align*}
 \end{itemize}
\end{Lemma}
\begin{Proof} \ First, we note that $g_0\geq0$ does not hold here. But this is not a problem, since we don't need positivity of the solution to define the quadratic entropy.\\

 (i): Since $f_\infty(x)=f_\infty(-x)$, we have
 \begin{align*}
  \ild v_0^TK^{-1}xf_\infty \d = 0,
 \end{align*}
 and $\ild g_0 \d = 1$ follows from the normalisation of $f_\infty$. We recall from (\ref{quadraticw}) that for quadratic $\psi$,
 \begin{align*}
  w &= \sqrt{2}(\frac{f_0}{f_\infty}-1) = \sqrt{2}v_0^TK^{-1}x\,.
 \end{align*}
Then
  \begin{align*}
   \nabla w &= \sqrt{2}K^{-1}v_0 \in L^2(\R^d,f_\infty),
  \end{align*}
 and thus $g_0$ is $\psi$-compatible for quadratic $\psi$ by Definition \ref{wdef}.\\
 
 (ii): We insert $g(t,x)$ into (\ref{linmasterequ}) and obtain
 \begin{align*}
  g_t(t,x) &= x^TK^{-1} v(t) f_\infty,\\
  \div(f_\infty (D+R)\nabla \frac{g(t,x)}{f_\infty}) &= \div(f_\infty Cv(t)) = x^TK^{-1}Cv(t)f_\infty,
 \end{align*}
 where we again used $D+R=CK$.\\
  
 (iii): The quadratic entropy satisfies
 \begin{align*}
  e_2(g(t)) &=  \ild  (\frac{g(t,x)}{f_\infty}-1)^2f_\infty \d =  \ild (x^TK^{-1}v(t))^2f_\infty \d.
 \end{align*}
 For fixed $t\geq0$, the directional derivative of $f_\infty$ satisfies
 \begin{align*}
 \partial_{v(t)} f_\infty &= -v(t)^TK^{-1}xf_\infty = -x^TK^{-1}v(t)f_\infty.
 \end{align*}
 Hence, it follows that
 \begin{align*}
  e_2(g(t)) &= - \ild (x^TK^{-1}v(t))\partial_{v(t)}f_\infty \d = \ild f_\infty\partial_{v(t)}(x^TK^{-1}v(t)) \d = v(t)^TK^{-1}v(t) \ild f_\infty \d \\
	    &= v(t)^TK^{-1}v(t).
 \end{align*}
\end{Proof}

{}From Proposition \ref{manifolds}(ii) and Lemma \ref{explicitquadentropydecay}, we see that we can reduce the discussion of sharp decay rates for relative entropies to discussing the decay of the term $v(t)^TK^{-1}v(t)$, where 
\begin{align}
\label{vequation} \dot v(t) &= -Cv(t),\quad v(t=0)=v_0\in\R^d.
\end{align}
A direct consequence is

\begin{Corollary}
\label{ZeroTangentStart}
 Let condition (A) hold, and let $t^*\in\R^+_0$. Then there is an initial condition $f_0$ [$g_0$] distinct from $f_\infty$ such that for the solution $f(t)$ [$g(t)$] to (\ref{linmasterequ}), the entropy dissipation $I_{\psi_1}$ [$I_{\psi_2}$] (see (\ref{entropyderivative})) for the logarithmic [quadratic] entropy vanishes at $t^*$, i.e. $I_{\psi_1}(f(t^*))=0$ [$I_{\psi_2}(g(t^*))=0$].
\end{Corollary}
\begin{Proof}
We take the time derivative of $v(t)^TK^{-1}v(t)$, where $v$ fulfils (\ref{vequation}), and obtain
\begin{align*}
 \ddt \left[v(t)^TK^{-1}v(t)\right] &= -v(t)^TC^TK^{-1}v(t)-v(t)^TK^{-1}Cv(t) = -2v(t)^TK^{-1}DK^{-1}v(t),
\end{align*}
where we have used (\ref{Dequ}). Let $0\neq w\in\ker D$. Setting $v_0:=e^{Ct^*}Kw$ implies $v(t^*)=Kw$, and hence:
$$
  \frac{d}{dt} \left(v^TK^{-1}v\right)\Big|_{t=t^*} = 0\,.
$$
This completes the proof.
\end{Proof}

We will now use Proposition \ref{manifolds}(ii) and Lemma \ref{explicitquadentropydecay} to prove Theorem \ref{sharpdecayrate}.\\

\begin{Proof}[ (of Theorem \ref{sharpdecayrate})]\\

 (i): There exists $0\neq v_0\in\R^d$ with $Cv_0=\mu v_0$. So the solution of (\ref{vequation}) is $v(t)=e^{-\mu t}v_0$, and thus
 \begin{align*}
  2e_1(f(t))=e_2(g(t))=v(t)^TK^{-1}v(t) &= e^{-2\mu t}v_0^TK^{-1}v_0\,.
 \end{align*}\\

 (ii): There exists $0\neq w\in\C^d$ with $Cw=\lambda w$, $\lambda\in\C$, $\Re\{\lambda\}=\mu>0$, $\Im\{\lambda\}=\omega\neq0$. Then $\overline{w}$ fulfils $C\overline{w}=\overline{\lambda}\overline{w}$, since $C$ is real. Moreover $v_0:=w+\overline{w}\in\R^d$, and $v_1:=i(\overline{w}-w)\in\R^d$. One easily verifies that $v(t):=e^{-\mu t}\left(\cos(\omega t)v_0+\sin(\omega t)v_1\right)$
 is the solution to (\ref{vequation}). We define
 \begin{align*}
  c:=\sup_{t\in\R^+_0} \left(\cos(\omega t)v_0+\sin(\omega t)v_1\right)^TK^{-1}\left(\cos(\omega t)v_0+\sin(\omega t)v_1\right) > 0,
 \end{align*}
 since $K$ is positive definite. Since $v(t)$ is $\frac{2\pi}{\omega}$-periodic, the function $v(t)^TK^{-1}v(t)$ takes the value $c$ for $t=t_0+k\frac{\pi}{\omega}$, with some $t_0\in\R^+_0$. It follows that
 \begin{align*}
  v(t)^TK^{-1}v(t) &= e^{-2\mu t}\left(\cos(\omega t)v_0+\sin(\omega t)v_1\right)^TK^{-1}\left(\cos(\omega t)v_0+\sin(\omega t)v_1\right) \leq ce^{-2\mu t},
 \end{align*}
 with equality for $t=t_0+k\frac{\pi}{\omega}$.\\
 
 (iii): We confine ourselves here to the case $\lambda=\mu\in\R$; the general case can be obtained by an extension of (ii). So, let $w,h\in\R^d$ with $Cw=\mu w$, $Ch=\mu h +w$. Let $v_0:=h$, then $v(t):=e^{-\mu t}(h-tw)$ is the solution to (\ref{vequation}), and
 \begin{align*}
  v(t)^TK^{-1}v(t) &= e^{-2\mu t}(h-tw)^TK^{-1}(h-tw)=e^{-2\mu t}(v_0^TK^{-1}v_0+c_1t+c_2t^2).
 \end{align*}
\end{Proof}

{}From the proof of Theorem \ref{sharpdecayrate}, we see that the constant $c$ in $e_\psi(f(t))\leq ce^{-2\mu t}$ does not derive from the initial entropy in a straightforward way, unless all eigenvalues of $C$ are real and non-defective. For case (ii), if $|v_1|\gg|v_0|$, then $c$ can be very large in comparison to $e_\psi(f_0)$; for case (iii), the same holds for $|w|\gg|h|$. \\

Next we shall discuss the sharpness of the leading multiplicative constant $c>1$ in the decay estimate of Theorem \ref{entropydecay} (for the non-defective case). The quest for these sharp constants for non-symmetric semigroups (particularly in $L^2$-estimates) is an active research area (cf.\ \cite{Ni14}).

Next we shall establish that, for any (admissible) choice of the matrix $P$, the leading constant in the entropy decay estimate \eqref{entropydecay-nondeg} is sharp in 2D. This also holds for regular diffusion matrices, as discussed in \S\ref{sec:nonsymmFP}. But in higher dimensions it does not hold in general.

\begin{Proposition}\label{sharp-constant}
 Let $d=2$ and let $L$ be non-symmetric on $L^2$, i.e. $L_{as}\neq0$ (cf. Theorem \ref{OPsplit}). Further assume that $C$ is not defective. Then for any matrix $P$ chosen according to Lemma \ref{Pdefinition} and for quadratic or logarithmic $\psi$, there exist initial data $f_0$ such that the estimate
 \begin{align}\label{e-decay2}
  e_\psi(f(t)) &\leq \frac{S_\psi(f_0)}{2\lambda_P} e^{-2\mu t}\,,\quad t\ge0
 \end{align}
 is optimal both with respect to the rate and the multiplicative constant.
\end{Proposition}
\begin{Proof}
The idea of the proof is to find an initial condition $f_0$ such that \eqref{e-decay2} is an equality at $t=0$. Hence, $f_0$ has to be chosen as an ``optimal function'' for the convex Sobolev inequality \eqref{convSobinequ}. But at the same time the trajectory $f(t)$ has to prove that $2\mu$ with $\mu=\min\{\Re\{\lambda\}\,|\,\lambda\in\sigma(C)\}$ is the sharp decay rate. 

Here we only give the proof for the logarithmic entropy, as the case of the quadratic entropy is very similar. The first requirement (sharp constant at $t=0$) holds iff $f_0$ is a shifted Gaussian of the form $f_0(x)=f_\infty(x-v_0)$, where $v_0\ne0$ satisfies the eigenvalue equation $PK^{-1} v_0=\lambda_P v_0$ (cf.\ \S3.5 in \cite{ArMaToUn01} and Remark \ref{const-remark}(ii)). For such an initial condition, Proposition \ref{manifolds}(ii) shows that 
 \begin{align*}
 e_1(f(t)) &= \frac{v(t)^TK^{-1}v(t)}2\,,\quad \mbox{with } v(t)=e^{-Ct}v_0\,.
 \end{align*}
With this explicit representation, it remains to show that $e_1(f(t))$ does not decay faster than $c\,e^{-2\mu t}$ with some $c>0$. 

Since we assumed that $C$ is non-defective, we have to discuss two cases: If $C$ has a complex conjugate eigenvalue pair (with real part $\mu$), $e^{-Ct}v$ decays for all $v\ne0$ exactly with rate $\mu$. And this proves the optimality statement.

It remains to discuss the case where $C$ has two different real eigenvalues, $0<\mu<\mu_2$. Here the decay rate is sharp iff $v_0$ is not an eigenvector of $C$ to the eigenvalue $\mu_2$ (as we would have $v(t)=e^{-\mu_2t} v_0$ otherwise). Equivalently, we want to rule out that $\tilde v_0:=\sqrt{K}^{\,-1}v_0$ is not an eigenvalue of $\tC:=\sqrt{K}^{\,-1}C\sqrt{K}$ pertaining to $\mu_2$.

The matrix $\tC$ can be diagonalised over $\R$: $\tC = A\cC A^{-1}$ for some $A\in\R^{2\times 2}$ and $\cC = \diag(\mu,\mu_2)$.
Inequality (\ref{matrixestimate}) then becomes
\begin{align*}
 \cC\cP + \cP\cC & \geq 2\mu\cP,
\end{align*}
where $\cP:=A^T\tP A$ is symmetric and positive definite, and $\tP:=\sqrt{K}^{\,-1}P \sqrt{K}^{\,-1}$. A short computation shows that this inequality can only hold if $\cP$ is diagonal. We write
\begin{align*}
 A &= \left(\begin{array}{cc} a & c \\ b & d \end{array}\right),
\end{align*}
where $w_1:=(a,b)^T$ and $w_2:=(c,d)^T$ are the eigenvectors of $\tC$ to $\mu$ and $\mu_2$, respectively. Assume now that $w_2$ is an eigenvector of $\tP$ pertaining to $\lambda_P$ (just as $\tilde v_0$ is). Using $(A^T)^{-1}\cP=\tP A$ we compute
\begin{align*}
 \tP w_2 &= \frac{\cP_{22}}{\det A} \left(\begin{array}{c} -b \\ a \end{array}\right).
\end{align*}
Then the assumption $\tP w_2=\lambda_P w_2$ implies that $w_1\perp w_2$. Thus $\tC$ is symmetric, i.e. (after multiplying the equality $\tC=\tC^T$ by $\sqrt{K}$ from left and right) $CK=KC^T$. Then \eqref{Dequ} implies $D=CK$. 

If $D$ is not regular, this is a contradiction. If $D$ is regular (as in \S\ref{sec:nonsymmFP}), then $C=DK^{-1}$ and thus $L$ can be written as
\begin{align*}
 Lf &= \div(D[\nabla f + K^{-1}xf]) = \div(D[\nabla f + f\,\nabla\frac{x^TK^{-1}x}{2}])\,,
\end{align*}
with a symmetric, positive definite $K^{-1}$. But this is a symmetric Fokker-Planck equation with $L_{as}=0$, which again contradicts our assumptions.
\end{Proof}

\tb{Remark:} For the case $d=3$, there are counterexamples to this result: For certain choices of $P$ one cannot have both a sharp rate and a sharp constant. This is the case in the example
\begin{align*}
 D=\left(\begin{array}{ccc} 1 & 0 & 0  \\ 0 & 1 & 0  \\ 0 & 0 & 0  \end{array}\right), &\quad C=\left(\begin{array}{ccc} 1 & 0 & 0  \\ 0 & 2 & 0  \\ 0 & 1 & 3  \end{array}\right), \quad K=\left(\begin{array}{ccc} 1 & 0 & 0  \\ 0 & 0.5 & -0.1  \\ 0 & -0.1 & 1/30  \end{array}\right), \quad P:=\left(\begin{array}{ccc} 2 & 0 & 0  \\ 0 & 61 & -11  \\ 0 & -11 & 2  \end{array}\right),
\end{align*}
with $\mu=1$. However, sharpness holds for ``better'' choices of $P$ (e.g., with the modification $P_{1\,1}=1$).

It remains an open question, whether one can always choose $P$ ``sufficiently careful'' such that both rate and constant are sharp.\\


\section{Kinetic Fokker-Planck equation}
\label{sec:kinFP}

In this section we shall illustrate how the modified entropy method from \S\ref{sec:modentmod} can be extended to kinetic Fokker-Planck equations \eqref{kinFP} with non-quadratic potentials (i.e.\ a drift term that is nonlinear in the position variable). Several proofs of the entropy and $L^2$--decay of this equation have already been obtained in the last few years: In \cite{DeVi01}, algebraic decay was proved for potentials that are asymptotically quadratic (as $|x|\to\infty$) and for initial conditions that are bounded below and above by Gaussians. The authors used logarithmic Sobolev inequalities and entropy methods. In \cite{HeNi04}, exponential decay was obtained also for faster growing potentials and more general initial conditions. That proof is based on hypoellipticity techniques. In \S2 of \cite{BaB13}, exponential convergence is proved with a modified $\Gamma_2$--approach for potentials with a bounded Hessian. In \cite{DoMoScH10} exponential decay in $L^2$ was proved, allowing for potentials with 
linear or super-linear growth.
This section will now provide an alternative proof of exponential entropy decay for \eqref{kinFP} with a certain class of non-quadratic potentials and for all admissible relative entropies $e_\psi$.

It is well known \cite{Vi02} that the unique normalized steady state of \eqref{kinFP} is given by
\begin{equation}\label{kinFP-steady}
  f_\infty(x,v)=\exp\left\{-\frac{\nu}{\sigma}[V(x)+\frac{|v|^2}{2}]\right\}\,,\quad x,v\in\R^n\,.
\end{equation}
Here we consider \eqref{kinFP} with $\displaystyle \lim_{|x|\to\infty} V(x)=\infty$ such that $f_\infty\in L^1(\R^{2n})$. For well-posedness and instantaneous smoothing results of the kinetic Fokker-Planck equation \eqref{kinFP} we refer to \cite{HeNi04,DeVi01,Vi02} as well as \S A.20, A.21 of \cite{ViH06}.

First we rewrite \eqref{kinFP} in the form of \eqref{masterequ}:
\begin{equation}\label{kinFP2}
   \partial_t f = Lf := \div_\xi [D\nabla_\xi f + G(\xi) f],
\end{equation}
with the notation $\xi:=(x,\,v)^T\in\R^d,\,d=2n$, the block diagonal diffusion matrix $D=\left(\begin{array}{cc} 0 & 0 \\ 0 & \sigma\,\Id \end{array}\right)$, and the drift vector field $G(x,v)=\left(\begin{array}{c} -v \\ \nabla_x V+\nu v  \end{array}\right)$. Moreover, we shall use the abbreviation $E(\xi):=\frac{\nu}{\sigma}[V(x)+\frac{|v|^2}{2}]$.\\

Concerning the positivity of the solution, we shall discuss here only the 1D case (i.e.\ $x,\,v\in\R$; $d=2$), using the interior maximum principle as in \S\ref{sec:solpos}:
\begin{Proposition}
\label{positivity-kinFP}
Let $V\in W^{2,\infty}_{loc}(\R)$ and $f_0\in L^1_+(\R^2)$ with $\int f_0(\xi) \d[\xi]=1$. Then the solution of \eqref{kinFP} satisfies
$$
  f(t,x,v)>0\quad \mbox{ for }t>0;\;\forall \,x,\,v\in\R\,.
$$
\end{Proposition}
\begin{Proof}
As for Theorem \ref{globalpos} we first rewrite the kinetic Fokker-Planck operator in degenerate elliptic form:
\begin{align*}
 \tilde Lf &:= \left[\left(\begin{array}{c} \partial_t \\ \nabla_\xi  \end{array}\right)^T\tilde D\left(\begin{array}{c} \partial_t \\ \nabla_\xi \end{array}\right)\right]f + b\cdot\left(\begin{array}{c} \partial_t \\ \nabla_\xi \end{array}\right)f,
\end{align*}
where
\begin{align}
 \tilde D &:= \left(\begin{array}{ccc} 0 & 0 & 0 \\ 0 & 0 & 0 \\ 0 & 0 & \sigma \end{array}\right)\in \R^{3\times3},\label{D-tilde}\\
 b(x) &:= \left(\begin{array}{c} -1 \\ -v \\ \nabla_x V+\nu v \end{array}\right) \in \R^3.\label{b-drift}
\end{align}
Comparing this with $L$, we have
\begin{align}
 \label{tildecomparison-kinFP}
 \tilde Lf &= Lf - f_t - \nu f.
\end{align}

Using the drift and diffusion trajectories of $\tilde L$ (cf.\ Definition \ref{DDtrajectories}) we find that the propagation set of each point $p=(t_1,x_1,v_1)\in\R^+\times\R^2$ contains an open layer ``before time $t_1$''. More precisely, there exists a continuous function $\tilde t:\,\R_x\to[0,t_1)$ with
\begin{equation}\label{prop-set}
  \left[\{(\tilde t(x),t_1)\times\{x\}\,|\,x\in\R\}\cup(t_1,x_1)\right]\times \R_v\subset S(p,\R^3)\,.
\end{equation}
The slightly technical proof of this statement is deferred to the Appendix.

We proceed as in the proof of Theorem \ref{globalpos}: Assuming $f(t_1,\xi_1)=0$ for some $t_1>0,\,\xi_1\in\R^2$ would imply $f=0$ on 
$\overline{S((t_1,\xi_1),\R^3)}$ and in particular $f(t_1,\cdot)\equiv0$. But this contradicts the mass conservation of  \eqref{kinFP}.
\end{Proof}\\

In analogy to Theorem \ref{OPsplit}, the operator $L$ from \eqref{kinFP2} can be decomposed on $L^2:=L^2(\R^d,f_\infty^{-1})$ in its symmetric and antisymmetric part as:
\begin{align*}
  \nonumber L_sf &= \div_\xi(D\nabla_\xi(\fract{f}{f_\infty})f_\infty), \\
  L_{as}f &= \div_\xi(R\nabla_\xi(\fract{f}{f_\infty})f_\infty),
\end{align*}
with the skew-symmetric (and $\xi$--independent!) matrix $R=\frac{\sigma}{\nu}\left(\begin{array}{cc} 0 & -\Id \\ \Id & 0 \end{array}\right) \in \R^{d\times d}$. \\

Next we introduce the modified entropy dissipation functional as in \eqref{Sdefinition}:
\begin{align*}
  S_\psi(f) &:=  \int\limits_{\fu>0} \psi''(\fu) \nabla(\fu)^TP\nabla(\fu)f_\infty\d[\xi], 
\end{align*}
with a positive definite and $\xi$--independent matrix $P\in \R^{d\times d}$ to be chosen later.
For the decay of $S_\psi(f(t))$, the computations from the proof of Proposition \ref{Sconvergence} carry over up to the following inequality:
\begin{align}\label{kinFP-Sineq}
  \ddt S_\psi(f(t)) &\le - \ild \psi''(\fract{f}{f_\infty}) u^T[(D-R)\frac{\partial^2 E}{\partial \xi^2}P+P\frac{\partial^2 E}{\partial \xi^2}(D+R)]u f_\infty \d[\xi]\,,
\end{align}
with the notation $u:=\nabla_\xi \frac{f}{f_\infty}$. In analogy to \S\ref{sec:modentmod} we define the matrix
\begin{equation}\label{matrixQ}
  Q(x):=(D-R)\frac{\partial^2 E}{\partial \xi^2}=\left(\begin{array}{cc} 0 & \Id \\ -\frac{\partial^2 V}{\partial x^2}(x) & \nu\,\Id \end{array}\right). 
\end{equation}
In order to estimate the r.h.s.\ of \eqref{kinFP-Sineq} we need to find an $x$--independent matrix $P>0$ and a constant $\kappa>0$, such that
$$
  Q(x)P+PQ^T(x)\ge2\kappa P\quad\forall\,x\in\R^n\,.
$$

In order to keep the presentation simple, we shall consider from now on only the 1D case, i.e.\ $x,\,v\in\R$ ($d=2$). More importantly, we shall consider potentials with bounded second derivatives. More precisely, we assume
\begin{equation}\label{potential}
   V(x)=\omega_0^2\,\frac{x^2}{2}+\tilde V(x)\quad\mbox{ with }\;|\tilde V''(x)|\le const.\; \forall\,x\in\R,
   \quad\mbox{ and }\;\omega_0\ne0.
\end{equation}
Corresponding to the ``unperturbed'' potential $\omega_0^2\,\frac{x^2}{2}$, we define the constant matrix 
$$
  Q_0:= \left(\begin{array}{cc} 0 & 1 \\ -\omega_0^2 & \nu \end{array}\right)\in\R^{2\times 2}\,,
$$
having the (real or complex) eigenvalues $\lambda_{1,2}=\frac{\nu}{2}\pm \sqrt{\frac{\nu^2}{4}-\omega_0^2}$.
Following the proof of Lemma \ref{Pdefinition} we choose the positive definite matrix $P$ corresponding to $Q_0$, using $b_j=1$ in \eqref{simpleP}. This choice of $b_j$ is for simplicity of the presentation only, and the final result could be optimised w.r.t.\ the quotient $b_1/b_2$. Let
\begin{equation}\label{Pdefinition1}
  P:=\left(\begin{array}{cc} 2 & \nu \\ \nu & \nu^2-2\omega_0^2 \end{array}\right)\quad \mbox{ if }4\omega_0^2<\nu^2\,,
\end{equation}
and 
\begin{equation}\label{Pdefinition2}
  P:=\left(\begin{array}{cc} 2 & \nu \\ \nu & 2\omega_0^2 \end{array}\right)\quad \mbox{ if }4\omega_0^2>\nu^2\,.
\end{equation}
Then, Lemma \ref{Pdefinition} implies
\begin{equation}\label{Q0-ineq}
  Q_0P+PQ_0^T\ge 2\kappa_0 P\,,
\end{equation}
with 
\begin{equation}\label{kappa-def}
  2\kappa_0 := \left\{
\begin{array}{c c l}
\nu-\sqrt{\nu^2-4\omega_0^2}, & & 4\omega_0^2<\nu^2\,, \\
\nu, & & 4\omega_0^2>\nu^2\,.
\end{array} \right. 
\end{equation}
We omit the defective case $4\omega_0^2=\nu^2$ here. But also in this case, a matrix $P=P(\eps)$ could easily be found from the proof of Lemma \ref{Pdefinition} (ii). 

In order to include the perturbative term 
$-\tilde V''$ from \eqref{matrixQ} we shall use the following lemma.
\begin{Lemma}
\label{Qperturbation}
Let $\lambda>0$ be fixed. Then, for any $0<P=P^T\in\R^{2\times 2}$ it holds:
\begin{equation}\label{Qperturb-inequ}
  \tilde P(\tau):=\left(\begin{array}{cc} 0 & 0 \\ \tau & 0 \end{array}\right)P + P\left(\begin{array}{cc} 0 & \tau \\ 0 & 0 \end{array}\right)
  +\lambda P\ge 0
\end{equation}
iff
\begin{equation}\label{Qperturb-cond}
  |\tau|\le \frac{\sqrt{\det P}}{p_{1\,1}}\,\lambda\,.
\end{equation}
\end{Lemma}
\begin{Proof}
By construction, $\tilde P(\tau=0)$ is positive definite. Since the eigenvalues of $\tilde P$ are continuous in $\tau\in \R$, we shall consider the zeros of $\det\tilde P(\tau)$. We have
$$
  \frac1\lambda\tilde P = \left(\begin{array}{cc} p_{1\,1} & p_{1\,2}+\tilde\tau p_{1\,1} \\ p_{1\,2}+\tilde\tau p_{1\,1} & p_{2\,2}+2\tilde\tau p_{1\,2} \end{array}\right)\,,
$$
with $\tilde\tau:=\tau/\lambda$.
Now, $\det(\tilde P/\lambda)=-{\tilde\tau}^2 p_{1\,1}^2 + (p_{1\,1}p_{2\,2}-p_{1\,2}^2)$. And this proves condition \eqref{Qperturb-cond}.
\end{Proof}\\

This allows us now to prove the exponential decay of $S_\psi(f(t))$, in analogy to Proposition \ref{Sconvergence}:
\begin{Proposition}
\label{Sconvergence-kinFP}
Let $4\omega_0^2\ne \nu^2$and let $\tilde V$ from \eqref{potential} satisfy for some fixed $\lambda\in(0,2\kappa_0)$ and $\forall\,x\in\R$:
$$
  |\tilde V''(x)| \le \frac{\sqrt{\det P}}{2}\,\lambda = \sqrt{|\omega_0^2-\nu^2/4|}\,\lambda\,
$$
for the matrix $P$ chosen in \eqref{Pdefinition1} or \eqref{Pdefinition2}. Then
\begin{align}\label{S-convergence-kinFP}
 S_\psi(f(t)) &\leq S_\psi(f_0)e^{-(2\kappa_0-\lambda)t},\quad t\geq 0,
\end{align}
with $\kappa_0$ defined in \eqref{kappa-def}.
\end{Proposition}
\begin{Proof}
{}From \eqref{Q0-ineq} and \eqref{Qperturb-inequ} with $\tau=-\tilde V''(x)$ we obtain
$$
  Q(x)P+PQ^T(x)\ge(2\kappa_0-\lambda) P\quad\forall\,x\in\R\,.
$$
Hence, \eqref{kinFP-Sineq} yields
$$
  \ddt S_\psi(f(t)) \le - (2\kappa_0-\lambda) \ild \psi''(\fract{f}{f_\infty}) u^T\,P\,u f_\infty \d[\xi]
  =- (2\kappa_0-\lambda)S_\psi(f(t))  \,,
$$
and the result follows.
\end{Proof}\\

As in \S\ref{sec:modentmod}, the decay of the modified entropy dissipation functional $S_\psi(f(t))$ implies the exponential decay of the relative entropy. But in contrast to Theorem \ref{convergencerate}, we shall refrain here from extending the regularisation Theorem \ref{entropyregularisation} to non-quadratic drift terms. 

\begin{Theorem}\label{kinFPdecay}
Let $\psi$ generate an admissible entropy and let $f$ be the solution to the kinetic Fokker-Planck equation (\ref{kinFP}) with a $\psi$-compatible initial state $f_0$ (in the sense of Definition \ref{wdefinition}). Under the assumptions of Proposition \ref{Sconvergence-kinFP} we then have:
\begin{equation}\label{e-convergence-kinFP}
  e_\psi(f(t)|f_\infty) \le c\,S_\psi(f_0)e^{-(2\kappa_0-\lambda)t},\quad t\geq 0\,,
\end{equation}
for some constant $c>0$ independent of $f_0$.
\end{Theorem} 
\begin{Proof}
For the case $4\omega_0^2>\nu^2$ we compute:
$$
  |\tilde V''(x)|\le\sqrt{\omega_0^2-\nu^2/4}\,\lambda< \sqrt{\omega_0^2-\nu^2/4}\,\nu\le\omega_0^2
  \quad\forall\,x\in\R\,,
$$
and the same estimate also holds for the case $4\omega_0^2<\nu^2$. Hence, $V$ from \eqref{potential} is uniformly convex on $\R$.

Thus, we can proceed as in the proof of Proposition \ref{Sfiniteefinite}: There exists a $\lambda_P>0$, such that the following Bakry-\'Emery condition for the operator $L_Pf:=\div_\xi(P\nabla_\xi(\fract{f}{f_\infty})f_\infty)$ holds uniformly in $x\in\R$:
$$
  \frac{\partial^2 E}{\partial\xi^2}(x) =\frac\nu\sigma \diag(V''(x),1)\ge \lambda_P P^{-1}\,.
$$
This implies the convex Sobolev inequality 
$$
  e_\psi(g|f_\infty) \leq \frac1{2\lambda_P} S_\psi(g)\,.
$$
And \eqref{e-convergence-kinFP} follows from \eqref{S-convergence-kinFP}.
\end{Proof}\\

The strategy of this section also applies to further examples of hypocoercive Fokker-Planck equations with nonlinear drift terms, see \S1.7.3 in \cite{Erb14}. E.g., this includes the following, generalized kinetic Fokker-Planck equation discussed in \cite{DeVi01}:
$$
  \partial_t f + \nabla_vW(v)\cdot\nabla_x f-\nabla_x V\cdot\nabla_v f= \nu \div_v (\nabla_vW(v)f)+\sigma\Delta_v f\,;\quad x,\,v\in\R^n;\,t>0,
$$
with $W(v)$ strictly convex and growing quadratically.


\section{Non-degenerate, non-symmetric Fokker-Planck equations}
\label{sec:nonsymmFP}

In this section we shall illustrate how the above developed method applies to non-symmetric Fokker-Planck equations that are \emph{non-degenerate}. We shall consider 
\begin{equation}\label{nsFP}
  \partial_t f=Lf:=\div (D\nabla f+ C xf)\,,
\end{equation}
with $D=D^T$ positive definite and $C$ positively stable. Its unique normalized steady state 
is still the (non-isotropic) Gaussian given in Theorem \ref{ssexistence}:
$$
  f_\infty(x)=c_K \exp(-\frac{x^TK^{-1}x}{2}) = c_K e^{-V(x)}\,,
$$
with the covariance matrix $K$ defined via \eqref{Dequ}.

With the coordinate transformation $x=\sqrt D\tilde x$ we can normalise the diffusion matrix and bring \eqref{nsFP} to the form analysed in \S\ref{sec:modentmod}:
\begin{equation}\label{nsFP-norm}
  \partial_t \tilde f=\div (\nabla \tilde f+\tilde C\tilde x\tilde f)\,,
\end{equation}
with the similarity transformation $\tilde C:={\sqrt{D}}^{\,-1}C\sqrt D$. Hence, 
\begin{equation}\label{mu-tilde}
  \mu:=\min\{\Re(\lambda)\,|\,\lambda\in\sigma(C)\}= \min\{\Re(\lambda)\,|\,\lambda\in\sigma(\tilde C)\}\,.
\end{equation}
Its steady state is
$$
  \tilde f_\infty(\tilde x)=c_{\tilde K} \exp(-\frac{\tilde x^T\tilde K^{-1}\tilde x}{2})\,,
$$
with $\tilde K=\sqrt{D}^{\,-1} K \sqrt{D}^{\,-1}$.
Clearly, the above computations of the \emph{hypocoercive entropy method} still apply without changes to the non-degenerate case. Therefore, the Theorems \ref{entropydecay}, \ref{convergencerate} and Remark \ref{const-remark} apply verbatim to the non-degenerate, non-symmetric Fokker-Planck equation \eqref{nsFP}. Here, the functionals $e_\psi$ and $S_\psi$ are expressed directly in the original variable $x$. Also, the scaling matrix $P$ from Lemma \ref{Pdefinition} is constructed directly from the original matrices $C$, $K$ (and not from $\tilde C$, $\tilde K$). Moreover, due to \eqref{mu-tilde}, the decay rate is independent of $D$!

Next we shall compare this new result to the known estimate from the standard entropy method. 
For \eqref{nsFP}, the standard entropy method from \S2.4 of \cite{ArMaToUn01} yields the decay estimate (with multiplicative constant equal to 1):
\begin{equation}\label{e-decay}
  e_\psi(f(t)|f_\infty)\le e^{-2\lambda_K t} e_\psi(f_0|f_\infty),\qquad t\ge0\,.
\end{equation}
Here, $\lambda_K$ is the largest constant to satisfy the Bakry-\'Emery condition $\frac{\partial^2 V}{\partial x^2}=K^{-1}\ge \lambda_K D^{-1}$, i.e.\ the smallest eigenvalue of $\sqrt D K^{-1}\sqrt D=\tilde K^{-1}$. For non-degenerate Fokker-Planck equations with Gaussian steady states, it is well known that this decay rate $\lambda_K$ is ``optimal'' (cf.\ \S3.5 of \cite{ArMaToUn01}, and the above sketched transformation for $D\ne\Id$). This also means that the \emph{non-symmetric entropy methods} from \cite{ArCaJuL08, BoGe10} cannot yield an improvement for this class of equations. In order to understand this ``optimality'' statement we first consider an example.\\

\begin{figure}[ht!]
 \includegraphics[scale=.41]{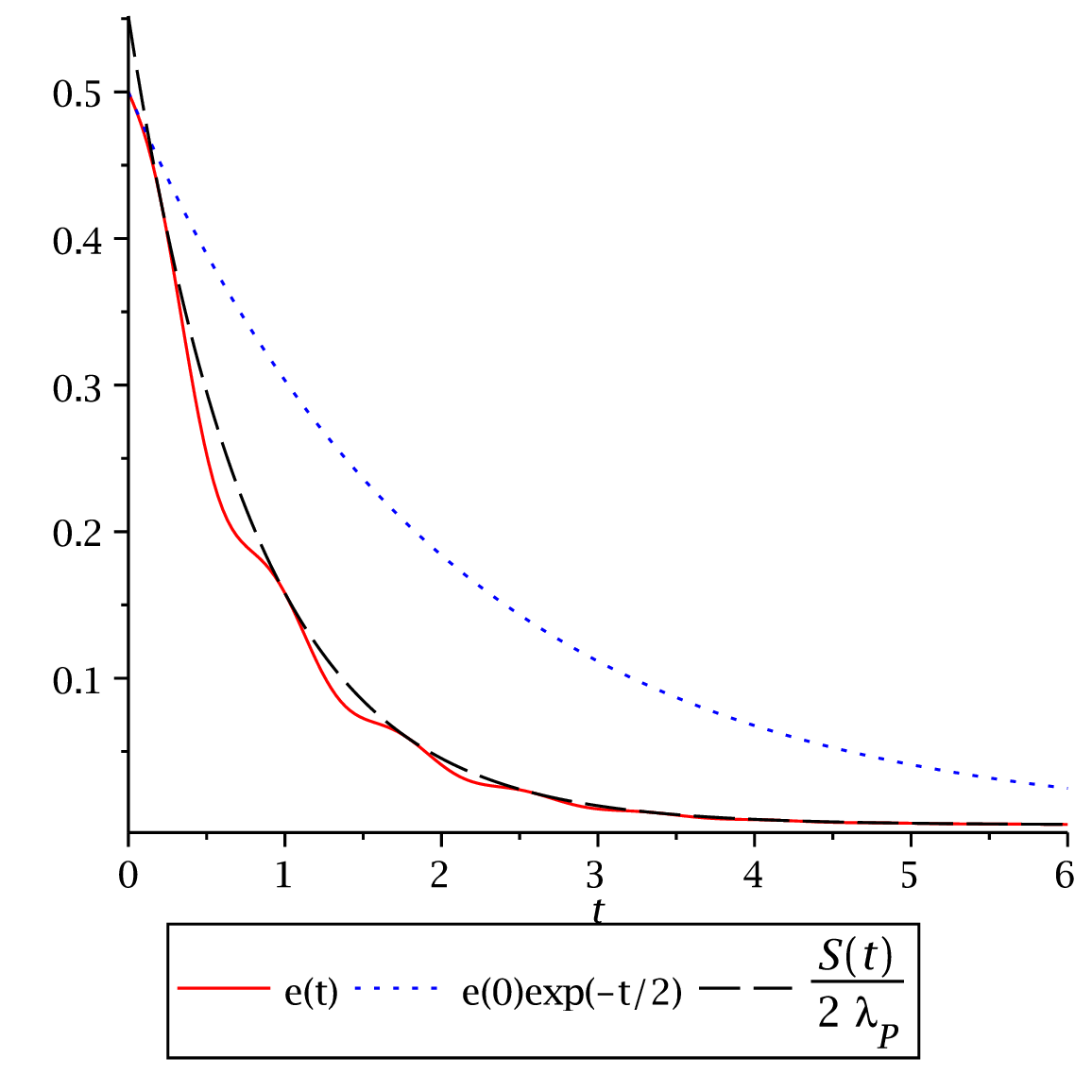}
 \caption{Entropy decay for the non-degenerate, non-symmetric Fokker-Planck equation \eqref{nsFP} 
 with $D=\diag(1/4,\;1)$, $C=[1/4 \;\, -4 \; ;\; 4 \;\, 1]$ 
 : \; \textcolor{red}{---} \, decay of the logarithmic entropy; \;\textcolor{blue}{$\cdots$} \, estimate of the \emph{local decay rate} from the standard entropy method; \;- - - \, estimate of the \emph{global decay rate} from the hypocoercive entropy method.}
 \label{fig:non-symmetric}
\end{figure}

We consider the non-symmetric Fokker-Planck equation \eqref{nsFP} 
with 
$$
  D=\diag(1/4,\;1)\,,\quad C=\left(\begin{array}{cc} 1/4 & -4 \\ 4 & 1 \end{array}\right)\,:
$$
The ``wavy'' decay of the logarithmic relative entropy (cf.\ the solid line in Fig. \ref{fig:non-symmetric}) is due to complex conjugate eigenvalue pairs of the operator $L$ and/or the non-orthogonality of its eigenfunctions in $L^2(\R^d,f_{\infty}^{-1})$ (as discussed in \S\ref{sec:spectrum}). The steady state is a Gaussian with $K^{-1}=\Id$. So the standard Bakry-\'Emery condition $\tilde K^{-1}=\diag(1/4,\,1)\ge\lambda_K \Id$ yields exponential entropy decay with the (optimal) \emph{local decay rate} $\lambda_K=1/4$ (see the dotted curve in Fig.\ \ref{fig:non-symmetric}). This reflects the (in absolute value) smallest slope of the relative entropy at any $t\ge0$. In Fig. \ref{fig:non-symmetric} this is realized, e.g., at $t=0$ with $f_0(x):=f_\infty(x-v_0)$, and $v_0=(1,\,0)^T$ is an eigenvector for the smallest eigenvalue of $\tilde K^{-1}$ (cf.\ \S3.5 of \cite{ArMaToUn01}). But the corresponding exponential function on the r.h.s.\ of \eqref{e-decay} is a crude estimate for large time.

The hypocoercive entropy method from \S\ref{sec:modentmod} yields the estimate
\begin{equation}\label{2d-decay}
  e_\psi(f(t)|f_\infty)\le \frac{1}{2\lambda_P}S_\psi(f_0)e^{-2\mu t},\qquad t\ge0\,,
\end{equation}
where $\mu=\min\{\Re(\lambda)\,|\,\lambda\in\sigma(C)\}=5/8$. The corresponding exponential function from the r.h.s.\ of \eqref{2d-decay} (see the dashed curve in Fig.\ \ref{fig:non-symmetric}) is here the sharp envelope of the relative entropy function; it accurately describes its \emph{global decay}. This was predicted in Proposition \ref{sharp-constant} for the case $d=2$. 

In contrast to the proof of Proposition \ref{sharp-constant} we did not choose here $v_0$ as an eigenvector of $PK^{-1}$. Hence this envelope does not touch the entropy function at $t=0$, but periodically at later times.
Since $C$ has a pair of complex conjugate eigenvalues, the 2D-trajectory $v(t)$ from \eqref{v-decay} converges to the origin in a spiral. Thus it will eventually be in the direction of the $\lambda_P$--eigenvector of $PK^{-1}$. This shows that any initial condition $f_0\in \M_1$ (cf.\ Proposition \ref{manifolds}) yields an entropy function $e_1(t)$ with the r.h.s.\ of \eqref{2d-decay} as its sharp envelope.\\

{}From the above discussion it is intuitively clear that the hypocoercive entropy method yields better decay rates. For general non-symmetric Fokker-Planck equations \eqref{nsFP} we have, in fact, the following comparison of the two decay estimates (standard entropy method vs.\ the new hypocoercive entropy method): 
\begin{Proposition}\label{decay-comp}
The decay rates $\lambda_K$ from \eqref{e-decay} and $\mu$ from Theorem \ref{entropydecay} satisfy:
 \begin{itemize}
  \item[(i)]  $\lambda_K\le\mu$ (and the strict inequality holds in many examples; see, e.g., Fig.\ \ref{fig:non-symmetric}).
  \item[(ii)]  If $C$ is diagonalizable (and hence also $\tilde C:={\sqrt{D}}^{\,-1}C\sqrt D$ with $\tilde C=\tilde A\Lambda \tilde A^{-1}$ and some diagonal matrix $\Lambda$), then 
$$
    \mu\le \kappa(\tilde A)^2\,\lambda_K\,,
$$
  where $\kappa(\tilde A):=\|\tilde A\|_2\|\tilde A^{-1}\|_2$ is the condition number of the matrix $\tilde A$.
\end{itemize}
\end{Proposition}
\begin{Proof}
(i): We first note that $\tilde K$ satisfies the continuous Lyapunov equation 
$$
  2\Id=\tilde C\tilde K + \tilde K\tilde C^T\,.
$$
Let $w$ be a (right) eigenvector of $\tilde C^T$ to an eigenvalue $\lambda_m$ with $\Re\{\lambda_m\}=\mu$. Hence
$$
  2|w|^2 = w^T\tilde C \tilde Kw+w^T\tilde K\tilde C^Tw = \bar\lambda_m w^T \tilde Kw+\lambda_m w^T\tilde Kw = 2\Re\{\lambda_m\}w^T\tilde K w\,.
$$
This yields the following estimate on the Rayleigh quotient of $\tilde K$:
$$
  \frac{1}{\mu} = \frac{w^T\tilde K w}{|w|^2}\le\lambda_{max}(\tilde K)
  =\frac{1}{\lambda_{min}(\tilde K^{-1})}=\frac{1}{\lambda_K}\,.\\
$$

(ii): For the upper bound on $\mu$ we consider again the Lyapunov equation for $\tilde K$.
Then Problem 9b of \S5.5, \cite{HoJoT91} gives the following bound on its solution:
$$
  \|\tilde K\|_2 \le \frac{\kappa(\tilde A)^2}{\mu}\,,
$$
and the result follows with $\|\tilde K\|_2=\lambda_{max}(\tilde K)=1/\lambda_K$. 
\end{Proof}\\

Finally, we remark that the hypocoercive entropy method cannot improve the standard decay estimate for symmetric Fokker-Planck equations: In that case, the matrix $D^{-1}C$ is symmetric positive definite in \eqref{nsFP}. Then, \eqref{Dequ} yields $K=C^{-1}D$. Hence, $\lambda_K:=\lambda_{\min}(\sqrt D K^{-1} \sqrt D)=\lambda_{\min}(\sqrt D^{\,-1} C \sqrt D)$, and $\mu=\lambda_K$ follows.


\section{Appendix}
\label{sec:appendix}

\begin{Proof}[ of Lemma \ref{PropSet} (propagation set of Fokker-Planck equations with linear drift)]\\

First, note that only drift-trajectories are non-constant in time, since the first row of $\tilde D$ is zero. A drift trajectory $\xi(s)=(t(s),v(s))$ starting at $\xi_0=(t_0,v_0)$ satisfies
\begin{align*}
 \dds \xi &= \left(\begin{array}{c} -1\\Cv\end{array}\right),\\
  \xi(0) &= \xi_0.
\end{align*}
The solution to this equation is
\begin{align*}
 \xi(s) &= \left(\begin{array}{c} t_0-s \\ e^{Cs}v_0 \end{array}\right).
\end{align*}
This means that drift trajectories move backwards in time linearly. Thus, for a point $q=(t',y)$ to be connected to $p=(t,x)$, it is necessary that $t'\leq t$. This is to be expected, as it is also the case for the classical maximum principle for parabolic equations.\\
Since the diffusion trajectories span the subspace $\R^k=\im D\subset\R^d$, we write $p=(t,x_D,x_0)$ and $q=(t',y_D,y_0)$, where $x_0$ and $y_0$ are the projection of $x$ and $y$ onto the kernel of $D$ (restricted to $\R^{d-k}$). Without moving backwards in time, we can only connect via diffusion trajectories. This implies
\begin{align*}
S(p,\R^{d+1})\cap\{(\tilde t,x)\in\R^{d+1}|\tilde t=t\}&=\{(t,x_0)\}\times\R^k.
\end{align*}
It remains to show that any point $q=(t',y)$ with $t'<t$ can be connected to $p$. The strategy here is the following: Since we can freely move around in $\im D$, we only need to connect $q$ and $p$ in the kernel of $D$ and in time. To achieve this, we employ Lemma \ref{Definiteness} (iv). We will proceed in a series of trajectories: A number of drift trajectories (equal to $\mu:=\dim\ker D+1=d-k+1$), each of them followed by up to $k=\operatorname{rank} D$ diffusion trajectories.
Starting at $\xi_0=(t,x)$, such a series of two drift and $2k$ diffusion trajectories will arrive at
\begin{align*}
 \big(t-s_1-s_2,e^{Cs_2}[e^{Cs_1}x+z_1]+z_2\big),
\end{align*}
where $z_1,z_2\in\im D$ are the results of shifts by diffusion trajectories and $0\leq s_1,s_2$. Thus, a series of $\mu$ trajectories will arrive at
\begin{align*}
 \big(t-\sum\limits_{j=1}^{\mu} s_j,\exp(C\sum\limits_{j=1}^{\mu}s_j)x+\sum\limits_{j=1}^{\mu-1} \exp(C\sum\limits_{l=1+j}^{\mu}s_l)z_j + z_{\mu}\big),
\end{align*}
where $z_j\in\im D$, $1\leq j\leq \mu$. Setting this equal to our target point $q=(t',y)$ and rearranging terms, we obtain the following requirements:
\begin{align}\label{system-s}
 \sum\limits_{j=1}^{\mu} s_j\stackrel{!}{=}t-t',\\
 y-e^{C(t-t')}x &\stackrel{!}{=} \sum\limits_{j=1}^{\mu-1} e^{Cr_j}z_j + z_\mu,\label{system-z}
\end{align}
with $r_j\in[0,t-t']$, $r_j:=\sum\limits_{l=j+1}^{\mu} s_l$, $s_j\geq 0$. The projection of  equation \eqref{system-z} (for $s_j,\,z_j$) onto $\im D$ can always be solved by choosing $z_\mu$ appropriately. For the projection onto $\ker D$, we get
\begin{align*}
 (\Id-D)\sum\limits_{j=1}^{\mu-1} e^{Cr_j}z_j &\stackrel{!}{=} (\Id-D)(y-e^{C(t-t')}x) =: v_0\in\ker D.
\end{align*}
The left hand side can be seen as a linear mapping from $(\im D)^{\mu-1}$ to $\ker D$, since each of the matrix exponentials can take an arbitrary argument $z_j\in\im D$. So we need to show that
\begin{align}
\label{bigmapping} ((\Id-D)e^{Cr_j})_{1\leq j\leq\mu} : (\im D)^{\mu-1} &\to \ker D,\\
\nonumber (z_j)_{1\leq j\leq \mu-1} &\mapsto (\Id-D)\sum\limits_{j=1}^{\mu-1} e^{Cr_j}z_j
\end{align}
is surjective for some choice of $0\leq r_{\mu-1}<r_{\mu-2}<\dots<r_1\leq t-t'$. Let $r_1\in[\fract{t-t'}2,t-t']$. Then either
\begin{align*}
 (\Id-D)e^{Cr_1}: \im D \to \ker D
\end{align*}
is surjective, or there is $\xi\in\ker D$ with $\xi\perp (\Id-D)e^{Cr_1}\im D$ (since the image of a linear map is always a linear subspace). But then, from Lemma \ref{Definiteness} (iv) there is $r_2\in(0,r_1)$ and $\eta\in\im D$ with
\begin{align*}
 \<(\Id-D)e^{Cr_2}\eta,\xi\>=\<\eta,e^{C^Tr_2}\xi\> &= 1.
\end{align*}
Now, since $\xi\not\perp (\Id-D)e^{Cr_2}\im D$, we have
\begin{align*}
 \dim\operatorname{span}\big[(\Id-D)e^{Cr_1}\im D,(\Id-D)e^{Cr_2}\im D\big] > \dim (\Id-D)e^{Cr_1}\im D.
\end{align*}
Then either
\begin{align*}
 \big((\Id-D)e^{Cr_1},(\Id-D)e^{Cr_2}\big): \im D\times\im D \to \ker D.
\end{align*}
is surjective, or we repeat the process. Each repetition increases by at least one the dimension of the reachable subspace of $\ker D$. Thus, we will need at most $\mu-1=\dim\ker D$ iterations, and hence the map (\ref{bigmapping}) is surjective.
\end{Proof}\\
\medskip


\begin{Proof}[ of the inclusion \eqref{prop-set} (propagation set of the kinetic Fokker-Planck equation)]\\

Here we shall consider points $q=(t_0,x_0,v_0)$ with $t_0<t_1$ that can be connected to the given point $p=(t_1,x_1,v_1)$ by a sequence of three trajectories (diffusion, drift, diffusion). Note that the form of $\tilde D$ in \eqref{D-tilde} lets the diffusion trajectories run purely in $v$--direction. Hence, we only need to connect $(t_0,x_0)$ to $(t_1,x_1)$ (both with arbitrary velocities) via a single drift trajectory.

For simplicity we set $s=-t$ in \eqref{p-ODE}. So, we consider the forward characteristic system
\begin{equation}\label{char-sys}
\begin{array}{rcl}
  \displaystyle\frac{d}{dt}x &=& v\,,\\[2mm]
  \displaystyle\frac{d}{dt}v &=& -V'(x)-\nu v\,,
\end{array}
\end{equation}
with the boundary data $x(t_0)=x_0$, $x(t_1)=x_1$. Here, $t_1$ and $x_1$ are fixed, and $V'$ is locally Lipschitz. Moreover, for all $x_0\in\R\setminus \{x_1\}$, we have to find a $\tilde t(x_0)\in[0,t_1)$ such that \eqref{char-sys} is solvable for all $t_0\in (\tilde t(x_0),t_1)$.

For the proof of Lemma \ref{PropSet} we just had to study the solvability of the linear system \eqref{system-s},  \eqref{system-z}. But for the nonlinear system \eqref{char-sys} this is not feasible explicitly. Hence we shall give an estimate on the propagation region of the characteristics.
To this end we define the total energy $H(x,v):=V(x)-V_{min}+\frac{v^2}{2}\ge0$ with $V_{min}:=\displaystyle\min_{x\in\R} V(x)$. Along a trajectory $\xi(t)$ of \eqref{char-sys}, it satisfies
$$
  -2\nu H(\xi(t)) \le \frac{d}{dt}H(\xi(t))=-\nu v^2 \le 0\,,
$$
and hence
\begin{equation}\label{H-inequ}
  H(\xi_1)\le H(\xi(t))\le H(\xi_1)\,e^{2\nu t_1}\,,\quad 0\le t\le t_1\,.
\end{equation}
Since $V\nearrow\infty$ as $|x|\to\infty$, all level curves $H(x,v)=const.$ are closed. We shall solve \eqref{char-sys} in the spirit of a (backward) shooting method starting at $t=t_1$ with some initial data $\xi_1=(x_1,\tilde v)$. Here, we shall choose $|\tilde v|$ large enough such that the trajectory passes ``above'' all local maxima of $V$ between $x_0$ and $x_1$. So, $\tilde v$ has to satisfy $H(x_1,\tilde v) > \displaystyle\max_{x_1\le x\le x_0} V(x)-V_{min}$ (w.r.o.g.\ we assumed here $x_1<x_0$).
For $|\tilde v|$ that large, the level curve $H=const:=H(x_1,\tilde v)$ (and hence also the trajectory $\xi(t),\,0\le t\le t_1$) crosses the line $x=x_0$ in the $(x,v)$ phase-plane (see Fig.\ \ref{characteristics}). 
The largest intersection time of $\xi(t)$ with $x=x_0$ can now be chosen as the desired time $\tilde t(x_0)<t_1$. 

\begin{figure}[htbp]
\hspace*{-6cm}
\includegraphics[width=25cm]{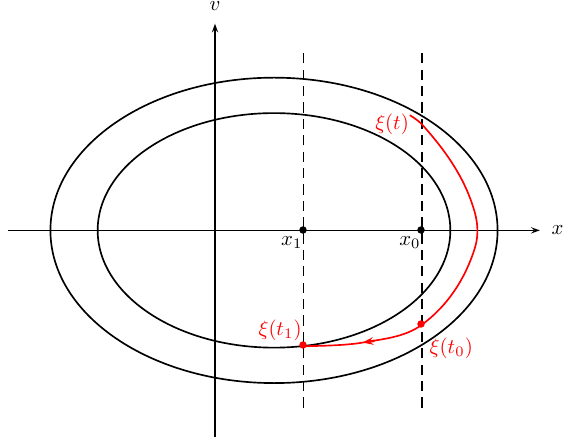} 
\vspace*{-24cm} 
\caption{\label{characteristics}The solution \textcolor{red}{$\xi(t)$} of the characteristic system \eqref{char-sys} lies in a domain bounded by two closed level curves, as estimated by \eqref{H-inequ}: $H=const:=H(x_1,\tilde v)$ (inner boundary), $H=const:=H(x_1,\tilde v)e^{2\nu t_1}$ (outer boundary).}
\end{figure}

Increasing $|\tilde v|$ further (and hence also $H(x_1,\tilde v)$) shows that the drift from $x_0$ to $x_1$ can be achieved in arbitrarily short time, i.e.\ for $t_0$ arbitrarily close to $t_1$. This proves the claim \eqref{prop-set}.
\end{Proof}\\
\medskip


\begin{Proof}[ of Lemma \ref{compactresolvent} (compactness of the resolvent of $L$)]\\

First we introduce the weighted $H^1$-space:
\begin{align*}
\H&:=\{f\in L^2|\nabla(\fract{f}{f_\infty})\in (L^2(\R^d,f_\infty))^d\},\\
\norm{f}_\H^2&:=\int\limits_{\R^d} |f|^2f_\infty^{-1}\d+\int\limits_{\R^d} |\nabla\fu|^2f_\infty\d=\norm{\fu}^2_{H^1(\R^d,f_\infty)}\,.
\end{align*}

For a (uniformly) elliptic operator, compactness of the resolvent can be shown by establishing that it maps $L^2$ into $\H$ and that the embedding $\H\hookrightarrow L^2$ is compact. However, for a degenerate elliptic operator, the resolvent will not map $L^2$ into $\H$, in general. So one has to work in spaces with fractional derivatives. For this proof, we shall therefore proceed in three steps. First we establish the space we work in, then we extend the regularisation result from Theorem \ref{entropyregularisation} for the solution semigroup $e^{Lt}$ 
on $L^2$. Finally, we use these two results to show compactness of the resolvent of $L$.\\
\\
\underline{Step 1 (interpolation spaces $\H_r$)}: We start by introducing the spaces $\H_r$, $0<r<1$, between $L^2$ and $\H$.
 An orthonormal basis $\{z_j|j\in\N_0^d\}$  of $L^2(\R^d,f_\infty)$ is given by the ``polynomial part'' of the eigenfunctions $z_j f_\infty$ of the (uniformly) elliptic Fokker-Planck operator
 \begin{align*}
  L_{\Id}f:=\div(\nabla(\fu)f_\infty)
 \end{align*}
 in $L^2$. They satisfy
 \begin{align*}
 L_{\Id}(z_jf_\infty) = -|j|z_jf_\infty,
 \end{align*}
 with $|j|$ the degree of the multi-index $j$.
 Then, for $f\in \H$ it holds
 \begin{align*}
  \norm{f}_{L^2}^2 = \sum\limits_{j\in\N_0^d} |c_j|^2,&\quad \norm{f}_{\H}^2 = \sum\limits_{j\in\N_0^d} (1+|j|)|c_j|^2,
 \end{align*}
 where $c_j$ is the coefficient of $\fu$ along $z_j$. We thus define
 \begin{align}
 \label{Halphadefinition} \H_r &:= \{f\in L^2|\sum\limits_{j\in\N_0^d} (1+|j|)^{r}|c_j|^2 <\infty\},
 \end{align}
 and have the interpolation inequality
 \begin{align}
  \label{Halphainterpolation} \norm{f}_{\H_r} &\leq \norm{f}_{\H}^{r}\norm{f}_{L^2}^{1-r}.
 \end{align}
 \underline{Step 2 (regularisation from $L^2$ to $\H_r$)}:
 Since $L$ generates a contraction semigroup on $L^2$, we have
 \begin{align}
  \label{csgproperty} \forall t\geq0:\quad \norm{e^{Lt}f}_{L^2} &\leq \norm{f}_{L^2}.
 \end{align}
 In the following estimate, we shall use the $L^2$-orthogonal decomposition $f=\tilde f + f_\infty\ild f\d$ with $\ild \tilde f\d =0$, and the scaled version of (\ref{regularityestimate}) for quadratic $\psi$:
 \begin{align*}
  \ild (\nabla\frac{f(t)}{f_\infty})^TP\nabla\frac{f(t)}{f_\infty}f_\infty \d &\leq ct^{-(2\tau+1)}\ild \Big(f-f_\infty\ild f\d\Big)^2 f_\infty^{-1} \d.
 \end{align*}
 We have:
 \begin{align*}
  \norm{e^{Lt}f}_{\H}^2 &= \norm{e^{Lt}f}_{L^2}^2 + \Big\|\nabla\frac{e^{Lt}f}{f_\infty}\Big\|_{L^2(\R^d,f_\infty)}^2 \\
	    &= \norm{e^{Lt}f}_{L^2}^2 + \Big\|\nabla\frac{e^{Lt}\tilde f}{f_\infty}\Big\|_{L^2(\R^d,f_\infty)}^2 \\
	    &\leq \norm{f}_{L^2}^2 + ct^{-(2\tau+1)}\norm{\tilde f}^2_{L^2} \\
	    &\leq (1+ct^{-(2\tau+1)})\norm{f}_{L^2}^2,
 \end{align*}
 where we have used the $L^2$-contractivity of $e^{Lt}$ and the positive definiteness of $P$. We thus obtain
  \begin{align}
 \label{sgregularisationcontrol} \forall 0< t\leq 1:\quad \norm{e^{Lt}f}_{\H} \leq \tilde c t^{-(\tau+\frac12)} \norm{f}_{L^2}
 \end{align}
 for all $f\in L^2$.
 By combining (\ref{Halphainterpolation}) -- (\ref{sgregularisationcontrol}), we obtain
 \begin{align}
  \label{Halphaestimate} \forall 0< t\leq 1:\quad \norm{e^{Lt}f}_{\H_r} &\leq \beta t^{-r(\tau+\frac12)}\norm{f}_{L^2}\,,
 \end{align}
 with $\beta:=\tilde c^r$.\\
  \underline{Step 3 (compact resolvent)}: For $r:=\frac1{\tau+1}>0$, we can integrate (\ref{Halphaestimate}) on $(0,1)$. This yields
 \begin{align}
  \label{compactresolvent1} \norm{\int\limits_0^1 e^{Lt}f\d[t]}_{\H_r} &\leq c \norm{f}_{L^2}.
 \end{align}
 By a well-known result for semigroups (see e.g. \cite{EnNaS06}, §II.1, Lemma 1.3 or \cite{PaS83}, §1.2, Theorem 2.4), for any $\lambda>0$ it holds that
 \begin{align*}
  \forall f\in D(L)\ \forall t>0:\quad \int\limits_0^t e^{(L-\lambda)s}(L-\lambda)f\d[s] &= e^{(L-\lambda)t}f - f.
 \end{align*}
 Due to (\ref{csgproperty}), $e^{(L-\lambda)t}$ decays exponentially and we conclude
 \begin{align}
 \label{restintcompact} \int\limits_1^\infty e^{(L-\lambda)t}(\lambda-L)f\d[t] &= e^{L-\lambda}f
 \end{align}
 for all $f\in D(L)$. Moreover (see e.g. \cite{EnNaS06}, §II.1, Theorem 1.10 or \cite{PaS83}, §1.3, Theorem 3.1), the resolvent $R(\lambda,L):= (\lambda-L)^{-1}$ has the representation
 \begin{align*}
  R(\lambda,L) &= \int\limits_0^\infty e^{(L-\lambda)t}\d[t] = \int\limits_0^1 e^{(L-\lambda)t}\d[t] + \int\limits_1^\infty e^{L-\lambda)t}\d[t].
 \end{align*}
 We apply this representation to (\ref{compactresolvent1}) and obtain
 \begin{align*}
  c\norm{f}_{L^2}\geq \norm{[R(\lambda,L)-\int\limits_1^\infty e^{(L-\lambda)t}]f\d[t]}_{\H_r},
 \end{align*}
 which yields
 \begin{align}
  \label{compactresolvent2} \norm{R(\lambda,L)f}_{\H_r} &\leq c\norm{f}_{L^2} + \norm{\int\limits_1^\infty e^{(L-\lambda)t}f\d[t]}_{\H_r}.
 \end{align}
 For $g\in D(L)$, we set in (\ref{compactresolvent2}) $f=(\lambda-L)g$ and obtain, using (\ref{restintcompact}),
 \begin{align*}
 \norm{g}_{\H_r}  &\leq c\norm{(\lambda-L)g}_{L^2} + \norm{\int\limits_1^\infty e^{(L-\lambda)t}(\lambda-L)g\d[t]}_{\H_r} \\
 &= c\norm{(\lambda-L)g}_{L^2} + e^{-\lambda}\norm{e^{L}g}_{\H_r}. 
 \end{align*}
 Applying (\ref{Halphaestimate}) with $t=1$ to the last term yields 
 \begin{align*}
  \norm{g}_{\H_r} &\leq c\norm{(\lambda-L)g}_{L^2}+\beta e^{-\lambda}\norm{g}_{L^2}.
 \end{align*}
Choosing $\lambda>\ln\beta$ allows to ``absorb'' the last term into the left-hand side, and hence
$$
  \|R(\lambda,L)f\|_{\H_r} \le c\|f\|_{L^2}\,.
$$
Due to the spectral representation of $\H_r$ in (\ref{Halphadefinition}), the embedding $\H_r\hookrightarrow L^2$ is compact for $r>0$. Hence, $R(\lambda,L)$ is compact for the chosen $\lambda$, and by the \emph{first resolvent formula} then also for all $\lambda$ in the resolvent set.
\end{Proof}\\
\medskip

\textit{Acknowledgement. The authors were supported by the FWF (project I 395-N16 and the doctoral school ``Dissipation and dispersion in non-linear partial differential equations'') and the ÖAD-project ``Long-time asymptotics for evolution equations in chemistry and biology.''}\\

\vspace{1cm}
\providecommand{\bysame}{\leavevmode\hbox to3em{\hrulefill}\thinspace}
\providecommand{\MR}{\relax\ifhmode\unskip\space\fi MR }
\providecommand{\MRhref}[2]{%
  \href{http://www.ams.org/mathscinet-getitem?mr=#1}{#2}
}
\providecommand{\href}[2]{#2}

\end{document}